\documentclass{IEEEtran}
\ifCLASSINFOpdf
\else
\fi
\usepackage[cmex10]{amsmath}
\usepackage{amsfonts}
\usepackage{xcolor}
\usepackage{rotating}
\usepackage{graphicx,multicol}

\DeclareMathOperator{\rank}{rank}

\DeclareMathOperator{\real}{Re}
\DeclareMathOperator{\imag}{Im}

\begin{document}
%
\title{Improving Optimal Power Flow Relaxations Using 3-Cycle Second-Order Cone Constraints

\author{Frederik~Geth,~\IEEEmembership{Member,~IEEE,}
        James Foster
\thanks{F. Geth and J. Foster are with the Energy Systems program, CSIRO Energy, Newcastle NSW, Australia (e-mail: frederik.geth@csiro.au, james.foster@csiro.au).}
}
}



\maketitle

%
%
%
%



\newcommand{\hermitiantranspose}{\text{H}}
\newcommand{\transpose}{\text{T}}
\newcommand{\complexconj}{*}

\newcommand{\Arrow}[1]{%
\parbox{#1}{\tikz{\draw[->](0,0)--(#1,0);}}
}
\newcommand{\Arrowright}[1]{%
\parbox{#1}{\tikz{\draw[<-](0,0)--(#1,0);}}
}

\newcommand{\indexScenario}{s}
\newcommand{\indexZone}{z}
\newcommand{\indexStakeholder}{b}
\newcommand{\indexUnitTech}{{x}_{\unitss}}
\newcommand{\indexPNTech}{{x}_{\pnss}}
\newcommand{\indexGITech}{{x}_{\giss}}
\newcommand{\indexGETech}{{x}_{\gess}}
\newcommand{\indexUnitRating}{{g}_{\unitss}}
\newcommand{\indexPNRating}{{g}_{\pnss}}
\newcommand{\indexGIRating}{{g}_{\giss}}
\newcommand{\indexGERating}{{g}_{\gess}}
\newcommand{\indexNodeRating}{{g}_{\nodess}}
\newcommand{\indexGridNode}{{i}}
\newcommand{\indexGridNodeTwo}{{j}}
\newcommand{\indexGridNodeThree}{{k}}
\newcommand{\indexGridLines}{{l}}
\newcommand{\indexGridLinesTwo}{{m}}
\newcommand{\indexGridLoad}{{h}}
\newcommand{\indexPhases}{p}
\newcommand{\indexPhasesTwo}{q}
\newcommand{\indexTimestep}{{k}}
\newcommand{\indexDay}{{d}}
\newcommand{\indexWeek}{{w}}
\newcommand{\indexUnit}{{u}}
\newcommand{\indexShunt}{{h}}
\newcommand{\indexPowerNode}{{e}}
\newcommand{\indexAdd}{{a}}
\newcommand{\indexRemove}{{r}}
\newcommand{\indexCirculation}{{\indexAdd\indexRemove}}
\newcommand{\indexFlows}{{f}}
\newcommand{\indexConductor}{n}
\newcommand{\indexdim}{n}
\newcommand{\indexdimtwo}{m}
\newcommand{\indexIteration}{k}
\newcommand{\indexContingency}{c}


\newcommand{\symPower}{P}
\newcommand{\symPowerToEnergy}{\gamma}
\newcommand{\symReactivePower}{Q}
\newcommand{\symApparentPower}{S}
\newcommand{\symVoltage}{U}
\newcommand{\symVoltageSOCP}{W}
\newcommand{\symVoltageAngle}{\theta}
\newcommand{\symPhaseDifference}{\varphi}
\newcommand{\symAnnualCost}{K}
\newcommand{\symPenalty}{Y}
\newcommand{\symPenaltyWeight}{N}
\newcommand{\symMass}{m}
\newcommand{\symSpecificMass}{\dot{\symMass}}
\newcommand{\symVolume}{V}
\newcommand{\symSpecificVolume}{v}
\newcommand{\symEnergy}{E}
\newcommand{\symEnergyFlow}{\dot{E}}
\newcommand{\symCurrent}{I}
\newcommand{\symCurrentSOCP}{L}
\newcommand{\symCurrentSquared}{L}
\newcommand{\symLineCurrent}{J}
\newcommand{\symObjectiveWeight}{w}
\newcommand{\symProbability}{\lambda}
\newcommand{\symBinary}{\alpha}
\newcommand{\symDecision}{{d}}
\newcommand{\symPhaseA}{{a}}
\newcommand{\symPhaseB}{{b}}
\newcommand{\symPhaseC}{{c}}
\newcommand{\symPhaseP}{{p}}
\newcommand{\symPhaseN}{{n}}
\newcommand{\symPhaseG}{{g}}
\newcommand{\symLossFactor}{{\rho}}
\newcommand{\symAvailability}{{b}}
\newcommand{\symDirectionality}{{d}}
\newcommand{\symPrice}{p}
\newcommand{\symInvestment}{C}
\newcommand{\symSpecificInvestment}{c}
\newcommand{\symLifetime}{\tau}
\newcommand{\symResistance}{r}
\newcommand{\symReactance}{x}
\newcommand{\symImpedance}{z}
\newcommand{\symAdmittance}{y}
\newcommand{\symSusceptance}{b}
\newcommand{\symConductance}{g}
\newcommand{\symLength}{l}
\newcommand{\symTime}{t}
\newcommand{\symRatio}{T}
\newcommand{\symLocation}{l}
\newcommand{\symEnthalpy}{H}
\newcommand{\symConnectivity}{C}

\newcommand{\symVolumeFlow}{\bar{\symVolume}}
\newcommand{\symMassFlow}{\bar{\symMass}}
\newcommand{\symDensity}{{\rho}}
\newcommand{\symGravitationalConstant}{\textcolor{\paramcolor}{g}}
\newcommand{\symElevation}{{h}}
\newcommand{\symElevationDifference}{\Delta{\symElevation}}
\newcommand{\symSetting}{\sigma}

\newcommand{\symRampRate}{\dot{\symPower}}
\newcommand{\symSpecRampRate}{{\zeta}}

\newcommand{\symPressure}[0]{p}
\newcommand{\symPressureSquared}[0]{\beta}
\newcommand{\symTemperature}[0]{   T }
\newcommand{\symHeatFlow}[0]{   Q }
\newcommand{\symViscosity}[0]{   \mu }
\newcommand{\symHeatResistance}[0]{   R }
\newcommand{\symFluidPipeDiameter}[0]{  \textcolor{\paramcolor}{D }    }
\newcommand{\symFluidPipeArea}[0]{  A     }
\newcommand{\symFluidPipeFriction}[0]{  \textcolor{\paramcolor}{f }    }
\newcommand{\symHeatCapacity}[0]{ c   }

\newcommand{\symRadiality}[0]{ \beta   }
\newcommand{\symSlack}{\epsilon}



\newcommand{\basess}{\text{{base}}}
\newcommand{\seriesss}{\text{{s}}}
\newcommand{\shuntss}{\text{{sh}}}
\newcommand{\invss}{\text{{inv}}}
\newcommand{\unbss}{\text{{unb}}}
\newcommand{\totss}{\text{{tot}}}
\newcommand{\batss}{\text{{bat}}}
\newcommand{\gridss}{\text{{grid}}}
\newcommand{\usabless}{\text{{usable}}}
\newcommand{\giss}{\text{{gi}}}
\newcommand{\pnss}{\text{{pn}}}
\newcommand{\lossss}{\text{{loss}}}
\newcommand{\cyclesss}{\text{{cycles}}}
\newcommand{\standbyss}{\text{{standby}}}
\newcommand{\operationss}{\text{{op}}}
\newcommand{\Pss}{\text{\sc{\symPower}}}
\newcommand{\Vss}{\text{\sc{\symVoltage}}}
\newcommand{\Iss}{\text{\sc{\symCurrent}}}
\newcommand{\Sss}{\text{\sc{\symApparentPower}}}
\newcommand{\Qss}{\text{\sc{\symReactivePower}}}
\newcommand{\Zss}{\text{\sc{\symImpedance}}}
\newcommand{\Ess}{\text{\sc{\symEnergy}}}
\newcommand{\nomss}{\text{nom}}
\newcommand{\ratedss}{\text{rated}}
\newcommand{\absss}{\text{abs}}
\newcommand{\maxss}{\text{max}}
\newcommand{\minss}{\text{min}}
\newcommand{\effss}{\text{eff}}
\newcommand{\chargess}{\text{c}}
\newcommand{\dischargess}{\text{d}}
\newcommand{\exoaddss}{\text{+}}
\newcommand{\exoremss}{\text{-}}
\newcommand{\orthss}{\perp}
\newcommand{\deprss}{\text{depr}}
\newcommand{\PFconvexss}{\text{PFconvex}}
\newcommand{\GIconvexss}{\text{GIconvex}}
\newcommand{\iterationss}{{(q)}}
\newcommand{\refss}{\text{ref}}
\newcommand{\addss}{\indexAdd,\text{add}}
\newcommand{\addtextss}{\text{add}}
\newcommand{\remss}{\indexRemove,\text{rem}}
\newcommand{\remtextss}{\text{rem}}
\newcommand{\circss}{\indexAdd\indexRemove,\text{circ}}
\newcommand{\flowss}{\indexFlows,\text{flow}}
\newcommand{\switchss}{\text{switch}}
\newcommand{\rotatedss}{\text{rot}}
\newcommand{\circumss}{\text{circum}}
\newcommand{\sbss}{\text{SB}}
\newcommand{\nonsbss}{\text{nonSB}}
\newcommand{\gess}{\text{ge}}
\newcommand{\unitss}{\text{unit}}
\newcommand{\nodess}{\text{node}}
\newcommand{\linss}{\text{lin}}
\newcommand{\realss}{\text{re}}
\newcommand{\imagss}{\text{im}}
\newcommand{\hagenposeuilless}[0]{  \text{ HP }}
\newcommand{\darcyweisbachss}[0]{ \text{  DC }}
\newcommand{\bentallinearization}[0]{ \text{  BenTal }}

\newcommand{\fortescuess}{\text{012}\!\cdot}

\newcommand{\plusss}{\text{+}}
\newcommand{\subsss}{\text{-}}
\newcommand{\combustionss}{\text{comb}}
\newcommand{\formationss}{\text{form}}



\newcommand{\PowerBase}{\symApparentPower_{\basess}}
\newcommand{\CurrentBase}{\symCurrent_{\basess}}
\newcommand{\VoltageBase}{\symVoltage_{\basess}^{\text{line}}}
\newcommand{\VoltageBasePhase}{\symVoltage_{\basess}^{\text{phase}}}
\newcommand{\TimeBase}{\symTime_{\basess}}
\newcommand{\EnergyBase}{\symEnergy_{\basess}}
\newcommand{\ImpedanceBase}{\symImpedance_{\basess}}
\newcommand{\AdmittanceBase}{\symAdmittance_{\basess}}
\newcommand{\CurrencyBase}{\euro_{\basess}}

\newcommand{\Knul}{\textcolor{\paramcolor}{\symAnnualCost_{0}}}
\newcommand{\Ynul}{\textcolor{\paramcolor}{\symPenalty_{0}}}



\newcommand{\yearunit}{a}
\newcommand{\VA}{VA}
\newcommand{\var}{var}
\newcommand{\kWh}{kWh}
\newcommand{\MWh}{MWh}
\newcommand{\si}[1]{#1}
\newcommand{\Si}[2]{\unit{#1}{#2}}
\newcommand{\percent}{\%}
\newcommand{\num}{}
\newcommand{\invtan}{\mathrm{atan2}}
\newcommand{\onen}[1]{\textcolor{\paramcolor}{\mathbf{1}_{#1}}}
\newcommand{\onenon}{\textcolor{\paramcolor}{\mathbf{1}}}

\newcommand{\nk}{n_{\setTimesteps}}
\newcommand{\nunits}{n_\setUnits}
\newcommand{\nn}{n_{\omegaGridNodes}}
\newcommand{\nc}{n_{\omegaGridLines}}
\newcommand{\nslackbuses}{n_{\omegaSlackbuses}}
\newcommand{\nnonslackbuses}{n_{\omegaNonSlackbuses}}

\newcommand{\setTimesteps}{\mathcal{K}}
\newcommand{\setTimestepsprev}{\mathcal{K}_0}
\newcommand{\omegaPNTech}{\mathcal{X}_{\pnss}}
\newcommand{\omegaGITech}{\mathcal{X}_{\giss}}
\newcommand{\omegaGETech}{\mathcal{X}_{\gess}}
\newcommand{\omegaUnitTech}{\mathcal{X}_{\unitss}}
\newcommand{\setPNRating}{\mathcal{G}_{\pnss}}
\newcommand{\setGIRating}{\mathcal{G}_{\giss}}
\newcommand{\setGERating}{\mathcal{G}_{\gess}}
\newcommand{\setUnitRating}{\mathcal{G}_{\unitss}}  
\newcommand{\setNodeRating}{\mathcal{G}_{\nodess}}  %
\newcommand{\omegaGridLines}{\textcolor{\paramcolor}{\mathcal{J}}}
\newcommand{\omegaGridNodes}{\textcolor{\paramcolor}{\mathcal{I}}}
\newcommand{\omegaGridNodesDomain}[1]{\omegaGridNodes_{#1}}
\newcommand{\omegaSlackbuses}{\omegaGridNodes_{\indexDomains, \sbss}}
\newcommand{\omegaNonSlackbuses}{\omegaGridNodes_{\indexDomains,\nonsbss}}
\newcommand{\setReferenceNodes}{\omegaGridNodes_{\text{ref}}}
\newcommand{\setReferenceNodesPrime}{\omegaGridNodes'_{\text{ref}}}
\newcommand{\setReferenceNodesRankOne}{\omegaGridNodes^{\text{R-1}}_{\text{ref}}}
\newcommand{\setReferenceNodesRankOnePrime}{\omegaGridNodes'^{\text{R-1}}_{\text{ref}}}
\newcommand{\setReferenceNodesBalanced}{\omegaGridNodes^{\text{bal}}_{\text{ref}}}
\newcommand{\setReferenceNodesFixed}{\omegaGridNodes^{\text{fix}}_{\text{ref}}}

\newcommand{\setPowerNodes}{\mathcal{E}}
\newcommand{\setAdd}{\mathcal{A}}
\newcommand{\setRemove}{\mathcal{R}}
\newcommand{\setCirculation}{\mathcal{C}}
\newcommand{\setFlows}{\mathcal{F}}
\newcommand{\setContingency}{\mathcal{C}}
\newcommand{\setTopology}{\textcolor{\paramcolor}{\mathcal{T}^{\Arrow{.1cm}}}}
\newcommand{\setTopologyBoth}{\textcolor{\paramcolor}{\mathcal{T}}}
\newcommand{\setTopologyReverse}{\textcolor{\paramcolor}{\mathcal{T}^{\Arrowright{.1cm}}}}
\newcommand{\setUnits}{\textcolor{\paramcolor}{\mathcal{U}}}
\newcommand{\setUnitTopology}{\textcolor{\paramcolor}{\mathcal{T^{\text{units}}}}}
\newcommand{\setShuntTopology}{\textcolor{\paramcolor}{\mathcal{T^{\text{shunts}}}}}

\newcommand{\setBinary}[0]{\textcolor{\setscolor}{ \left\{ 0,1 \right\} }}
\newcommand{\setComplex}[0]{\textcolor{\setscolor}{ \mathbb{C}   }}
\newcommand{\setReal}[0]{\textcolor{\setscolor}{ \mathbb{R}   }}
\newcommand{\setIntegers}[0]{\textcolor{\setscolor}{ \mathbb{Z}   }}
\newcommand{\setNaturalNumbers}[0]{\textcolor{\setscolor}{ \mathbb{N}   }}
\newcommand{\setNaturalNumbersPos}[0]{\textcolor{\setscolor}{ \mathbb{N}^{+}   }}
\newcommand{\setRealMatrix}[2]{\textcolor{\setscolor}{ \mathbb{R}^{#1 \times #2}   }}
\newcommand{\setComplexMatrix}[2]{\textcolor{\setscolor}{ \mathbb{C}^{#1\times #2}   }}
\newcommand{\setHermitianMatrix}[1]{\textcolor{\setscolor}{ \mathbb{H}^{#1}   }}
\newcommand{\setRealPos}[0]{\textcolor{\setscolor}{  \setReal_{+}    }}

\newcommand{\setPhases}{\textcolor{\paramcolor}{\mathcal{P}}}
\newcommand{\setPhasesTwo}{{\Phi}}
\newcommand{\setFortescue}{\textcolor{\paramcolor}{\mathcal{F}}}
\newcommand{\setCircuit}{\textcolor{\paramcolor}{\mathcal{C}}}
\newcommand{\setConductors}{\textcolor{\paramcolor}{\mathcal{N}}}
\newcommand{\setBuspairs}{\textcolor{\paramcolor}{\mathcal{B}}}
\newcommand{\setZIP}{\mathcal{Z}}

\newcommand{\setPSDConen}[1]{\textcolor{\setscolor}{  \mathbb{S}_{+}^{#1}    }}
\newcommand{\setNNOrthant}[0]{\textcolor{\setscolor}{  \setReal_{+}^{\indexdim}    }}
\newcommand{\setSOCone}[0]{\textcolor{\setscolor}{  \mathbb{Q}_{}^{\indexdim}    }}
\newcommand{\setPSDCone}[0]{\textcolor{\setscolor}{  \mathbb{S}_{+}^{\indexdim}    }}
\newcommand{\setPDCone}[0]{\textcolor{\setscolor}{  \mathbb{S}_{++}^{\indexdim}    }}
\newcommand{\setScenario}{\mathcal{S}}
\newcommand{\setZone}{\mathcal{Z}}
\newcommand{\setStakeholder}{\mathcal{B}}

\newcommand{\ineuro}[1]{\EUR{\ensuremath{#1}}}
\newcommand{\eurounit}{\euro}
\newcommand{\ineuroperyear}[1]{\ensuremath{#1}\,\eurounit \si{\per\yearunit}}

\newcommand{\variables}{\textcolor{black}{real-valued variables}}
\newcommand{\parameters}{\textcolor{\paramcolor}{real-valued parameters}}
\newcommand{\binaryparameters}{\textcolor{\decisioncolor}{binary parameters}}
\newcommand{\bounds}{\textcolor{\boundscolor}{generic bounds}}
\newcommand{\sizingvariables}{\textcolor{\sizingcolor}{sizing parameters}}
\newcommand{\sizingbounds}{\textcolor{\ratingboundscolor}{sizing bounds}}
\newcommand{\sizingparameters}{\textcolor{\investmentcolor}{sizing parameters}}
\newcommand{\binaryvariables}{\textcolor{\binarycolor}{binary variables}}
\newcommand{\timeparameters}{\textcolor{\timecolor}{time parameters}}
\newcommand{\priceparameters}{\textcolor{\pricecolor}{price parameters}}
\newcommand{\complexvariables}{\textcolor{\complexcolor}{complex variables}}
\newcommand{\complexparameters}{\textcolor{\complexparamcolor}{complex parameters}}

\newcommand{\paramcolor}{red}
\newcommand{\sizingcolor}{red}
\newcommand{\investmentcolor}{teal}
\newcommand{\decisioncolor}{brown}
\newcommand{\boundscolor}{violet}
\newcommand{\binarycolor}{darkgray}
\newcommand{\timecolor}{cyan}
\newcommand{\pricecolor}{brown}
\newcommand{\ratingboundscolor}{olive}
\newcommand{\setscolor}{darkgray}
\newcommand{\complexcolor}{blue}  
\newcommand{\complexparamcolor}{brown}  


\newcommand{\costheta}{{c}_{\indexGridNode \indexGridNodeTwo}   }
\newcommand{\sintheta}{{s}_{\indexGridNode \indexGridNodeTwo}   }
\newcommand{\usquared}{{u}_{\indexGridNode \indexGridNode}   }
\newcommand{\umult}{{u}_{\indexGridNode \indexGridNodeTwo}   }
\newcommand{\uucos}{{r}_{\indexGridNode \indexGridNodeTwo}   }
\newcommand{\uusin}{{t}_{\indexGridNode \indexGridNodeTwo}   }

\newcommand{\phasesi}{\textcolor{\paramcolor}{{\setPhasesTwo}_{\indexGridNode }   }}
\newcommand{\phasesj}{\textcolor{\paramcolor}{{\setPhasesTwo}_{\indexGridNodeTwo }   }}
\newcommand{\phasesij}{\textcolor{\paramcolor}{{\setPhasesTwo}_{\indexGridNode \indexGridNodeTwo}   }}

\newcommand{\fortescue}{\textcolor{\complexparamcolor}{\mathbf{F}}}
\newcommand{\alpharot}{\textcolor{\complexparamcolor}{\alpha}}

\newcommand{\Ts}{\textcolor{\timecolor}{T_{\text{s}} }}
\newcommand{\thor}{\textcolor{\timecolor}{t_{\text{h}}}}
\newcommand{\etac}{\textcolor{\paramcolor}{\eta_{\chargess} }}
\newcommand{\etad}{\textcolor{\paramcolor}{\eta_{\dischargess} }}
\newcommand{\etap}{\textcolor{\paramcolor}{\eta_{\exoaddss} }}
\newcommand{\etam}{\textcolor{\paramcolor}{\eta_{\exoremss} }}

\newcommand{\etaadd}{\textcolor{\paramcolor}{\eta_{\addss} }}
\newcommand{\etarem}{\textcolor{\paramcolor}{\eta_{\remss} }}


\newcommand{\Pck}[0]{\symPower_{\chargess} }
\newcommand{\Pdk}[0]{\symPower_{\dischargess} }
\newcommand{\Ppk}[0]{\symPower_{\exoaddss} }
\newcommand{\Pmk}[0]{\symPower_{\exoremss} }

\newcommand{\Paddk}[0]{\symPower_{\addss} }
\newcommand{\Premk}[0]{\symPower_{\remss} }
\newcommand{\Pflowk}[0]{\symPower_{\flowss} }
\newcommand{\Pflowkprev}[0]{\symPower_{\flowss} prev}
\newcommand{\Pflowrampmax}[0]{\textcolor{\paramcolor}{\symRampRate_{\flowss}^{\maxss}}}
\newcommand{\Pflowrampspecmax}[0]{\textcolor{\paramcolor}{\symSpecRampRate_{\flowss}}}

\newcommand{\Pclossk}[0]{\symPower_{\chargess, \lossss} }
\newcommand{\Pdlossk}[0]{\symPower_{\dischargess, \lossss} }
\newcommand{\Pplossk}[0]{\symPower_{\exoaddss, \lossss} }
\newcommand{\Pmlossk}[0]{\symPower_{\exoremss, \lossss} }
\newcommand{\Pgilossk}[0]{\symPower_{\exoremss, \lossss} }
\newcommand{\PPNlossk}[0]{\symPower_{\indexPowerNode, \lossss} }
\newcommand{\PPNexolossk}[0]{\symPower_{\indexPowerNode, \text{exo},\lossss} }
\newcommand{\PPNeleclossk}[0]{\symPower_{\indexPowerNode, \text{elec}, \lossss} }
\newcommand{\Punitlossk}[0]{\symPower_{\indexUnit, \lossss} }

\newcommand{\Paddlossk}[0]{\symPower_{\addss, \lossss} }
\newcommand{\Premlossk}[0]{\symPower_{\remss, \lossss} }

\newcommand{\Kswitch}[0]{\symAnnualCost^{\switchss}_{}}
\newcommand{\Kswitchl}[0]{\symAnnualCost^{\switchss}_{\indexGridLines}}

\newcommand{\Kgridloss}[0]{\symAnnualCost^{\gridss\lossss}_{}}
\newcommand{\Kgridlossl}[0]{\symAnnualCost^{\gridss\lossss}_{\indexGridLines}}

\newcommand{\Kc}[0]{\symAnnualCost^{\chargess}_{}}
\newcommand{\Kd}[0]{\symAnnualCost^{\dischargess}_{}}
\newcommand{\Kp}[0]{\symAnnualCost^{\exoaddss}_{}}
\newcommand{\Km}[0]{\symAnnualCost^{\exoremss}_{}}
\newcommand{\KE}[0]{\symAnnualCost^{\Ess}_{}}
\newcommand{\Kgi}[0]{\symAnnualCost^{\giss}_{}}
\newcommand{\KPN}[0]{\symAnnualCost_{\indexPowerNode}}
\newcommand{\Kunit}[0]{\symAnnualCost_{\indexUnit}}
\newcommand{\Kunits}[0]{\symAnnualCost_{\setUnits}}

\newcommand{\Kflow}[0]{\symAnnualCost^{\flowss}_{}}

\newcommand{\Kopex}[0]{\symAnnualCost^{\operationss}_{}}

\newcommand{\Kcoper}[0]{\symAnnualCost^{\chargess,\operationss}}
\newcommand{\Kdoper}[0]{\symAnnualCost^{\dischargess,\operationss}}
\newcommand{\Kctotoper}[0]{\symAnnualCost^{\chargess,\totss,\operationss}}
\newcommand{\Kdtotoper}[0]{\symAnnualCost^{\dischargess,\totss,\operationss}}
\newcommand{\Kpoper}[0]{\symAnnualCost^{\exoaddss,\operationss}}
\newcommand{\Kmoper}[0]{\symAnnualCost^{\exoremss,\operationss}}
\newcommand{\KPNoper}[0]{\symAnnualCost_{\indexPowerNode}^{\operationss}}
\newcommand{\Kunitoper}[0]{\symAnnualCost_{\indexUnit}^{\operationss}}
\newcommand{\Kunitsoper}[0]{\symAnnualCost_{\setUnits}^{\operationss}}

\newcommand{\Kflowoper}[0]{\symAnnualCost^{\flowss,\operationss}}

\newcommand{\Ec}[0]{\symEnergy_{\chargess}}
\newcommand{\Ed}[0]{\symEnergy_{\dischargess}}
\newcommand{\Ep}[0]{\symEnergy_{\exoaddss}}
\newcommand{\Em}[0]{\symEnergy_{\exoremss}}
\newcommand{\Egi}[0]{\symEnergy_{\giss}}
\newcommand{\EE}[0]{n_{\cyclesss}}

\newcommand{\Eflow}[0]{\symEnergy_{\flowss}}

\newcommand{\Ecflow}[0]{\symEnergyFlow_{\chargess}}
\newcommand{\Edflow}[0]{\symEnergyFlow_{\dischargess}}
\newcommand{\Epflow}[0]{\symEnergyFlow_{\exoaddss}}
\newcommand{\Emflow}[0]{\symEnergyFlow_{\exoremss}}
\newcommand{\Egiflow}[0]{\symEnergyFlow_{\giss}}
\newcommand{\EEflow}[0]{\dot{n}_{\cyclesss}}

\newcommand{\Eflowflow}[0]{\symEnergyFlow_{\flowss}}

\newcommand{\Ecmax}[0]{\textcolor{\paramcolor}{\symEnergy_{\chargess}^{\maxss}}}
\newcommand{\Edmax}[0]{\textcolor{\paramcolor}{\symEnergy_{\dischargess}^{\maxss}}}
\newcommand{\Epmax}[0]{\textcolor{\paramcolor}{\symEnergy_{\exoaddss}^{\maxss}}}
\newcommand{\Emmax}[0]{\textcolor{\paramcolor}{\symEnergy_{\exoremss}^{\maxss}}}
\newcommand{\Egimax}[0]{\textcolor{\paramcolor}{\symEnergy_{\giss}^{\maxss}}}
\newcommand{\EEmax}[0]{\textcolor{\paramcolor}{n_{\cyclesss}^{\maxss}}}

\newcommand{\Eflowmax}[0]{\textcolor{\paramcolor}{\symEnergy_{\flowss}^{\maxss}}}

\newcommand{\Ecmaxflow}[0]{\textcolor{\paramcolor}{\dot{\symEnergy}_{\chargess}^{\maxss}}}
\newcommand{\Edmaxflow}[0]{\textcolor{\paramcolor}{\dot{\symEnergy}_{\dischargess}^{\maxss}}}
\newcommand{\Epmaxflow}[0]{\textcolor{\paramcolor}{\dot{\symEnergy}_{\exoaddss}^{\maxss}}}
\newcommand{\Emmaxflow}[0]{\textcolor{\paramcolor}{\dot{\symEnergy}_{\exoremss}^{\maxss}}}
\newcommand{\Egimaxflow}[0]{\textcolor{\paramcolor}{\dot{\symEnergy}_{\giss}^{\maxss}}}
\newcommand{\EEmaxflow}[0]{\textcolor{\paramcolor}{\dot{n}_{\cyclesss}^{\maxss}}}

\newcommand{\Eflowmaxflow}[0]{\textcolor{\paramcolor}{\dot{\symEnergy}_{\flowss}^{\maxss}}}

\newcommand{\Kcdepr}[0]{\symAnnualCost^{\chargess,\deprss}}
\newcommand{\Kddepr}[0]{\symAnnualCost^{\dischargess,\deprss}}
\newcommand{\Kpdepr}[0]{\symAnnualCost^{\exoaddss,\deprss}}
\newcommand{\Kmdepr}[0]{\symAnnualCost^{\exoremss,\deprss}}
\newcommand{\KEdepr}[0]{\symAnnualCost^{\Ess,\deprss}}
\newcommand{\Kgidepr}[0]{\symAnnualCost^{\giss,\deprss}}
\newcommand{\KPNdepr}[0]{\symAnnualCost_{\indexPowerNode}^{\deprss}}
\newcommand{\Kunitdepr}[0]{\symAnnualCost_{\indexUnit}^{\deprss}}
\newcommand{\Kunitsdepr}[0]{\symAnnualCost_{\setUnits}^{\deprss}}

\newcommand{\Kflowdepr}[0]{\symAnnualCost^{\flowss,\deprss}}

\newcommand{\Cc}[0]{\textcolor{\investmentcolor}{\symInvestment_{\chargess}}}
\newcommand{\Cd}[0]{\textcolor{\investmentcolor}{\symInvestment_{\dischargess}}}
\newcommand{\Cp}[0]{\textcolor{\investmentcolor}{\symInvestment_{\exoaddss}}}
\newcommand{\Cm}[0]{\textcolor{\investmentcolor}{\symInvestment_{\exoremss}}}
\newcommand{\Cgi}[0]{\textcolor{\investmentcolor}{\symInvestment_{\Sss}}}
\newcommand{\CE}[0]{\textcolor{\investmentcolor}{\symInvestment_{\Ess}}}
\newcommand{\Cunit}[0]{\textcolor{\investmentcolor}{\symInvestment_{\indexUnit}}}
\newcommand{\CPN}[0]{\textcolor{\investmentcolor}{\symInvestment_{\indexPowerNode}}}
\newcommand{\Cunitmax}[0]{\textcolor{\boundscolor}{\symInvestment^{\maxss}_{\indexUnit}}}

\newcommand{\Cflow}[0]{\textcolor{\investmentcolor}{\symInvestment_{\flowss}}}

\newcommand{\cpc}[0]{\textcolor{\pricecolor}{\symSpecificInvestment_{\chargess}}}
\newcommand{\cpd}[0]{\textcolor{\pricecolor}{\symSpecificInvestment_{\dischargess}}}
\newcommand{\cpp}[0]{\textcolor{\pricecolor}{\symSpecificInvestment_{\exoaddss}}}
\newcommand{\cpm}[0]{\textcolor{\pricecolor}{\symSpecificInvestment_{\exoremss}}}
\newcommand{\cpgi}[0]{\textcolor{\pricecolor}{\symSpecificInvestment_{\Sss}}}
\newcommand{\cpE}[0]{\textcolor{\pricecolor}{\symSpecificInvestment_{\Ess}}}

\newcommand{\cpflow}[0]{\textcolor{\pricecolor}{\symSpecificInvestment_{\flowss}}}

\newcommand{\Vc}[0]{\textcolor{\investmentcolor}{\symVolume_{\chargess}}}
\newcommand{\Vd}[0]{\textcolor{\investmentcolor}{\symVolume_{\dischargess}}}
\newcommand{\Vp}[0]{\textcolor{\investmentcolor}{\symVolume_{\exoaddss}}}
\newcommand{\Vm}[0]{\textcolor{\investmentcolor}{\symVolume_{\exoremss}}}
\newcommand{\Vgi}[0]{\textcolor{\investmentcolor}{\symVolume_{\Sss}}}
\newcommand{\VE}[0]{\textcolor{\investmentcolor}{\symVolume_{\Ess}}}
\newcommand{\Vunit}[0]{\textcolor{\investmentcolor}{\symVolume_{\indexUnit}}}
\newcommand{\Vunitmax}[0]{\textcolor{\boundscolor}{\symVolume^{\maxss}_{\indexUnit}}}
\newcommand{\VPN}[0]{\textcolor{\investmentcolor}{\symVolume_{\indexPowerNode}}}

\newcommand{\Vflow}[0]{\textcolor{\investmentcolor}{\symVolume_{\flowss}}}

\newcommand{\vpc}[0]{\textcolor{\pricecolor}{\symSpecificVolume_{\chargess}}}
\newcommand{\vpd}[0]{\textcolor{\pricecolor}{\symSpecificVolume_{\dischargess}}}
\newcommand{\vpp}[0]{\textcolor{\pricecolor}{\symSpecificVolume_{\exoaddss}}}
\newcommand{\vpm}[0]{\textcolor{\pricecolor}{\symSpecificVolume_{\exoremss}}}
\newcommand{\vpgi}[0]{\textcolor{\pricecolor}{\symSpecificVolume_{\Sss}}}
\newcommand{\vpE}[0]{\textcolor{\pricecolor}{\symSpecificVolume_{\Ess}}}

\newcommand{\vpflow}[0]{\textcolor{\pricecolor}{\symSpecificVolume_{\flowss}}}

\newcommand{\Mc}[0]{\textcolor{\investmentcolor}{\symMass_{\chargess}}}
\newcommand{\Md}[0]{\textcolor{\investmentcolor}{\symMass_{\dischargess}}}
\newcommand{\Mp}[0]{\textcolor{\investmentcolor}{\symMass_{\exoaddss}}}
\newcommand{\Mm}[0]{\textcolor{\investmentcolor}{\symMass_{\exoremss}}}
\newcommand{\Mgi}[0]{\textcolor{\investmentcolor}{\symMass_{\Sss}}}
\newcommand{\ME}[0]{\textcolor{\investmentcolor}{\symMass_{\Ess}}}
\newcommand{\Munit}[0]{\textcolor{\investmentcolor}{\symMass_{\indexUnit}}}
\newcommand{\Munitmax}[0]{\textcolor{\boundscolor}{\symMass^{\maxss}_{\indexUnit}}}
\newcommand{\MPN}[0]{\textcolor{\investmentcolor}{\symMass_{\indexPowerNode}}}

\newcommand{\Mflow}[0]{\textcolor{\investmentcolor}{\symMass_{\flowss}}}

\newcommand{\mpc}[0]{\textcolor{\pricecolor}{\symSpecificMass_{\chargess}}}
\newcommand{\mpd}[0]{\textcolor{\pricecolor}{\symSpecificMass_{\dischargess}}}
\newcommand{\mpp}[0]{\textcolor{\pricecolor}{\symSpecificMass_{\exoaddss}}}
\newcommand{\mpm}[0]{\textcolor{\pricecolor}{\symSpecificMass_{\exoremss}}}
\newcommand{\mpgi}[0]{\textcolor{\pricecolor}{\symSpecificMass_{\Sss}}}
\newcommand{\mpE}[0]{\textcolor{\pricecolor}{\symSpecificMass_{\Ess}}}

\newcommand{\mpflow}[0]{\textcolor{\pricecolor}{\symSpecificMass_{\flowss}}}

\newcommand{\tauc}[0]{\textcolor{black}{\symLifetime_{\chargess}}}
\newcommand{\taud}[0]{\textcolor{black}{\symLifetime_{\dischargess}}}
\newcommand{\taup}[0]{\textcolor{black}{\symLifetime_{\exoaddss}}}
\newcommand{\taum}[0]{\textcolor{black}{\symLifetime_{\exoremss}}}
\newcommand{\taugi}[0]{\textcolor{black}{\symLifetime_{\Sss}}}
\newcommand{\tauE}[0]{\textcolor{black}{\symLifetime{_\Ess}}}

\newcommand{\tauflow}[0]{\textcolor{black}{\symLifetime_{\flowss}}}

\newcommand{\taucmax}[0]{\textcolor{\paramcolor}{\symLifetime^{\maxss}_{\chargess}}}
\newcommand{\taudmax}[0]{\textcolor{\paramcolor}{\symLifetime{^{\maxss}_\dischargess}}}
\newcommand{\taupmax}[0]{\textcolor{\paramcolor}{\symLifetime^{\maxss}_{\exoaddss}}}
\newcommand{\taummax}[0]{\textcolor{\paramcolor}{\symLifetime^{\maxss}_{\exoremss}}}
\newcommand{\taugimax}[0]{\textcolor{\paramcolor}{\symLifetime^{\maxss}_{\Sss}}}
\newcommand{\tauEmax}[0]{\textcolor{\paramcolor}{\symLifetime^{\maxss}_{\Ess}}}

\newcommand{\tauflowmax}[0]{\textcolor{\paramcolor}{\symLifetime^{\maxss}_{\flowss}}}

\newcommand{\ppc}[0]{\textcolor{\pricecolor}{\symPrice_{\chargess} }}
\newcommand{\ppd}[0]{\textcolor{\pricecolor}{\symPrice_{\dischargess} }}
\newcommand{\ppp}[0]{\textcolor{\pricecolor}{\symPrice_{\exoaddss} }}
\newcommand{\ppm}[0]{\textcolor{\pricecolor}{\symPrice_{\exoremss} }}

\newcommand{\ppflow}[0]{\textcolor{\pricecolor}{\symPrice_{\flowss} }}

\newcommand{\ppswitch}[0]{\textcolor{\pricecolor}{\symPrice_{\switchss} }}
\newcommand{\ppswitchl}[0]{\textcolor{\pricecolor}{\symPrice_{\indexGridLines,\switchss} }}

\newcommand{\ppPgridlossl}[0]{\textcolor{\pricecolor}{\symPrice_{\indexGridLines\gridss\lossss}^{\Pss}    }}
\newcommand{\ppQgridlossl}[0]{\textcolor{\pricecolor}{\symPrice_{\indexGridLines\gridss\lossss}^{\Qss} }}
\newcommand{\ppQabsgridlossl}[0]{\textcolor{\pricecolor}{\symPrice_{\indexGridLines\gridss\lossss}^{|\Qss|} }}

\newcommand{\Pckref}[0]{\textcolor{\paramcolor}{\symPower^{\refss}_{\chargess} }}
\newcommand{\Pdkref}[0]{\textcolor{\paramcolor}{\symPower^{\refss}_{\dischargess} }}
\newcommand{\Ppkref}[0]{\textcolor{\paramcolor}{\symPower^{\refss}_{\exoaddss} }}
\newcommand{\Pmkref}[0]{\textcolor{\paramcolor}{\symPower^{\refss}_{\exoremss} }}

\newcommand{\Pckmax}[0]{\textcolor{\boundscolor}{\symPower^{\maxss}_{\chargess} }}
\newcommand{\Pdkmax}[0]{\textcolor{\boundscolor}{\symPower^{\maxss}_{\dischargess} }}
\newcommand{\Ppkmax}[0]{\textcolor{\boundscolor}{\symPower^{\maxss}_{\exoaddss} }}
\newcommand{\Pmkmax}[0]{\textcolor{\boundscolor}{\symPower^{\maxss}_{\exoremss} }}

\newcommand{\Pckmin}[0]{\textcolor{\boundscolor}{\symPower^{\minss}_{\chargess} }}
\newcommand{\Pdkmin}[0]{\textcolor{\boundscolor}{\symPower^{\minss}_{\dischargess} }}
\newcommand{\Ppkmin}[0]{\textcolor{\boundscolor}{\symPower^{\minss}_{\exoaddss} }}
\newcommand{\Pmkmin}[0]{\textcolor{\boundscolor}{\symPower^{\minss}_{\exoremss} }}
\newcommand{\Pflowkref}[0]{\textcolor{\paramcolor}{\symPower^{\refss}_{\flowss} }}

\newcommand{\Paddkmin}[0]{\textcolor{\boundscolor}{\symPower^{\minss}_{\addss} }}
\newcommand{\Paddkmax}[0]{\textcolor{\boundscolor}{\symPower^{\maxss}_{\addss} }}
\newcommand{\Premkmin}[0]{\textcolor{\boundscolor}{\symPower^{\minss}_{\remss} }}
\newcommand{\Premkmax}[0]{\textcolor{\boundscolor}{\symPower^{\maxss}_{\remss} }}
\newcommand{\Pflowkmin}[0]{\textcolor{\boundscolor}{\symPower^{\minss}_{\flowss} }}
\newcommand{\Pflowkmax}[0]{\textcolor{\boundscolor}{\symPower^{\maxss}_{\flowss} }}

\newcommand{\Pcratedmax}[0]{\textcolor{\ratingboundscolor}{\symPower'^{\maxss}_{\chargess}}}
\newcommand{\Pdratedmax}[0]{\textcolor{\ratingboundscolor}{\symPower'^{\maxss}_{\dischargess}}}
\newcommand{\Ppratedmax}[0]{\textcolor{\ratingboundscolor}{\symPower'^{\maxss}_{\exoaddss}}}
\newcommand{\Pmratedmax}[0]{\textcolor{\ratingboundscolor}{\symPower'^{\maxss}_{\exoremss}}}

\newcommand{\Pcratedmin}[0]{\textcolor{\ratingboundscolor}{\symPower'^{\minss}_{\chargess}}}
\newcommand{\Pdratedmin}[0]{\textcolor{\ratingboundscolor}{\symPower'^{\minss}_{\dischargess}}}
\newcommand{\Ppratedmin}[0]{\textcolor{\ratingboundscolor}{\symPower'^{\minss}_{\exoaddss}}}
\newcommand{\Pmratedmin}[0]{\textcolor{\ratingboundscolor}{\symPower'^{\minss}_{\exoremss}}}

\newcommand{\Paddratedmin}[0]{\textcolor{\ratingboundscolor}{\symPower'^{\minss}_{\addss}}}
\newcommand{\Paddratedmax}[0]{\textcolor{\ratingboundscolor}{\symPower'^{\maxss}_{\addss}}}
\newcommand{\Premratedmin}[0]{\textcolor{\ratingboundscolor}{\symPower'^{\minss}_{\remss}}}
\newcommand{\Premratedmax}[0]{\textcolor{\ratingboundscolor}{\symPower'^{\maxss}_{\remss}}}

\newcommand{\Pflowrated}[0]{\textcolor{\sizingcolor}{\symPower'^{}_{\flowss}}}
\newcommand{\Pflowratedmin}[0]{\textcolor{\ratingboundscolor}{\symPower'^{\minss}_{\flowss}}}
\newcommand{\Pflowratedmax}[0]{\textcolor{\ratingboundscolor}{\symPower'^{\maxss}_{\flowss}}}

\newcommand{\Eratedmaxmax}[0]{\textcolor{\ratingboundscolor}{\symEnergy'^{\maxss\maxss}}}
\newcommand{\Eratedmaxmin}[0]{\textcolor{\ratingboundscolor}{\symEnergy'^{\maxss\minss}}}
\newcommand{\Eratedmax}[0]{\textcolor{\sizingcolor}{\symEnergy'^{\maxss}}}
\newcommand{\Eratedminmax}[0]{\textcolor{\ratingboundscolor}{\symEnergy'^{\minss\maxss}}}
\newcommand{\Eratedminmin}[0]{\textcolor{\ratingboundscolor}{\symEnergy'^{\minss\minss}}}
\newcommand{\Eratedmin}[0]{\textcolor{\sizingcolor}{\symEnergy'^{\minss}}}

\newcommand{\Pcrated}[0]{\textcolor{\sizingcolor}{\symPower'^{}_{\chargess}}}
\newcommand{\Pdrated}[0]{\textcolor{\sizingcolor}{\symPower'^{}_{\dischargess}}}
\newcommand{\Pprated}[0]{\textcolor{\sizingcolor}{\symPower'^{}_{\exoaddss}}}
\newcommand{\Pmrated}[0]{\textcolor{\sizingcolor}{\symPower'^{}_{\exoremss}}}

\newcommand{\Paddrated}[0]{\textcolor{\sizingcolor}{\symPower'^{}_{\addss}}}
\newcommand{\Premrated}[0]{\textcolor{\sizingcolor}{\symPower'^{}_{\remss}}}

\newcommand{\Pcorthk}[0]{\symPower^{\orthss}_{\chargess} }
\newcommand{\Pdorthk}[0]{\symPower^{\orthss}_{\dischargess} }
\newcommand{\Pporthk}[0]{\symPower^{\orthss}_{\exoaddss} }
\newcommand{\Pmorthk}[0]{\symPower^{\orthss}_{\exoremss} }
\newcommand{\Porthk}[0]{\symPower^{\orthss}_{} }

\newcommand{\Paddorthk}[0]{\symPower^{\orthss}_{\addss} }
\newcommand{\Premorthk}[0]{\symPower^{\orthss}_{\remss} }
\newcommand{\Pfloworthk}[0]{\symPower^{\orthss}_{\flowss} }
\newcommand{\Paddorthtotk}[0]{\symPower^{\orthss}_{\addtextss} }
\newcommand{\Premorthtotk}[0]{\symPower^{\orthss}_{\remtextss} }

\newcommand{\Pcorthkmax}[0]{\textcolor{\boundscolor}{\symPower^{\orthss, \maxss}_{\chargess} }}
\newcommand{\Pdorthkmax}[0]{\textcolor{\boundscolor}{\symPower^{\orthss, \maxss}_{\dischargess} }}
\newcommand{\Pporthkmax}[0]{\textcolor{\boundscolor}{\symPower^{\orthss, \maxss}_{\exoaddss} }}
\newcommand{\Pmorthkmax}[0]{\textcolor{\boundscolor}{\symPower^{\orthss, \maxss}_{\exoremss} }}

\newcommand{\Paddorthkmax}[0]{\textcolor{\boundscolor}{\symPower^{\orthss, \maxss}_{\addss} }}
\newcommand{\Premorthkmax}[0]{\textcolor{\boundscolor}{\symPower^{\orthss, \maxss}_{\remss} }}
\newcommand{\Pfloworthkmax}[0]{\textcolor{\boundscolor}{\symPower^{\orthss, \maxss}_{\flowss} }}

\newcommand{\Pcdk}[0]{\symPower_{\chargess\dischargess} }
\newcommand{\Pcmk}[0]{\symPower_{\chargess\exoremss} }
\newcommand{\Ppmk}[0]{\symPower_{\exoaddss\exoremss} }
\newcommand{\Ppdk}[0]{\symPower_{\exoaddss\dischargess} }

\newcommand{\Pcirck}[0]{\symPower_{\circss} }
\newcommand{\Pcircmax}[0]{\symPower^{\maxss}_{\circss}}

\newcommand{\Pcdmax}[0]{\symPower^{\maxss}_{\chargess\dischargess}}
\newcommand{\Pcmmax}[0]{\symPower^{\maxss}_{\chargess\exoremss}}
\newcommand{\Ppmmax}[0]{\symPower^{\maxss}_{\exoaddss\exoremss}}
\newcommand{\Ppdmax}[0]{\symPower^{\maxss}_{\exoaddss\dischargess}}

\newcommand{\Pcpk}[0]{\symPower_{\chargess\exoaddss} }
\newcommand{\Pdmk}[0]{\symPower_{\dischargess\exoremss} }
\newcommand{\gammacp}[0]{\textcolor{\paramcolor}{\symPowerToEnergy_{\chargess\exoaddss} }}
\newcommand{\gammadm}[0]{\textcolor{\paramcolor}{\symPowerToEnergy_{\dischargess\exoremss} }}

\newcommand{\gammaadd}[0]{\textcolor{\paramcolor}{\symPowerToEnergy_{\addss} }}
\newcommand{\gammarem}[0]{\textcolor{\paramcolor}{\symPowerToEnergy_{\remss} }}

\newcommand{\dcd}[0]{\textcolor{\decisioncolor}{\symDecision_{\chargess\dischargess} }}
\newcommand{\dcm}[0]{\textcolor{\decisioncolor}{\symDecision_{\chargess\exoremss} }}
\newcommand{\dpm}[0]{\textcolor{\decisioncolor}{\symDecision_{\exoaddss\exoremss} }}
\newcommand{\dpd}[0]{\textcolor{\decisioncolor}{\symDecision_{\exoaddss\dischargess} }}

\newcommand{\dcirc}[0]{\textcolor{\decisioncolor}{\symDecision_{\circss} }}

\newcommand{\dcorth}[0]{\textcolor{\decisioncolor}{\symDecision^{\orthss}_{\chargess} }}
\newcommand{\ddorth}[0]{\textcolor{\decisioncolor}{\symDecision^{\orthss}_{\dischargess} }}
\newcommand{\dporth}[0]{\textcolor{\decisioncolor}{\symDecision^{\orthss}_{\exoaddss} }}
\newcommand{\dmorth}[0]{\textcolor{\decisioncolor}{\symDecision^{\orthss}_{\exoremss} }}

\newcommand{\daddorth}[0]{\textcolor{\decisioncolor}{\symDecision^{\orthss}_{\addss} }}
\newcommand{\dremorth}[0]{\textcolor{\decisioncolor}{\symDecision^{\orthss}_{\remss} }}
\newcommand{\dfloworth}[0]{\textcolor{\decisioncolor}{\symDecision^{\orthss}_{\flowss} }}

\newcommand{\dclevels}[0]{\textcolor{\decisioncolor}{\symDecision^{}_{\chargess}  }}
\newcommand{\ddlevels}[0]{\textcolor{\decisioncolor}{\symDecision^{}_{\dischargess}  }}
\newcommand{\dplevels}[0]{\textcolor{\decisioncolor}{\symDecision^{}_{\exoaddss}  }}
\newcommand{\dmlevels}[0]{\textcolor{\decisioncolor}{\symDecision^{}_{\exoremss}  }}
\newcommand{\dflowlevels}[0]{\textcolor{\decisioncolor}{\symDecision^{}_{\flowss}  }}

\newcommand{\dcmaxlevels}[0]{\textcolor{\paramcolor}{\symDecision^{\maxss}_{\chargess}  }}
\newcommand{\ddmaxlevels}[0]{\textcolor{\paramcolor}{\symDecision^{\maxss}_{\dischargess}  }}
\newcommand{\dpmaxlevels}[0]{\textcolor{\paramcolor}{\symDecision^{\maxss}_{\exoaddss}  }}
\newcommand{\dmmaxlevels}[0]{\textcolor{\paramcolor}{\symDecision^{\maxss}_{\exoremss}  }}
\newcommand{\dflowmaxlevels}[0]{\textcolor{\paramcolor}{\symDecision^{\maxss}_{\flowss}  }}


\newcommand{\dE}[0]{\textcolor{\decisioncolor}{\symDecision_{\Ess}}}

\newcommand{\Emax}[0]{\textcolor{\sizingcolor}{\symEnergy'^{\maxss}}}
\newcommand{\Emin}[0]{\textcolor{\sizingcolor}{\symEnergy'^{\minss}}}
\newcommand{\Eeff}[0]{\textcolor{\sizingcolor}{\symEnergy'^{\usabless}}}
\newcommand{\Ekmax}[0]{\textcolor{\boundscolor}{\symEnergy^{\maxss} }}
\newcommand{\Ekmin}[0]{\textcolor{\boundscolor}{\symEnergy^{\minss} }}
\newcommand{\Ekeff}[0]{\textcolor{\boundscolor}{\symEnergy^{\usabless} }}
\newcommand{\Ek}[0]{\symEnergy }
\newcommand{\Ekprev}[0]{\symEnergy prev}
\newcommand{\iEk}[0]{\textcolor{\binarycolor}{\symBinary_{\Ess} }}
\newcommand{\iEkprev}[0]{\textcolor{\binarycolor}{\symBinary_{\Ess} prev}}

\newcommand{\iijk}[0]{\textcolor{\binarycolor}{\symBinary_{\indexGridNode \indexGridNodeTwo} }}
\newcommand{\ijik}[0]{\textcolor{\binarycolor}{\symBinary_{\indexGridNodeTwo \indexGridNode} }}
\newcommand{\ilk}[0]{\textcolor{\binarycolor}{\symBinary_{\indexGridLines} }}
\newcommand{\ilkprev}[0]{\textcolor{\binarycolor}{\symBinary_{\indexGridLines} prev}}
\newcommand{\ilkmin}[0]{\textcolor{\boundscolor}{\symBinary^{\minss}_{\indexGridLines} }}
\newcommand{\ilkmax}[0]{\textcolor{\boundscolor}{\symBinary^{\maxss}_{\indexGridLines} }}

\newcommand{\betaijk}[0]{\textcolor{\binarycolor}{\symRadiality_{\indexGridNode \indexGridNodeTwo} }}
\newcommand{\betajik}[0]{\textcolor{\binarycolor}{\symRadiality_{\indexGridNodeTwo \indexGridNode} }}


\newcommand{\Sgik}[0]{\textcolor{\complexcolor}{\symApparentPower_{\giss} }}
\newcommand{\Sauxk}[0]{\symApparentPower_{\text{aux}} }
\newcommand{\SauxkA}[0]{\symApparentPower_{\text{aux},\symPhaseA} }
\newcommand{\SauxkB}[0]{\symApparentPower_{\text{aux},\symPhaseB} }
\newcommand{\SauxkC}[0]{\symApparentPower_{\text{aux},\symPhaseC} }
\newcommand{\Sgirated}[0]{\textcolor{\sizingcolor}{\symApparentPower_{\giss}^{'}}}
\newcommand{\lgik}[0]{\textcolor{\paramcolor}{\symLocation_{\giss} }}
\newcommand{\Pgik}[0]{\symPower_{\giss} }
\newcommand{\Qgik}[0]{\symReactivePower_{\giss} }
\newcommand{\Qgikref}[0]{\textcolor{\boundscolor}{\symReactivePower^{\refss}_{\giss} }}
\newcommand{\Qgikmax}[0]{\textcolor{\boundscolor}{\symReactivePower^{\maxss}_{\giss} }}
\newcommand{\Qgikmin}[0]{\textcolor{\boundscolor}{\symReactivePower^{\minss}_{\giss} }}
\newcommand{\Pctotk}[0]{\symPower_{\chargess,\totss} }
\newcommand{\Pdtotk}[0]{\symPower_{\dischargess,\totss} }
\newcommand{\Pstandbyk}[0]{\symPower_{\standbyss} }
\newcommand{\PSlossk}[0]{\symPower_{\Sss,\lossss} }
\newcommand{\Punblossk}[0]{\symPower_{\unbss,\lossss} }
\newcommand{\etagi}{\textcolor{\paramcolor}{\eta_\giss  }}

\newcommand{\Pstandbykmin}[0]{\textcolor{\boundscolor}{\symPower_{\standbyss}^{\minss}  }}
\newcommand{\Pstandbykmax}[0]{\textcolor{\boundscolor}{\symPower_{\standbyss}^{\maxss}  }}

\newcommand{\aPZk}[0]{\textcolor{\paramcolor}{a_{  }^{ \Zss\Pss }  }}
\newcommand{\aPIk}[0]{\textcolor{\paramcolor}{a_{  }^{ \Iss\Pss }  }}
\newcommand{\aPPk}[0]{\textcolor{\paramcolor}{a_{  }^{ \Pss\Pss }  }}
\newcommand{\aQZk}[0]{\textcolor{\paramcolor}{a_{  }^{ \Zss\Qss }  }}
\newcommand{\aQIk}[0]{\textcolor{\paramcolor}{a_{  }^{ \Iss\Qss }  }}
\newcommand{\aQPk}[0]{\textcolor{\paramcolor}{a_{  }^{ \Pss\Qss }  }}

\newcommand{\aPZ}[0]{\textcolor{\paramcolor}{a_{ \indexZIP }^{ \Zss \Pss} }}
\newcommand{\aPI}[0]{\textcolor{\paramcolor}{a_{ \indexZIP }^{ \Iss \Pss} }}
\newcommand{\aPP}[0]{\textcolor{\paramcolor}{a_{\indexZIP  }^{ \Pss\Pss } }}
\newcommand{\aQZ}[0]{\textcolor{\paramcolor}{a_{  \indexZIP}^{ \Zss\Qss } }}
\newcommand{\aQI}[0]{\textcolor{\paramcolor}{a_{ \indexZIP }^{ \Iss \Qss} }}
\newcommand{\aQP}[0]{\textcolor{\paramcolor}{a_{ \indexZIP }^{ \Pss\Qss } }}

\newcommand{\PgikA}[0]{\symPower_{\giss,\symPhaseA} }
\newcommand{\QgikA}[0]{\symReactivePower_{\giss,\symPhaseA} }
\newcommand{\SgikA}[0]{\textcolor{\complexcolor}{\symApparentPower_{\giss,\symPhaseA} }}
\newcommand{\PgikB}[0]{\symPower_{\giss,\symPhaseB} }
\newcommand{\QgikB}[0]{\symReactivePower_{\giss,\symPhaseB} }
\newcommand{\SgikB}[0]{\textcolor{\complexcolor}{\symApparentPower_{\giss,\symPhaseB} }}
\newcommand{\PgikC}[0]{\symPower_{\giss,\symPhaseC} }
\newcommand{\QgikC}[0]{\symReactivePower_{\giss,\symPhaseC} }
\newcommand{\SgikC}[0]{\textcolor{\complexcolor}{\symApparentPower_{\giss,\symPhaseC} }}
\newcommand{\SgiratedA}[0]{\textcolor{\sizingcolor}{\symApparentPower_{\giss,\symPhaseA}^{'}}}
\newcommand{\SgiratedB}[0]{\textcolor{\sizingcolor}{\symApparentPower_{\giss,\symPhaseB}^{'}}}
\newcommand{\SgiratedC}[0]{\textcolor{\sizingcolor}{\symApparentPower_{\giss,\symPhaseC}^{'}}}

\newcommand{\Sgiratedmin}[0]{\textcolor{\ratingboundscolor}{\symApparentPower_{\giss}^{\minss}}}
\newcommand{\SgiratedminA}[0]{\textcolor{\ratingboundscolor}{\symApparentPower_{\giss,\symPhaseA}^{\minss}}}
\newcommand{\SgiratedminB}[0]{\textcolor{\ratingboundscolor}{\symApparentPower_{\giss,\symPhaseB}^{\minss}}}
\newcommand{\SgiratedminC}[0]{\textcolor{\ratingboundscolor}{\symApparentPower_{\giss,\symPhaseC}^{\minss}}}
\newcommand{\Sgiratedmax}[0]{\textcolor{\ratingboundscolor}{\symApparentPower_{\giss}^{\maxss}}}
\newcommand{\SgiratedmaxA}[0]{\textcolor{\ratingboundscolor}{\symApparentPower_{\giss,\symPhaseA}^{\maxss}}}
\newcommand{\SgiratedmaxB}[0]{\textcolor{\ratingboundscolor}{\symApparentPower_{\giss,\symPhaseB}^{\maxss}}}
\newcommand{\SgiratedmaxC}[0]{\textcolor{\ratingboundscolor}{\symApparentPower_{\giss,\symPhaseC}^{\maxss}}}

\newcommand{\Pgikunb}[0]{\symPower_{\giss,\unbss} }
\newcommand{\Qgikunb}[0]{\symReactivePower_{\giss,\unbss} }

\newcommand{\rhoPgikunb}[0]{\textcolor{\paramcolor}{\symLossFactor_{\giss,\Pss, \unbss} }}
\newcommand{\rhoQgikunb}[0]{\textcolor{\paramcolor}{\symLossFactor_{\giss,\Qss,\unbss} }}

\newcommand{\bgik}[0]{\textcolor{\paramcolor}{\symAvailability_{\giss} }}
\newcommand{\bgikA}[0]{\textcolor{\paramcolor}{\symAvailability_{\giss, \symPhaseA} }}
\newcommand{\bgikB}[0]{\textcolor{\paramcolor}{\symAvailability_{\giss, \symPhaseB} }}
\newcommand{\bgikC}[0]{\textcolor{\paramcolor}{\symAvailability_{\giss, \symPhaseC} }}

\newcommand{\dPtogik}[0]{\textcolor{\paramcolor}{\symDirectionality_{\giss}^{\Pss\plusss}   }}
\newcommand{\dPtogikA}[0]{\textcolor{\paramcolor}{\symDirectionality_{\giss, \symPhaseA}^{\Pss\plusss} }}
\newcommand{\dPtogikB}[0]{\textcolor{\paramcolor}{\symDirectionality_{\giss, \symPhaseB}^{\Pss\plusss} }}
\newcommand{\dPtogikC}[0]{\textcolor{\paramcolor}{\symDirectionality_{\giss, \symPhaseC}^{\Pss\plusss} }}

\newcommand{\dQtogik}[0]{\textcolor{\paramcolor}{\symDirectionality_{\giss}^{\Qss\plusss}  }}
\newcommand{\dQtogikA}[0]{\textcolor{\paramcolor}{\symDirectionality_{\giss, \symPhaseA}^{\Qss\plusss}  }}
\newcommand{\dQtogikB}[0]{\textcolor{\paramcolor}{\symDirectionality_{\giss, \symPhaseB}^{\Qss\plusss}  }}
\newcommand{\dQtogikC}[0]{\textcolor{\paramcolor}{\symDirectionality_{\giss, \symPhaseC}^{\Qss\plusss}  }}

\newcommand{\dPfromgik}[0]{\textcolor{\paramcolor}{\symDirectionality_{\giss}^{\Pss\subsss} }}
\newcommand{\dPfromgikA}[0]{\textcolor{\paramcolor}{\symDirectionality_{\giss, \symPhaseA}^{\Pss\subsss} }}
\newcommand{\dPfromgikB}[0]{\textcolor{\paramcolor}{\symDirectionality_{\giss, \symPhaseB}^{\Pss\subsss} }}
\newcommand{\dPfromgikC}[0]{\textcolor{\paramcolor}{\symDirectionality_{\giss, \symPhaseC}^{\Pss\subsss} }}

\newcommand{\dQfromgik}[0]{\textcolor{\paramcolor}{\symDirectionality_{\giss}^{\Qss\subsss}  }}
\newcommand{\dQfromgikA}[0]{\textcolor{\paramcolor}{\symDirectionality_{\giss, \symPhaseA}^{\Qss\subsss}  }}
\newcommand{\dQfromgikB}[0]{\textcolor{\paramcolor}{\symDirectionality_{\giss, \symPhaseB}^{\Qss\subsss}  }}
\newcommand{\dQfromgikC}[0]{\textcolor{\paramcolor}{\symDirectionality_{\giss, \symPhaseC}^{\Qss\subsss}  }}

\newcommand{\rhoDODmax}[0]{\textcolor{\paramcolor}{\symLossFactor_{\Ess}^{\maxss}}}
\newcommand{\rhoDODmin}[0]{\textcolor{\paramcolor}{\symLossFactor_{\Ess}^{\minss}}}
\newcommand{\rhoDODeff}[0]{\textcolor{\paramcolor}{\symLossFactor_{\Ess}^{\usabless}}}

\newcommand{\Eijk}[0]{\textcolor{\complexcolor}{\symEnergy_{\indexGridNode \indexGridNodeTwo}^{}  }}
\newcommand{\Eijkre}[0]{\textcolor{black}{\symEnergy_{\indexGridNode \indexGridNodeTwo}^{\realss}  }}

\newcommand{\Sij}[0]{\textcolor{\complexcolor}{\symApparentPower_{\indexGridNode \indexGridNodeTwo}^{} }}

\newcommand{\Pzref}[0]{\textcolor{\paramcolor}{\symPower_{\indexZIP }^{\refss} }}
\newcommand{\Pinode}[0]{\textcolor{black}{\symPower_{\indexGridNode}^{} }}
\newcommand{\Pz}[0]{\textcolor{black}{\symPower_{\indexZIP}^{} }}

\newcommand{\Pij}[0]{\textcolor{black}{\symPower_{\indexGridNode \indexGridNodeTwo}^{} }}
\newcommand{\Pji}[0]{\textcolor{black}{\symPower_{\indexGridNodeTwo \indexGridNode }^{} }}
\newcommand{\Pijloss}[0]{\textcolor{black}{\symPower_{ \indexGridLines}^{\lossss} }}
\newcommand{\Pijmax}[0]{\textcolor{\paramcolor}{\symPower_{\indexGridNode \indexGridNodeTwo}^{\maxss} }}
\newcommand{\Pijacc}[0]{\textcolor{black}{\symPower_{\indexGridNode \indexGridNodeTwo}^{\star} }}
\newcommand{\Pijrated}[0]{\textcolor{\sizingcolor}{\symPower_{\indexGridNode \indexGridNodeTwo}^{\ratedss} }}
\newcommand{\Pjirated}[0]{\textcolor{\sizingcolor}{\symPower_{\indexGridNodeTwo\indexGridNode }^{\ratedss} }}

\newcommand{\Qij}[0]{\textcolor{black}{\symReactivePower_{\indexGridNode \indexGridNodeTwo}^{} }}
\newcommand{\Qji}[0]{\textcolor{black}{\symReactivePower_{\indexGridNodeTwo \indexGridNode }^{} }}
\newcommand{\Qijrated}[0]{\textcolor{\sizingcolor}{\symReactivePower_{\indexGridNode \indexGridNodeTwo}^{\ratedss} }}
\newcommand{\Qjirated}[0]{\textcolor{\sizingcolor}{\symReactivePower_{\indexGridNodeTwo\indexGridNode }^{\ratedss} }}
\newcommand{\Qijmax}[0]{\textcolor{\paramcolor}{\symReactivePower_{\indexGridNode \indexGridNodeTwo}^{\maxss} }}
\newcommand{\Qijmin}[0]{\textcolor{\paramcolor}{\symReactivePower_{\indexGridNode \indexGridNodeTwo}^{\minss} }}
\newcommand{\Qjimax}[0]{\textcolor{\paramcolor}{\symReactivePower_{\indexGridNodeTwo \indexGridNode }^{\maxss} }}
\newcommand{\Qjimin}[0]{\textcolor{\paramcolor}{\symReactivePower_{\indexGridNodeTwo\indexGridNode }^{\minss} }}
\newcommand{\Qijloss}[0]{\textcolor{black}{\symReactivePower_{ \indexGridLines}^{\lossss} }}

\newcommand{\Sjiks}[0]{\textcolor{\complexcolor}{\symApparentPower_{\indexGridNodeTwo \indexGridNode , \seriesss}^{}  }}
\newcommand{\Sijks}[0]{\textcolor{\complexcolor}{\symApparentPower_{\indexGridNode \indexGridNodeTwo, \seriesss}^{}  }}
\newcommand{\Pijks}[0]{\textcolor{black}{\symPower_{\indexGridNode \indexGridNodeTwo, \seriesss}^{}  }}
\newcommand{\Qijks}[0]{\textcolor{black}{\symReactivePower_{\indexGridNode \indexGridNodeTwo, \seriesss}^{}  }}
\newcommand{\Pjiks}[0]{\textcolor{black}{\symPower_{\indexGridNodeTwo \indexGridNode , \seriesss}^{}  }}
\newcommand{\Qjiks}[0]{\textcolor{black}{\symReactivePower_{\indexGridNodeTwo \indexGridNode , \seriesss}^{}  }}

\newcommand{\Sijlossk}[0]{\textcolor{\complexcolor}{\symApparentPower_{ \indexGridLines}^{\lossss}  }}
\newcommand{\Slossk}[0]{\textcolor{black}{\symApparentPower_{}^{\lossss}  }}

\newcommand{\Sijlosssk}[0]{\textcolor{\complexcolor}{\symApparentPower_{ \indexGridLines, \seriesss}^{\lossss}  }}

\newcommand{\Plossk}[0]{\textcolor{black}{\symPower_{}^{\lossss}  }}
\newcommand{\Pijk}[0]{\textcolor{black}{\symPower_{\indexGridNode \indexGridNodeTwo}^{}  }}
\newcommand{\Pijlossk}[0]{\textcolor{black}{\symPower_{ \indexGridLines}^{\lossss}  }}
\newcommand{\Pijlosssk}[0]{\textcolor{black}{\symPower_{ \indexGridLines, \seriesss}^{\lossss}  }}
\newcommand{\Pijlossshk}[0]{\textcolor{black}{\symPower_{ \indexGridLines, \shuntss}^{\lossss}  }}
\newcommand{\Pilossshk}[0]{\textcolor{black}{\symPower_{\indexGridNode \indexGridNodeTwo, \shuntss }^{\lossss}  }}
\newcommand{\Pjlossshk}[0]{\textcolor{black}{\symPower_{ \indexGridNodeTwo\indexGridNode, \shuntss}^{\lossss}  }}
\newcommand{\Pjik}[0]{\textcolor{black}{\symPower_{\indexGridNodeTwo \indexGridNode}^{}    }}
\newcommand{\Pinodek}[0]{\textcolor{black}{\symPower_{ \indexGridNode}^{}    }}
\newcommand{\Qijk}[0]{\textcolor{black}{\symReactivePower_{\indexGridNode \indexGridNodeTwo}^{}  }}
\newcommand{\Qijkdelta}[0]{\textcolor{black}{\symReactivePower_{\indexGridNode \indexGridNodeTwo}^{\Delta}  }}
\newcommand{\Qijlossk}[0]{\textcolor{black}{\symReactivePower_{ \indexGridLines}^{\lossss}  }}
\newcommand{\Qijlosssk}[0]{\textcolor{black}{\symReactivePower_{ \indexGridLines,\seriesss}^{\lossss}  }}
\newcommand{\Qijlossshk}[0]{\textcolor{black}{\symReactivePower_{ \indexGridLines,\shuntss}^{\lossss}  }}
\newcommand{\Qilossshk}[0]{\textcolor{black}{\symReactivePower_{\indexGridNode\indexGridNodeTwo,\shuntss }^{\lossss}  }}
\newcommand{\Qjlossshk}[0]{\textcolor{black}{\symReactivePower_{ \indexGridNodeTwo\indexGridNode,\shuntss}^{\lossss}  }}
\newcommand{\Qjik}[0]{\textcolor{black}{\symReactivePower_{\indexGridNodeTwo \indexGridNode}^{}    }}
\newcommand{\Qinodek}[0]{\textcolor{black}{\symReactivePower_{ \indexGridNode}^{}    }}
\newcommand{\Sinodek}[0]{\textcolor{\complexcolor}{\symApparentPower_{ \indexGridNode}^{}    }}

\newcommand{\Wijk}[0]{\textcolor{\complexcolor}{W_{\indexGridNode \indexGridNodeTwo}^{}  }} 
\newcommand{\WijkSDP}[0]{\textcolor{\complexcolor}{\mathbf{W}_{}^{}  }} 
\newcommand{\Wiik}[0]{\textcolor{black}{W_{\indexGridNode \indexGridNode}^{}  }} 
\newcommand{\Wjjk}[0]{\textcolor{black}{W_{\indexGridNodeTwo \indexGridNodeTwo}^{}  }} 
\newcommand{\Wjik}[0]{\textcolor{\complexcolor}{W_{\indexGridNodeTwo \indexGridNode}^{}  }} 

\newcommand{\Rijk}[0]{\textcolor{black}{R_{\indexGridNode \indexGridNodeTwo}^{}  }} 
\newcommand{\Tijk}[0]{\textcolor{black}{T_{\indexGridNode \indexGridNodeTwo}^{}  }}  
\newcommand{\Rjik}[0]{\textcolor{black}{R_{\indexGridNodeTwo\indexGridNode }^{}  }} 
\newcommand{\Tjik}[0]{\textcolor{black}{T_{\indexGridNodeTwo\indexGridNode }^{}  }} 
\newcommand{\Rlk}[0]{\textcolor{black}{R_{\indexGridLines }^{}  }}    
\newcommand{\rijk}[0]{\textcolor{black}{r_{\indexGridNode \indexGridNodeTwo}^{}  }}   

\newcommand{\SikA}[0]{\textcolor{\complexcolor}{\symApparentPower_{\indexGridNode , \symPhaseA}^{}  }}
\newcommand{\SikB}[0]{\textcolor{\complexcolor}{\symApparentPower_{\indexGridNode , \symPhaseB}^{}  }}
\newcommand{\SikC}[0]{\textcolor{\complexcolor}{\symApparentPower_{\indexGridNode , \symPhaseC}^{}  }}

\newcommand{\PijkA}[0]{\textcolor{black}{\symPower_{\indexGridNode \indexGridNodeTwo, \symPhaseA}^{}  }}
\newcommand{\PijkB}[0]{\textcolor{black}{\symPower_{\indexGridNode \indexGridNodeTwo, \symPhaseB}^{}  }}
\newcommand{\PijkC}[0]{\textcolor{black}{\symPower_{\indexGridNode \indexGridNodeTwo, \symPhaseC}^{}  }}
\newcommand{\QijkA}[0]{\textcolor{black}{\symReactivePower_{\indexGridNode \indexGridNodeTwo, \symPhaseA}^{}  }}
\newcommand{\QijkB}[0]{\textcolor{black}{\symReactivePower_{\indexGridNode \indexGridNodeTwo, \symPhaseB}^{}  }}
\newcommand{\QijkC}[0]{\textcolor{black}{\symReactivePower_{\indexGridNode \indexGridNodeTwo, \symPhaseC}^{}  }}

\newcommand{\PjikA}[0]{\textcolor{black}{\symPower_{ \indexGridNodeTwo\indexGridNode, \symPhaseA}^{}  }}
\newcommand{\PjikB}[0]{\textcolor{black}{\symPower_{ \indexGridNodeTwo\indexGridNode, \symPhaseB}^{}  }}
\newcommand{\PjikC}[0]{\textcolor{black}{\symPower_{ \indexGridNodeTwo\indexGridNode, \symPhaseC}^{}  }}
\newcommand{\QjikA}[0]{\textcolor{black}{\symReactivePower_{ \indexGridNodeTwo\indexGridNode, \symPhaseA}^{}  }}
\newcommand{\QjikB}[0]{\textcolor{black}{\symReactivePower_{ \indexGridNodeTwo\indexGridNode, \symPhaseB}^{}  }}
\newcommand{\QjikC}[0]{\textcolor{black}{\symReactivePower_{ \indexGridNodeTwo\indexGridNode, \symPhaseC}^{}  }}

\newcommand{\PijlosskA}[0]{\textcolor{black}{\symPower_{\indexGridNode \indexGridNodeTwo,\symPhaseA}^{\lossss}  }}
\newcommand{\PijlosskB}[0]{\textcolor{black}{\symPower_{\indexGridNode \indexGridNodeTwo,\symPhaseB}^{\lossss}  }}
\newcommand{\PijlosskC}[0]{\textcolor{black}{\symPower_{\indexGridNode \indexGridNodeTwo,\symPhaseC}^{\lossss}  }}
\newcommand{\QijlosskA}[0]{\textcolor{black}{\symReactivePower_{\indexGridNode \indexGridNodeTwo,\symPhaseA}^{\lossss}  }}
\newcommand{\QijlosskB}[0]{\textcolor{black}{\symReactivePower_{\indexGridNode \indexGridNodeTwo,\symPhaseB}^{\lossss}  }}
\newcommand{\QijlosskC}[0]{\textcolor{black}{\symReactivePower_{\indexGridNode \indexGridNodeTwo,\symPhaseC}^{\lossss}  }}

\newcommand{\PijlosskminA}[0]{\textcolor{\boundscolor}{\symPower_{\indexGridNode \indexGridNodeTwo,\symPhaseA}^{\lossss,\minss}  }}
\newcommand{\PijlosskminB}[0]{\textcolor{\boundscolor}{\symPower_{\indexGridNode \indexGridNodeTwo,\symPhaseB}^{\lossss,\minss}  }}
\newcommand{\PijlosskminC}[0]{\textcolor{\boundscolor}{\symPower_{\indexGridNode \indexGridNodeTwo,\symPhaseC}^{\lossss,\minss}  }}
\newcommand{\QijlosskminA}[0]{\textcolor{\boundscolor}{\symReactivePower_{\indexGridNode \indexGridNodeTwo,\symPhaseA}^{\lossss,\minss}  }}
\newcommand{\QijlosskminB}[0]{\textcolor{\boundscolor}{\symReactivePower_{\indexGridNode \indexGridNodeTwo,\symPhaseB}^{\lossss,\minss}  }}
\newcommand{\QijlosskminC}[0]{\textcolor{\boundscolor}{\symReactivePower_{\indexGridNode \indexGridNodeTwo,\symPhaseC}^{\lossss,\minss}  }}

\newcommand{\PijlosskmaxA}[0]{\textcolor{\boundscolor}{\symPower_{\indexGridNode \indexGridNodeTwo,\symPhaseA}^{\lossss,\maxss}  }}
\newcommand{\PijlosskmaxB}[0]{\textcolor{\boundscolor}{\symPower_{\indexGridNode \indexGridNodeTwo,\symPhaseB}^{\lossss,\maxss}  }}
\newcommand{\PijlosskmaxC}[0]{\textcolor{\boundscolor}{\symPower_{\indexGridNode \indexGridNodeTwo,\symPhaseC}^{\lossss,\maxss}  }}
\newcommand{\QijlosskmaxA}[0]{\textcolor{\boundscolor}{\symReactivePower_{\indexGridNode \indexGridNodeTwo,\symPhaseA}^{\lossss,\maxss}  }}
\newcommand{\QijlosskmaxB}[0]{\textcolor{\boundscolor}{\symReactivePower_{\indexGridNode \indexGridNodeTwo,\symPhaseB}^{\lossss,\maxss}  }}
\newcommand{\QijlosskmaxC}[0]{\textcolor{\boundscolor}{\symReactivePower_{\indexGridNode \indexGridNodeTwo,\symPhaseC}^{\lossss,\maxss}  }}

\newcommand{\Iijkabs}[0]{\textcolor{black}{|\symCurrent_{\indexGridLines\indexGridNode \indexGridNodeTwo}^{}  |}}
\newcommand{\Sijkabs}[0]{\textcolor{black}{|\symApparentPower_{\indexGridLines\indexGridNode \indexGridNodeTwo}^{}  |}}
\newcommand{\Sijk}[0]{\textcolor{\complexcolor}{\symApparentPower_{\indexGridLines\indexGridNode \indexGridNodeTwo}^{}  }}

\newcommand{\SijkSDPrated}[0]{\textcolor{\paramcolor}{\mathbf{\symApparentPower}_{\indexGridLines\indexGridNode \indexGridNodeTwo}^{\ratedss}  }}
\newcommand{\IijkSDPrated}[0]{\textcolor{\paramcolor}{\mathbf{\symCurrent}_{\indexGridLines\indexGridNode \indexGridNodeTwo}^{\ratedss}  }}
\newcommand{\IijksSDPrated}[0]{\textcolor{\paramcolor}{\mathbf{\symCurrent}_{\indexGridLines }^{\seriesss, \ratedss}  }}

\newcommand{\SijkSDPseq}[0]{\textcolor{\complexcolor}{\mathbf{\symApparentPower}_{\indexGridLines\indexGridNode \indexGridNodeTwo}^{\fortescuess}  }}
\newcommand{\SjikSDPseq}[0]{\textcolor{\complexcolor}{\mathbf{\symApparentPower}_{\indexGridLines\indexGridNodeTwo \indexGridNode }^{\fortescuess} }}
\newcommand{\SijkSDPH}[0]{\textcolor{\complexcolor}{(\mathbf{\symApparentPower}_{\indexGridLines\indexGridNode \indexGridNodeTwo})^{\hermitiantranspose}  }}
\newcommand{\SijkSDP}[0]{\textcolor{\complexcolor}{\mathbf{\symApparentPower}_{\indexGridLines\indexGridNode \indexGridNodeTwo}^{}  }}
\newcommand{\SjikSDP}[0]{\textcolor{\complexcolor}{\mathbf{\symApparentPower}_{\indexGridLines\indexGridNodeTwo \indexGridNode }^{}  }}
\newcommand{\SjkkSDP}[0]{\textcolor{\complexcolor}{\mathbf{\symApparentPower}_{ \indexGridLines\indexGridNodeTwo \indexGridNodeThree}^{}  }}
\newcommand{\PjkkSDP}[0]{\textcolor{black}{\mathbf{\symPower}_{ \indexGridLines\indexGridNodeTwo \indexGridNodeThree}^{}  }}
\newcommand{\QjkkSDP}[0]{\textcolor{black}{\mathbf{\symReactivePower}_{ \indexGridLines\indexGridNodeTwo \indexGridNodeThree}^{}  }}

\newcommand{\SijkSDPAdot}[0]{\textcolor{\complexcolor}{\mathbf{\symApparentPower}_{\indexGridLines\indexGridNode \indexGridNodeTwo}^{\symPhaseA \! \cdot}  }}
\newcommand{\SijkSDPBdot}[0]{\textcolor{\complexcolor}{\mathbf{\symApparentPower}_{\indexGridLines\indexGridNode \indexGridNodeTwo}^{\symPhaseB \! \cdot}  }}
\newcommand{\SijkSDPCdot}[0]{\textcolor{\complexcolor}{\mathbf{\symApparentPower}_{\indexGridLines\indexGridNode \indexGridNodeTwo}^{\symPhaseC \! \cdot}  }}
\newcommand{\SijkSDPNdot}[0]{\textcolor{\complexcolor}{\mathbf{\symApparentPower}_{\indexGridLines\indexGridNode \indexGridNodeTwo}^{\symPhaseN \! \cdot}  }}
\newcommand{\SijkSDPGdot}[0]{\textcolor{\complexcolor}{\mathbf{\symApparentPower}_{\indexGridLines\indexGridNode \indexGridNodeTwo}^{\symPhaseG \! \cdot}  }}

\newcommand{\SjikSDPAdot}[0]{\textcolor{\complexcolor}{\mathbf{\symApparentPower}_{\indexGridLines\indexGridNodeTwo\indexGridNode }^{\symPhaseA \! \cdot}  }}
\newcommand{\SjikSDPBdot}[0]{\textcolor{\complexcolor}{\mathbf{\symApparentPower}_{\indexGridLines\indexGridNodeTwo\indexGridNode }^{\symPhaseB \! \cdot}  }}
\newcommand{\SjikSDPCdot}[0]{\textcolor{\complexcolor}{\mathbf{\symApparentPower}_{\indexGridLines\indexGridNodeTwo\indexGridNode }^{\symPhaseC \! \cdot}  }}
\newcommand{\SjikSDPNdot}[0]{\textcolor{\complexcolor}{\mathbf{\symApparentPower}_{\indexGridLines\indexGridNodeTwo\indexGridNode }^{\symPhaseN \! \cdot}  }}

\newcommand{\PijkSDPAdot}[0]{\textcolor{black}{\mathbf{\symPower}_{\indexGridLines\indexGridNode \indexGridNodeTwo}^{\symPhaseA \! \cdot}  }}
\newcommand{\PijkSDPBdot}[0]{\textcolor{black}{\mathbf{\symPower}_{\indexGridLines\indexGridNode \indexGridNodeTwo}^{\symPhaseB \! \cdot}  }}
\newcommand{\PijkSDPCdot}[0]{\textcolor{black}{\mathbf{\symPower}_{\indexGridLines\indexGridNode \indexGridNodeTwo}^{\symPhaseC \! \cdot}  }}
\newcommand{\PijkSDPNdot}[0]{\textcolor{black}{\mathbf{\symPower}_{\indexGridLines\indexGridNode \indexGridNodeTwo}^{\symPhaseN \! \cdot}  }}

\newcommand{\QijkSDPAdot}[0]{\textcolor{black}{\mathbf{\symReactivePower}_{\indexGridLines\indexGridNode \indexGridNodeTwo}^{\symPhaseA \! \cdot}  }}
\newcommand{\QijkSDPBdot}[0]{\textcolor{black}{\mathbf{\symReactivePower}_{\indexGridLines\indexGridNode \indexGridNodeTwo}^{\symPhaseB \! \cdot}  }}
\newcommand{\QijkSDPCdot}[0]{\textcolor{black}{\mathbf{\symReactivePower}_{\indexGridLines\indexGridNode \indexGridNodeTwo}^{\symPhaseC \! \cdot}  }}
\newcommand{\QijkSDPNdot}[0]{\textcolor{black}{\mathbf{\symReactivePower}_{\indexGridLines\indexGridNode \indexGridNodeTwo}^{\symPhaseN \! \cdot}  }}

\newcommand{\PijksSDPAdot}[0]{\textcolor{black}{\mathbf{\symPower}_{\indexGridLines\indexGridNode \indexGridNodeTwo}^{\seriesss,\symPhaseA \! \cdot}  }}
\newcommand{\PijksSDPBdot}[0]{\textcolor{black}{\mathbf{\symPower}_{\indexGridLines\indexGridNode \indexGridNodeTwo}^{\seriesss,\symPhaseB \! \cdot}  }}
\newcommand{\PijksSDPCdot}[0]{\textcolor{black}{\mathbf{\symPower}_{\indexGridLines\indexGridNode \indexGridNodeTwo}^{\seriesss,\symPhaseC \! \cdot}  }}
\newcommand{\PijksSDPNdot}[0]{\textcolor{black}{\mathbf{\symPower}_{\indexGridLines\indexGridNode \indexGridNodeTwo}^{\seriesss,\symPhaseN \! \cdot}  }}

\newcommand{\QijksSDPAdot}[0]{\textcolor{black}{\mathbf{\symReactivePower}_{\indexGridLines\indexGridNode \indexGridNodeTwo}^{\seriesss,\symPhaseA \! \cdot}  }}
\newcommand{\QijksSDPBdot}[0]{\textcolor{black}{\mathbf{\symReactivePower}_{\indexGridLines\indexGridNode \indexGridNodeTwo}^{\seriesss,\symPhaseB \! \cdot}  }}
\newcommand{\QijksSDPCdot}[0]{\textcolor{black}{\mathbf{\symReactivePower}_{\indexGridLines\indexGridNode \indexGridNodeTwo}^{\seriesss,\symPhaseC \! \cdot}  }}
\newcommand{\QijksSDPNdot}[0]{\textcolor{black}{\mathbf{\symReactivePower}_{\indexGridLines\indexGridNode \indexGridNodeTwo}^{\seriesss,\symPhaseN \! \cdot}  }}

\newcommand{\SijkSDPdotA}[0]{\textcolor{\complexcolor}{\mathbf{\symApparentPower}_{\indexGridLines\indexGridNode \indexGridNodeTwo}^{\cdot\!\symPhaseA}  }}
\newcommand{\SijkSDPdotB}[0]{\textcolor{\complexcolor}{\mathbf{\symApparentPower}_{\indexGridLines\indexGridNode \indexGridNodeTwo}^{\cdot\!\symPhaseB }  }}
\newcommand{\SijkSDPdotC}[0]{\textcolor{\complexcolor}{\mathbf{\symApparentPower}_{\indexGridLines\indexGridNode \indexGridNodeTwo}^{\cdot\!\symPhaseC }  }}
\newcommand{\SijkSDPdotN}[0]{\textcolor{\complexcolor}{\mathbf{\symApparentPower}_{\indexGridLines\indexGridNode \indexGridNodeTwo}^{\cdot \!\symPhaseN }  }}

\newcommand{\PijkSDPdotA}[0]{\textcolor{black}{\mathbf{\symPower}_{\indexGridLines\indexGridNode \indexGridNodeTwo}^{\cdot\!\symPhaseA}  }}
\newcommand{\PijkSDPdotB}[0]{\textcolor{black}{\mathbf{\symPower}_{\indexGridLines\indexGridNode \indexGridNodeTwo}^{\cdot\!\symPhaseB }  }}
\newcommand{\PijkSDPdotC}[0]{\textcolor{black}{\mathbf{\symPower}_{\indexGridLines\indexGridNode \indexGridNodeTwo}^{\cdot\!\symPhaseC }  }}
\newcommand{\PijkSDPdotN}[0]{\textcolor{black}{\mathbf{\symPower}_{\indexGridLines\indexGridNode \indexGridNodeTwo}^{\cdot \!\symPhaseN }  }}

\newcommand{\QijkSDPdotA}[0]{\textcolor{black}{\mathbf{\symReactivePower}_{\indexGridLines\indexGridNode \indexGridNodeTwo}^{\cdot\!\symPhaseA}  }}
\newcommand{\QijkSDPdotB}[0]{\textcolor{black}{\mathbf{\symReactivePower}_{\indexGridLines\indexGridNode \indexGridNodeTwo}^{\cdot\!\symPhaseB }  }}
\newcommand{\QijkSDPdotC}[0]{\textcolor{black}{\mathbf{\symReactivePower}_{\indexGridLines\indexGridNode \indexGridNodeTwo}^{\cdot\!\symPhaseC }  }}
\newcommand{\QijkSDPdotN}[0]{\textcolor{black}{\mathbf{\symReactivePower}_{\indexGridLines\indexGridNode \indexGridNodeTwo}^{\cdot \!\symPhaseN }  }}

\newcommand{\PijksSDPdotA}[0]{\textcolor{black}{\mathbf{\symPower}_{\indexGridLines\indexGridNode \indexGridNodeTwo}^{\seriesss,\cdot\!\symPhaseA}  }}
\newcommand{\PijksSDPdotB}[0]{\textcolor{black}{\mathbf{\symPower}_{\indexGridLines\indexGridNode \indexGridNodeTwo}^{\seriesss,\cdot\!\symPhaseB }  }}
\newcommand{\PijksSDPdotC}[0]{\textcolor{black}{\mathbf{\symPower}_{\indexGridLines\indexGridNode \indexGridNodeTwo}^{\seriesss,\cdot\!\symPhaseC }  }}
\newcommand{\PijksSDPdotN}[0]{\textcolor{black}{\mathbf{\symPower}_{\indexGridLines\indexGridNode \indexGridNodeTwo}^{\seriesss,\cdot \!\symPhaseN }  }}

\newcommand{\QijksSDPdotA}[0]{\textcolor{black}{\mathbf{\symReactivePower}_{\indexGridLines\indexGridNode \indexGridNodeTwo}^{\seriesss,\cdot\!\symPhaseA}  }}
\newcommand{\QijksSDPdotB}[0]{\textcolor{black}{\mathbf{\symReactivePower}_{\indexGridLines\indexGridNode \indexGridNodeTwo}^{\seriesss,\cdot\!\symPhaseB }  }}
\newcommand{\QijksSDPdotC}[0]{\textcolor{black}{\mathbf{\symReactivePower}_{\indexGridLines\indexGridNode \indexGridNodeTwo}^{\seriesss,\cdot\!\symPhaseC }  }}
\newcommand{\QijksSDPdotN}[0]{\textcolor{black}{\mathbf{\symReactivePower}_{\indexGridLines\indexGridNode \indexGridNodeTwo}^{\seriesss,\cdot \!\symPhaseN }  }}

\newcommand{\SijskSDPAdot}[0]{\textcolor{\complexcolor}{\mathbf{\symApparentPower}_{\indexGridLines\indexGridNode \indexGridNodeTwo}^{\seriesss, \symPhaseA \! \cdot}  }}
\newcommand{\SijskSDPBdot}[0]{\textcolor{\complexcolor}{\mathbf{\symApparentPower}_{\indexGridLines\indexGridNode \indexGridNodeTwo}^{\seriesss, \symPhaseB \! \cdot}  }}
\newcommand{\SijskSDPCdot}[0]{\textcolor{\complexcolor}{\mathbf{\symApparentPower}_{\indexGridLines\indexGridNode \indexGridNodeTwo}^{\seriesss, \symPhaseC \! \cdot}  }}
\newcommand{\SijskSDPNdot}[0]{\textcolor{\complexcolor}{\mathbf{\symApparentPower}_{\indexGridLines\indexGridNode \indexGridNodeTwo}^{\seriesss, \symPhaseN \! \cdot}  }}
\newcommand{\SijskSDPGdot}[0]{\textcolor{\complexcolor}{\mathbf{\symApparentPower}_{\indexGridLines\indexGridNode \indexGridNodeTwo}^{\seriesss, \symPhaseG \! \cdot}  }}

\newcommand{\SijskSDPdotA}[0]{\textcolor{\complexcolor}{\mathbf{\symApparentPower}_{\indexGridLines\indexGridNode \indexGridNodeTwo}^{\seriesss, \cdot\!\symPhaseA}  }}
\newcommand{\SijskSDPdotB}[0]{\textcolor{\complexcolor}{\mathbf{\symApparentPower}_{\indexGridLines\indexGridNode \indexGridNodeTwo}^{\seriesss, \cdot\!\symPhaseB }  }}
\newcommand{\SijskSDPdotC}[0]{\textcolor{\complexcolor}{\mathbf{\symApparentPower}_{\indexGridLines\indexGridNode \indexGridNodeTwo}^{\seriesss, \cdot\!\symPhaseC }  }}
\newcommand{\SijskSDPdotN}[0]{\textcolor{\complexcolor}{\mathbf{\symApparentPower}_{\indexGridLines\indexGridNode \indexGridNodeTwo}^{\seriesss, \cdot \!\symPhaseN }  }}

\newcommand{\SijksSDPmax}[0]{\textcolor{\paramcolor}{\mathbf{\symApparentPower}_{\indexGridLines \indexGridNode \indexGridNodeTwo}^{\seriesss, \maxss}  }}

\newcommand{\SijksSDP}[0]{\textcolor{\complexcolor}{\mathbf{\symApparentPower}_{\indexGridLines \indexGridNode \indexGridNodeTwo}^{\seriesss}  }}
\newcommand{\SjiksSDP}[0]{\textcolor{\complexcolor}{\mathbf{\symApparentPower}_{\indexGridLines \indexGridNodeTwo\indexGridNode }^{\seriesss}  }}

\newcommand{\SijksSDPH}[0]{\textcolor{\complexcolor}{\mathbf{(\symApparentPower}_{\indexGridLines \indexGridNode \indexGridNodeTwo}^{\seriesss})^{\hermitiantranspose}  }}

\newcommand{\SijksSDPdiag}[0]{\textcolor{\complexcolor}{\mathbf{\symApparentPower}_{\indexGridLines\indexGridNode \indexGridNodeTwo}^{\seriesss, \text{diag}}  }}
\newcommand{\PijksSDPdiag}[0]{\textcolor{black}{\mathbf{\symPower}_{\indexGridLines\indexGridNode \indexGridNodeTwo}^{\seriesss, \text{diag}}  }}
\newcommand{\QijksSDPdiag}[0]{\textcolor{black}{\mathbf{\symReactivePower}_{\indexGridLines\indexGridNode \indexGridNodeTwo}^{\seriesss, \text{diag}}  }}

\newcommand{\PijkSDP}[0]{\textcolor{black}{\mathbf{\symPower}_{\indexGridLines\indexGridNode \indexGridNodeTwo}^{}  }}
\newcommand{\QijkSDP}[0]{\textcolor{black}{\mathbf{\symReactivePower}_{\indexGridLines\indexGridNode \indexGridNodeTwo}^{}  }}
\newcommand{\PijkSDPH}[0]{\textcolor{black}{\mathbf{\symPower}_{\indexGridLines\indexGridNode \indexGridNodeTwo}^{\hermitiantranspose}  }}
\newcommand{\QijkSDPH}[0]{\textcolor{black}{\mathbf{\symReactivePower}_{\indexGridLines\indexGridNode \indexGridNodeTwo}^{\hermitiantranspose}  }}
\newcommand{\PjikSDP}[0]{\textcolor{black}{\mathbf{\symPower}_{\indexGridLines\indexGridNodeTwo\indexGridNode }^{}  }}
\newcommand{\QjikSDP}[0]{\textcolor{black}{\mathbf{\symReactivePower}_{\indexGridLines\indexGridNodeTwo\indexGridNode }^{}  }}

\newcommand{\PijksSDPH}[0]{\textcolor{black}{\mathbf{\symPower}_{\indexGridLines\indexGridNode \indexGridNodeTwo, \seriesss}^{\hermitiantranspose}  }}
\newcommand{\QijksSDPH}[0]{\textcolor{black}{\mathbf{\symReactivePower}_{\indexGridLines\indexGridNode \indexGridNodeTwo, \seriesss}^{\hermitiantranspose}  }}

\newcommand{\PijksSDP}[0]{\textcolor{black}{\mathbf{\symPower}_{\indexGridLines\indexGridNode \indexGridNodeTwo}^{\seriesss}  }}
\newcommand{\QijksSDP}[0]{\textcolor{black}{\mathbf{\symReactivePower}_{\indexGridLines\indexGridNode \indexGridNodeTwo}^{\seriesss}  }}
\newcommand{\PjiksSDP}[0]{\textcolor{black}{\mathbf{\symPower}_{\indexGridLines\indexGridNodeTwo\indexGridNode }^{\seriesss}  }}
\newcommand{\QjiksSDP}[0]{\textcolor{black}{\mathbf{\symReactivePower}_{\indexGridLines\indexGridNodeTwo\indexGridNode }^{\seriesss}  }}

\newcommand{\SijklossSDP}[0]{\textcolor{\complexcolor}{\mathbf{\symApparentPower}_{\indexGridNode \indexGridNodeTwo}^{\lossss}  }}
\newcommand{\SjiklossSDP}[0]{\textcolor{\complexcolor}{\mathbf{\symApparentPower}_{\indexGridNodeTwo\indexGridNode }^{\lossss}  }}
\newcommand{\SlklossSDP}[0]{\textcolor{\complexcolor}{\mathbf{\symApparentPower}_{\indexGridLines}^{\lossss}  }}
\newcommand{\PijklossSDP}[0]{\textcolor{black}{\mathbf{\symPower}_{\indexGridNode \indexGridNodeTwo}^{\lossss}  }}
\newcommand{\QijklossSDP}[0]{\textcolor{black}{\mathbf{\symReactivePower}_{\indexGridNode \indexGridNodeTwo}^{\lossss}  }}

\newcommand{\PijkslossSDP}[0]{\textcolor{black}{\mathbf{\symPower}_{\indexGridNode \indexGridNodeTwo, \seriesss}^{\lossss}  }}
\newcommand{\QijkslossSDP}[0]{\textcolor{black}{\mathbf{\symReactivePower}_{\indexGridNode \indexGridNodeTwo, \seriesss}^{\lossss}  }}

\newcommand{\SijkshlossSDP}[0]{\textcolor{\complexcolor}{\mathbf{\symApparentPower}_{\indexGridLines\indexGridNode \indexGridNodeTwo}^{\lossss, \shuntss}  }}
\newcommand{\SjikshlossSDP}[0]{\textcolor{\complexcolor}{\mathbf{\symApparentPower}_{\indexGridLines\indexGridNodeTwo\indexGridNode}^{\lossss, \shuntss}  }}
\newcommand{\SlkshlossSDP}[0]{\textcolor{\complexcolor}{\mathbf{\symApparentPower}_{\indexGridLines}^{\lossss, \shuntss}  }}

\newcommand{\SijkslossSDP}[0]{\textcolor{\complexcolor}{\mathbf{\symApparentPower}_{\indexGridLines }^{\lossss,\seriesss}  }}
\newcommand{\SjikslossSDP}[0]{\textcolor{\complexcolor}{\mathbf{\symApparentPower}_{\indexGridLines}^{\lossss,\seriesss}  }}
\newcommand{\SlkslossSDP}[0]{\textcolor{\complexcolor}{\mathbf{\symApparentPower}_{\indexGridLines}^{\lossss,\seriesss}  }}

\newcommand{\Sinode}[0]{\textcolor{black}{\symApparentPower_{\indexGridNode }^{}  }}
\newcommand{\Sjnode}[0]{\textcolor{black}{\symApparentPower_{ \indexGridNodeTwo}^{}  }}
\newcommand{\SjnodeSDP}[0]{\textcolor{\complexcolor}{\mathbf{\symApparentPower}_{ \indexGridNodeTwo}^{}  }}
\newcommand{\QjnodeSDP}[0]{\textcolor{black}{\mathbf{\symReactivePower}_{ \indexGridNodeTwo}^{}  }}
\newcommand{\PjnodeSDP}[0]{\textcolor{black}{\mathbf{\symPower}_{ \indexGridNodeTwo}^{}  }}

\newcommand{\SinodeSDP}[0]{\textcolor{\complexcolor}{\mathbf{\symApparentPower}_{ \indexGridNode}^{}  }}
\newcommand{\QinodeSDP}[0]{\textcolor{black}{\mathbf{\symReactivePower}_{ \indexGridNode}^{}  }}
\newcommand{\PinodeSDP}[0]{\textcolor{black}{\mathbf{\symPower}_{ \indexGridNode}^{}  }}

\newcommand{\IinodeSDP}[0]{\textcolor{\complexcolor}{\mathbf{\symCurrent}_{ \indexGridNode}^{}  }}

\newcommand{\IuunitSDPrated}[0]{\textcolor{\paramcolor}{\mathbf{\symCurrent}_{ \indexUnit}^{\ratedss}  }}
\newcommand{\IuunitSDPdeltarated}[0]{\textcolor{\paramcolor}{\mathbf{\symCurrent}_{ \indexUnit}^{\Delta\ratedss}  }}

\newcommand{\IuunitSDPref}[0]{\textcolor{\complexparamcolor}{\mathbf{\symCurrent}_{ \indexUnit}^{\refss}  }}

\newcommand{\IuunitSDP}[0]{\textcolor{\complexcolor}{\mathbf{\symCurrent}_{ \indexUnit}^{}  }}
\newcommand{\IuunitSDPreal}[0]{\textcolor{black}{\mathbf{\symCurrent}_{ \indexUnit}^{\text{re}}  }}
\newcommand{\IuunitSDPimag}[0]{\textcolor{black}{\mathbf{\symCurrent}_{ \indexUnit}^{\text{im}}  }}

\newcommand{\IuunitSDPdelta}[0]{\textcolor{\complexcolor}{\mathbf{\symCurrent}_{ \indexUnit}^{\Delta}  }}

\newcommand{\UuunitSDP}[0]{\textcolor{\complexcolor}{\mathbf{\symVoltage}_{ \indexUnit}^{}  }}

\newcommand{\IbSDP}[0]{\textcolor{\complexcolor}{\mathbf{\symCurrent}_{ \indexShunt}^{}  }}
\newcommand{\IbSDPreal}[0]{\textcolor{black}{\mathbf{\symCurrent}_{ \indexShunt}^{\text{re}}  }}
\newcommand{\IbSDPimag}[0]{\textcolor{black}{\mathbf{\symCurrent}_{ \indexShunt}^{\text{im}}  }}

\newcommand{\SbSDP}[0]{\textcolor{\complexcolor}{\mathbf{\symApparentPower}_{ \indexShunt} }}
\newcommand{\PbSDP}[0]{\textcolor{black}{\mathbf{\symPower}_{ \indexShunt} }}
\newcommand{\QbSDP}[0]{\textcolor{black}{\mathbf{\symReactivePower}_{ \indexShunt} }}

\newcommand{\SuunitSDPint}[0]{\textcolor{\complexcolor}{\mathbf{\symApparentPower}_{ \indexUnit}^{'}  }}

\newcommand{\SuunitSDPref}[0]{\textcolor{\complexparamcolor}{\mathbf{\symApparentPower}_{ \indexUnit}^{\refss}  }}
\newcommand{\SuunitSDPrefdelta}[0]{\textcolor{\complexparamcolor}{\mathbf{\symApparentPower}_{ \indexUnit}^{\Delta\refss}  }}

\newcommand{\PuunitSDPmin}[0]{\textcolor{\paramcolor}{\mathbf{\symPower}_{ \indexUnit}^{\minss}  }}
\newcommand{\QuunitSDPmin}[0]{\textcolor{\paramcolor}{\mathbf{\symReactivePower}_{ \indexUnit}^{\minss}  }}
\newcommand{\PuunitSDPmax}[0]{\textcolor{\paramcolor}{\mathbf{\symPower}_{ \indexUnit}^{\maxss}  }}
\newcommand{\QuunitSDPmax}[0]{\textcolor{\paramcolor}{\mathbf{\symReactivePower}_{ \indexUnit}^{\maxss}  }}

\newcommand{\PuunitSDPmindelta}[0]{\textcolor{\paramcolor}{\mathbf{\symPower}_{ \indexUnit}^{\Delta\minss}  }}
\newcommand{\QuunitSDPmindelta}[0]{\textcolor{\paramcolor}{\mathbf{\symReactivePower}_{ \indexUnit}^{\Delta\minss}  }}
\newcommand{\PuunitSDPmaxdelta}[0]{\textcolor{\paramcolor}{\mathbf{\symPower}_{ \indexUnit}^{\Delta\maxss}  }}
\newcommand{\QuunitSDPmaxdelta}[0]{\textcolor{\paramcolor}{\mathbf{\symReactivePower}_{ \indexUnit}^{\Delta\maxss}  }}

\newcommand{\SuunitSDP}[0]{\textcolor{\complexcolor}{\mathbf{\symApparentPower}_{ \indexUnit}^{}  }}
\newcommand{\PuunitSDP}[0]{\textcolor{black}{\mathbf{\symPower}_{ \indexUnit}^{}  }}
\newcommand{\QuunitSDP}[0]{\textcolor{black}{\mathbf{\symReactivePower}_{ \indexUnit}^{}  }}

\newcommand{\SuunitSDPdelta}[0]{\textcolor{\complexcolor}{\mathbf{\symApparentPower}_{ \indexUnit}^{\Delta}  }}
\newcommand{\PuunitSDPdelta}[0]{\textcolor{black}{\mathbf{\symPower}_{ \indexUnit}^{\Delta}  }}
\newcommand{\QuunitSDPdelta}[0]{\textcolor{black}{\mathbf{\symReactivePower}_{ \indexUnit}^{\Delta}  }}

\newcommand{\Tuunitdelta}[0]{\textcolor{\paramcolor}{\mathbf{T}^{\Delta}  }}
\newcommand{\Nwye}[0]{\textcolor{\paramcolor}{\mathbf{N}^{\text{pn}}  }}

\newcommand{\XuunitSDPdelta}[0]{\textcolor{\complexcolor}{\mathbf{X}_{ \indexUnit}^{\Delta}  }}
\newcommand{\XuunitSDPdeltareal}[0]{\textcolor{black}{\mathbf{X}_{ \indexUnit}^{\Delta\realss}  }}
\newcommand{\XuunitSDPdeltaimag}[0]{\textcolor{black}{\mathbf{X}_{ \indexUnit}^{\Delta\imagss}  }}

\newcommand{\IuunitAN}[0]{\textcolor{\complexcolor}{{\symCurrent}_{ \indexUnit, \symPhaseA}^{}  }}
\newcommand{\IuunitBN}[0]{\textcolor{\complexcolor}{{\symCurrent}_{ \indexUnit, \symPhaseB}^{}  }}
\newcommand{\IuunitCN}[0]{\textcolor{\complexcolor}{{\symCurrent}_{ \indexUnit, \symPhaseC}^{}  }}
\newcommand{\IuunitN}[0]{\textcolor{\complexcolor}{{\symCurrent}_{ \indexUnit, \symPhaseN}^{}  }}
\newcommand{\IuunitG}[0]{\textcolor{\complexcolor}{{\symCurrent}_{ \indexUnit, \symPhaseG}^{}  }}

\newcommand{\IuunitArated}[0]{\textcolor{\paramcolor}{{\symCurrent}_{ \indexUnit, \symPhaseA}^{\ratedss}  }}
\newcommand{\IuunitBrated}[0]{\textcolor{\paramcolor}{{\symCurrent}_{ \indexUnit, \symPhaseB}^{\ratedss}  }}
\newcommand{\IuunitCrated}[0]{\textcolor{\paramcolor}{{\symCurrent}_{ \indexUnit, \symPhaseC}^{\ratedss}  }}
\newcommand{\IuunitNrated}[0]{\textcolor{\paramcolor}{{\symCurrent}_{ \indexUnit, \symPhaseN}^{\ratedss}  }}

\newcommand{\IuunitABrated}[0]{\textcolor{\paramcolor}{{\symCurrent}_{ \indexUnit, \symPhaseA\symPhaseB}^{\ratedss}  }}
\newcommand{\IuunitBCrated}[0]{\textcolor{\paramcolor}{{\symCurrent}_{ \indexUnit, \symPhaseB\symPhaseC}^{\ratedss}  }}
\newcommand{\IuunitCArated}[0]{\textcolor{\paramcolor}{{\symCurrent}_{ \indexUnit, \symPhaseC\symPhaseA}^{\ratedss}  }}

\newcommand{\IbA}[0]{\textcolor{\complexcolor}{{\symCurrent}_{ \indexShunt, \symPhaseA}^{}  }}
\newcommand{\IbB}[0]{\textcolor{\complexcolor}{{\symCurrent}_{ \indexShunt, \symPhaseB}^{}  }}
\newcommand{\IbC}[0]{\textcolor{\complexcolor}{{\symCurrent}_{ \indexShunt, \symPhaseC}^{}  }}
\newcommand{\IbN}[0]{\textcolor{\complexcolor}{{\symCurrent}_{ \indexShunt, \symPhaseN}^{}  }}

\newcommand{\IuunitNreal}[0]{\textcolor{black}{{\symCurrent}_{ \indexUnit, \symPhaseN}^{\text{re}}  }}
\newcommand{\IuunitNimag}[0]{\textcolor{black}{{\symCurrent}_{ \indexUnit, \symPhaseN}^{\text{im}}  }}

\newcommand{\IuunitAB}[0]{\textcolor{\complexcolor}{{\symCurrent}_{ \indexUnit, \symPhaseA\symPhaseB}^{\Delta}  }}
\newcommand{\IuunitBC}[0]{\textcolor{\complexcolor}{{\symCurrent}_{ \indexUnit, \symPhaseB\symPhaseC}^{\Delta}  }}
\newcommand{\IuunitCA}[0]{\textcolor{\complexcolor}{{\symCurrent}_{ \indexUnit, \symPhaseC\symPhaseA}^{\Delta}  }}

\newcommand{\IuunitBA}[0]{\textcolor{\complexcolor}{{\symCurrent}_{ \indexUnit,\symPhaseB \symPhaseA}^{\Delta}  }}
\newcommand{\IuunitCB}[0]{\textcolor{\complexcolor}{{\symCurrent}_{ \indexUnit, \symPhaseC\symPhaseB}^{\Delta}  }}
\newcommand{\IuunitAC}[0]{\textcolor{\complexcolor}{{\symCurrent}_{ \indexUnit, \symPhaseA\symPhaseC}^{\Delta}  }}

\newcommand{\SuunitAAdelta}[0]{\textcolor{\complexcolor}{{\symApparentPower}_{ \indexUnit, \symPhaseA\symPhaseA}^{\Delta}  }}
\newcommand{\SuunitBBdelta}[0]{\textcolor{\complexcolor}{{\symApparentPower}_{ \indexUnit, \symPhaseB\symPhaseB}^{\Delta}  }}
\newcommand{\SuunitCCdelta}[0]{\textcolor{\complexcolor}{{\symApparentPower}_{ \indexUnit, \symPhaseC\symPhaseC}^{\Delta}  }}

\newcommand{\PuunitAAdelta}[0]{\textcolor{black}{{\symPower}_{ \indexUnit, \symPhaseA\symPhaseA}^{\Delta}  }}
\newcommand{\PuunitBBdelta}[0]{\textcolor{black}{{\symPower}_{ \indexUnit, \symPhaseB\symPhaseB}^{\Delta}  }}
\newcommand{\PuunitCCdelta}[0]{\textcolor{black}{{\symPower}_{ \indexUnit, \symPhaseC\symPhaseC}^{\Delta}  }}

\newcommand{\QuunitAAdelta}[0]{\textcolor{black}{{\symReactivePower}_{ \indexUnit, \symPhaseA\symPhaseA}^{\Delta}  }}
\newcommand{\QuunitBBdelta}[0]{\textcolor{black}{{\symReactivePower}_{ \indexUnit, \symPhaseB\symPhaseB}^{\Delta}  }}
\newcommand{\QuunitCCdelta}[0]{\textcolor{black}{{\symReactivePower}_{ \indexUnit, \symPhaseC\symPhaseC}^{\Delta}  }}

\newcommand{\SuunitABdelta}[0]{\textcolor{\complexcolor}{{\symApparentPower}_{ \indexUnit, \symPhaseA\symPhaseB}^{\Delta}  }}
\newcommand{\SuunitBCdelta}[0]{\textcolor{\complexcolor}{{\symApparentPower}_{ \indexUnit, \symPhaseB\symPhaseC}^{\Delta}  }}
\newcommand{\SuunitCAdelta}[0]{\textcolor{\complexcolor}{{\symApparentPower}_{ \indexUnit, \symPhaseC\symPhaseA}^{\Delta}  }}

\newcommand{\SuunitBAdelta}[0]{\textcolor{\complexcolor}{{\symApparentPower}_{ \indexUnit,\symPhaseB \symPhaseA}^{\Delta}  }}
\newcommand{\SuunitCBdelta}[0]{\textcolor{\complexcolor}{{\symApparentPower}_{ \indexUnit, \symPhaseC\symPhaseB}^{\Delta}  }}
\newcommand{\SuunitACdelta}[0]{\textcolor{\complexcolor}{{\symApparentPower}_{ \indexUnit, \symPhaseA\symPhaseC}^{\Delta}  }}

\newcommand{\SuunitABCdelta}[0]{\textcolor{\complexcolor}{{\symApparentPower}_{ \indexUnit, \symPhaseA\symPhaseB\symPhaseC}^{\Delta}  }}
\newcommand{\SuunitBCAdelta}[0]{\textcolor{\complexcolor}{{\symApparentPower}_{ \indexUnit, \symPhaseB\symPhaseC\symPhaseA}^{\Delta}  }}
\newcommand{\SuunitCABdelta}[0]{\textcolor{\complexcolor}{{\symApparentPower}_{ \indexUnit, \symPhaseC\symPhaseA\symPhaseB}^{\Delta}  }}

\newcommand{\PuunitABdeltamax}[0]{\textcolor{\paramcolor}{{\symPower}_{ \indexUnit, \symPhaseA\symPhaseB}^{\Delta, \maxss}  }}
\newcommand{\PuunitBCdeltamax}[0]{\textcolor{\paramcolor}{{\symPower}_{ \indexUnit, \symPhaseB\symPhaseC}^{\Delta, \maxss}  }}
\newcommand{\PuunitCAdeltamax}[0]{\textcolor{\paramcolor}{{\symPower}_{ \indexUnit, \symPhaseC\symPhaseA}^{\Delta, \maxss}  }}

\newcommand{\QuunitABdeltamax}[0]{\textcolor{\paramcolor}{{\symReactivePower}_{ \indexUnit, \symPhaseA\symPhaseB}^{\Delta, \maxss}  }}
\newcommand{\QuunitBCdeltamax}[0]{\textcolor{\paramcolor}{{\symReactivePower}_{ \indexUnit, \symPhaseB\symPhaseC}^{\Delta, \maxss}  }}
\newcommand{\QuunitCAdeltamax}[0]{\textcolor{\paramcolor}{{\symReactivePower}_{ \indexUnit, \symPhaseC\symPhaseA}^{\Delta, \maxss}  }}

\newcommand{\PuunitABdeltamin}[0]{\textcolor{\paramcolor}{{\symPower}_{ \indexUnit, \symPhaseA\symPhaseB}^{\Delta, \minss}  }}
\newcommand{\PuunitBCdeltamin}[0]{\textcolor{\paramcolor}{{\symPower}_{ \indexUnit, \symPhaseB\symPhaseC}^{\Delta, \minss}  }}
\newcommand{\PuunitCAdeltamin}[0]{\textcolor{\paramcolor}{{\symPower}_{ \indexUnit, \symPhaseC\symPhaseA}^{\Delta, \minss}  }}

\newcommand{\QuunitABdeltamin}[0]{\textcolor{\paramcolor}{{\symReactivePower}_{ \indexUnit, \symPhaseA\symPhaseB}^{\Delta, \minss}  }}
\newcommand{\QuunitBCdeltamin}[0]{\textcolor{\paramcolor}{{\symReactivePower}_{ \indexUnit, \symPhaseB\symPhaseC}^{\Delta, \minss}  }}
\newcommand{\QuunitCAdeltamin}[0]{\textcolor{\paramcolor}{{\symReactivePower}_{ \indexUnit, \symPhaseC\symPhaseA}^{\Delta, \minss}  }}

\newcommand{\PuunitAAmax}[0]{\textcolor{\paramcolor}{{\symPower}_{ \indexUnit, \symPhaseA\symPhaseA}^{\maxss}  }}
\newcommand{\PuunitBBmax}[0]{\textcolor{\paramcolor}{{\symPower}_{ \indexUnit, \symPhaseB\symPhaseB}^{ \maxss}  }}
\newcommand{\PuunitCCmax}[0]{\textcolor{\paramcolor}{{\symPower}_{ \indexUnit, \symPhaseC\symPhaseC}^{ \maxss}  }}

\newcommand{\QuunitAAmax}[0]{\textcolor{\paramcolor}{{\symReactivePower}_{ \indexUnit, \symPhaseA\symPhaseA}^{\maxss}  }}
\newcommand{\QuunitBBmax}[0]{\textcolor{\paramcolor}{{\symReactivePower}_{ \indexUnit, \symPhaseB\symPhaseB}^{ \maxss}  }}
\newcommand{\QuunitCCmax}[0]{\textcolor{\paramcolor}{{\symReactivePower}_{ \indexUnit, \symPhaseC\symPhaseC}^{ \maxss}  }}

\newcommand{\PuunitAAmin}[0]{\textcolor{\paramcolor}{{\symPower}_{ \indexUnit, \symPhaseA\symPhaseA}^{\minss}  }}
\newcommand{\PuunitBBmin}[0]{\textcolor{\paramcolor}{{\symPower}_{ \indexUnit, \symPhaseB\symPhaseB}^{ \minss}  }}
\newcommand{\PuunitCCmin}[0]{\textcolor{\paramcolor}{{\symPower}_{ \indexUnit, \symPhaseC\symPhaseC}^{ \minss}  }}

\newcommand{\QuunitAAmin}[0]{\textcolor{\paramcolor}{{\symReactivePower}_{ \indexUnit, \symPhaseA\symPhaseA}^{\minss}  }}
\newcommand{\QuunitBBmin}[0]{\textcolor{\paramcolor}{{\symReactivePower}_{ \indexUnit, \symPhaseB\symPhaseB}^{ \minss}  }}
\newcommand{\QuunitCCmin}[0]{\textcolor{\paramcolor}{{\symReactivePower}_{ \indexUnit, \symPhaseC\symPhaseC}^{ \minss}  }}

\newcommand{\PuunitABdelta}[0]{\textcolor{black}{{\symPower}_{ \indexUnit, \symPhaseA\symPhaseB}^{\Delta}  }}
\newcommand{\PuunitBCdelta}[0]{\textcolor{black}{{\symPower}_{ \indexUnit, \symPhaseB\symPhaseC}^{\Delta}  }}
\newcommand{\PuunitCAdelta}[0]{\textcolor{black}{{\symPower}_{ \indexUnit, \symPhaseC\symPhaseA}^{\Delta}  }}

\newcommand{\PuunitBAdelta}[0]{\textcolor{black}{{\symPower}_{ \indexUnit,\symPhaseB \symPhaseA}^{\Delta}  }}
\newcommand{\PuunitCBdelta}[0]{\textcolor{black}{{\symPower}_{ \indexUnit, \symPhaseC\symPhaseB}^{\Delta}  }}
\newcommand{\PuunitACdelta}[0]{\textcolor{black}{{\symPower}_{ \indexUnit, \symPhaseA\symPhaseC}^{\Delta}  }}

\newcommand{\QuunitABdelta}[0]{\textcolor{black}{{\symReactivePower}_{ \indexUnit, \symPhaseA\symPhaseB}^{\Delta}  }}
\newcommand{\QuunitBCdelta}[0]{\textcolor{black}{{\symReactivePower}_{ \indexUnit, \symPhaseB\symPhaseC}^{\Delta}  }}
\newcommand{\QuunitCAdelta}[0]{\textcolor{black}{{\symReactivePower}_{ \indexUnit, \symPhaseC\symPhaseA}^{\Delta}  }}

\newcommand{\QuunitBAdelta}[0]{\textcolor{black}{{\symReactivePower}_{ \indexUnit,\symPhaseB \symPhaseA}^{\Delta}  }}
\newcommand{\QuunitCBdelta}[0]{\textcolor{black}{{\symReactivePower}_{ \indexUnit, \symPhaseC\symPhaseB}^{\Delta}  }}
\newcommand{\QuunitACdelta}[0]{\textcolor{black}{{\symReactivePower}_{ \indexUnit, \symPhaseA\symPhaseC}^{\Delta}  }}

\newcommand{\PuunitAN}[0]{\textcolor{black}{{\symPower}_{ \indexUnit, \symPhaseA}^{}  }}
\newcommand{\PuunitBN}[0]{\textcolor{black}{{\symPower}_{ \indexUnit, \symPhaseB}^{}  }}
\newcommand{\PuunitCN}[0]{\textcolor{black}{{\symPower}_{ \indexUnit, \symPhaseC}^{}  }}

\newcommand{\SuunitrefAA}[0]{\textcolor{\complexparamcolor}{{\symApparentPower}_{ \indexUnit, \symPhaseA\symPhaseA}^{\refss}  }}
\newcommand{\SuunitrefBB}[0]{\textcolor{\complexparamcolor}{{\symApparentPower}_{ \indexUnit, \symPhaseB\symPhaseB}^{\refss}  }}
\newcommand{\SuunitrefCC}[0]{\textcolor{\complexparamcolor}{{\symApparentPower}_{ \indexUnit, \symPhaseC\symPhaseC}^{\refss}  }}

\newcommand{\SuunitrefAB}[0]{\textcolor{\complexparamcolor}{{\symApparentPower}_{ \indexUnit, \symPhaseA\symPhaseB}^{\refss}  }}
\newcommand{\SuunitrefBC}[0]{\textcolor{\complexparamcolor}{{\symApparentPower}_{ \indexUnit, \symPhaseB\symPhaseC}^{\refss}  }}
\newcommand{\SuunitrefCA}[0]{\textcolor{\complexparamcolor}{{\symApparentPower}_{ \indexUnit, \symPhaseC\symPhaseA}^{\refss}  }}

\newcommand{\PuunitrefAB}[0]{\textcolor{\paramcolor}{{\symPower}_{ \indexUnit, \symPhaseA\symPhaseB}^{\refss}  }}
\newcommand{\PuunitrefBC}[0]{\textcolor{\paramcolor}{{\symPower}_{ \indexUnit, \symPhaseB\symPhaseC}^{\refss}  }}
\newcommand{\PuunitrefCA}[0]{\textcolor{\paramcolor}{{\symPower}_{ \indexUnit, \symPhaseC\symPhaseA}^{\refss}  }}

\newcommand{\QuunitrefAB}[0]{\textcolor{\paramcolor}{{\symReactivePower}_{ \indexUnit, \symPhaseA\symPhaseB}^{\refss}  }}
\newcommand{\QuunitrefBC}[0]{\textcolor{\paramcolor}{{\symReactivePower}_{ \indexUnit, \symPhaseB\symPhaseC}^{\refss}  }}
\newcommand{\QuunitrefCA}[0]{\textcolor{\paramcolor}{{\symReactivePower}_{ \indexUnit, \symPhaseC\symPhaseA}^{\refss}  }}

\newcommand{\PuunitrefAA}[0]{\textcolor{\paramcolor}{{\symPower}_{ \indexUnit, \symPhaseA\symPhaseA}^{\refss}  }}
\newcommand{\PuunitrefBB}[0]{\textcolor{\paramcolor}{{\symPower}_{ \indexUnit, \symPhaseB\symPhaseB}^{\refss}  }}
\newcommand{\PuunitrefCC}[0]{\textcolor{\paramcolor}{{\symPower}_{ \indexUnit, \symPhaseC\symPhaseC}^{\refss}  }}

\newcommand{\QuunitrefAA}[0]{\textcolor{\paramcolor}{{\symReactivePower}_{ \indexUnit, \symPhaseA\symPhaseA}^{\refss}  }}
\newcommand{\QuunitrefBB}[0]{\textcolor{\paramcolor}{{\symReactivePower}_{ \indexUnit, \symPhaseB\symPhaseB}^{\refss}  }}
\newcommand{\QuunitrefCC}[0]{\textcolor{\paramcolor}{{\symReactivePower}_{ \indexUnit, \symPhaseC\symPhaseC}^{\refss}  }}

\newcommand{\SuunitAdot}[0]{\textcolor{\complexcolor}{\mathbf{\symApparentPower}_{ \indexUnit}^{\symPhaseA\!\cdot}  }}
\newcommand{\SuunitBdot}[0]{\textcolor{\complexcolor}{\mathbf{\symApparentPower}_{ \indexUnit}^{\symPhaseB\!\cdot}  }}
\newcommand{\SuunitCdot}[0]{\textcolor{\complexcolor}{\mathbf{\symApparentPower}_{ \indexUnit}^{\symPhaseC\!\cdot}  }}
\newcommand{\SuunitNdot}[0]{\textcolor{\complexcolor}{\mathbf{\symApparentPower}_{ \indexUnit}^{\symPhaseN\!\cdot}  }}
\newcommand{\SuunitGdot}[0]{\textcolor{\complexcolor}{\mathbf{\symApparentPower}_{ \indexUnit}^{\symPhaseG\!\cdot}  }}

\newcommand{\SuunitdotN}[0]{\textcolor{\complexcolor}{\mathbf{\symApparentPower}_{ \indexUnit}^{\cdot\! \symPhaseN}  }}

\newcommand{\PuunitNdot}[0]{\textcolor{black}{\mathbf{\symPower}_{ \indexUnit}^{\symPhaseN\!\cdot}  }}
\newcommand{\QuunitNdot}[0]{\textcolor{black}{\mathbf{\symReactivePower}_{ \indexUnit}^{\symPhaseN\!\cdot}  }}

\newcommand{\XuunitAA}[0]{\textcolor{\complexcolor}{{X}_{ \indexUnit, \symPhaseA\symPhaseA}^{}  }}
\newcommand{\XuunitAB}[0]{\textcolor{\complexcolor}{{X}_{ \indexUnit, \symPhaseA\symPhaseB}^{}  }}
\newcommand{\XuunitAC}[0]{\textcolor{\complexcolor}{{X}_{ \indexUnit, \symPhaseA\symPhaseC}^{}  }}

\newcommand{\XuunitBA}[0]{\textcolor{\complexcolor}{{X}_{ \indexUnit, \symPhaseB\symPhaseA}^{}  }}
\newcommand{\XuunitBB}[0]{\textcolor{\complexcolor}{{X}_{ \indexUnit, \symPhaseB\symPhaseB}^{}  }}
\newcommand{\XuunitBC}[0]{\textcolor{\complexcolor}{{X}_{ \indexUnit, \symPhaseB\symPhaseC}^{}  }}

\newcommand{\XuunitCA}[0]{\textcolor{\complexcolor}{{X}_{ \indexUnit, \symPhaseC\symPhaseA}^{}  }}
\newcommand{\XuunitCB}[0]{\textcolor{\complexcolor}{{X}_{ \indexUnit, \symPhaseC\symPhaseB}^{}  }}
\newcommand{\XuunitCC}[0]{\textcolor{\complexcolor}{{X}_{ \indexUnit, \symPhaseC\symPhaseC}^{}  }}

\newcommand{\XuunitAAreal}[0]{\textcolor{black}{{X}_{ \indexUnit, \symPhaseA\symPhaseA}^{\realss}  }}
\newcommand{\XuunitABreal}[0]{\textcolor{black}{{X}_{ \indexUnit, \symPhaseA\symPhaseB}^{\realss}  }}
\newcommand{\XuunitACreal}[0]{\textcolor{black}{{X}_{ \indexUnit, \symPhaseA\symPhaseC}^{\realss}  }}

\newcommand{\XuunitBAreal}[0]{\textcolor{black}{{X}_{ \indexUnit, \symPhaseB\symPhaseA}^{\realss}  }}
\newcommand{\XuunitBBreal}[0]{\textcolor{black}{{X}_{ \indexUnit, \symPhaseB\symPhaseB}^{\realss}  }}
\newcommand{\XuunitBCreal}[0]{\textcolor{black}{{X}_{ \indexUnit, \symPhaseB\symPhaseC}^{\realss}  }}

\newcommand{\XuunitCAreal}[0]{\textcolor{black}{{X}_{ \indexUnit, \symPhaseC\symPhaseA}^{\realss}  }}
\newcommand{\XuunitCBreal}[0]{\textcolor{black}{{X}_{ \indexUnit, \symPhaseC\symPhaseB}^{\realss}  }}
\newcommand{\XuunitCCreal}[0]{\textcolor{black}{{X}_{ \indexUnit, \symPhaseC\symPhaseC}^{\realss}  }}

\newcommand{\XuunitAAimag}[0]{\textcolor{black}{{X}_{ \indexUnit, \symPhaseA\symPhaseA}^{\imagss}  }}
\newcommand{\XuunitABimag}[0]{\textcolor{black}{{X}_{ \indexUnit, \symPhaseA\symPhaseB}^{\imagss}  }}
\newcommand{\XuunitACimag}[0]{\textcolor{black}{{X}_{ \indexUnit, \symPhaseA\symPhaseC}^{\imagss}  }}

\newcommand{\XuunitBAimag}[0]{\textcolor{black}{{X}_{ \indexUnit, \symPhaseB\symPhaseA}^{\imagss}  }}
\newcommand{\XuunitBBimag}[0]{\textcolor{black}{{X}_{ \indexUnit, \symPhaseB\symPhaseB}^{\imagss}  }}
\newcommand{\XuunitBCimag}[0]{\textcolor{black}{{X}_{ \indexUnit, \symPhaseB\symPhaseC}^{\imagss}  }}

\newcommand{\XuunitCAimag}[0]{\textcolor{black}{{X}_{ \indexUnit, \symPhaseC\symPhaseA}^{\imagss}  }}
\newcommand{\XuunitCBimag}[0]{\textcolor{black}{{X}_{ \indexUnit, \symPhaseC\symPhaseB}^{\imagss}  }}
\newcommand{\XuunitCCimag}[0]{\textcolor{black}{{X}_{ \indexUnit, \symPhaseC\symPhaseC}^{\imagss}  }}

\newcommand{\SuunitAA}[0]{\textcolor{\complexcolor}{{\symApparentPower}_{ \indexUnit, \symPhaseA\symPhaseA}^{}  }}
\newcommand{\SuunitAB}[0]{\textcolor{\complexcolor}{{\symApparentPower}_{ \indexUnit, \symPhaseA\symPhaseB}^{}  }}
\newcommand{\SuunitAC}[0]{\textcolor{\complexcolor}{{\symApparentPower}_{ \indexUnit, \symPhaseA\symPhaseC}^{}  }}
\newcommand{\SuunitAN}[0]{\textcolor{\complexcolor}{{\symApparentPower}_{ \indexUnit, \symPhaseA\symPhaseN}^{}  }}

\newcommand{\SuunitBA}[0]{\textcolor{\complexcolor}{{\symApparentPower}_{ \indexUnit, \symPhaseB\symPhaseA}^{}  }}
\newcommand{\SuunitBB}[0]{\textcolor{\complexcolor}{{\symApparentPower}_{ \indexUnit, \symPhaseB\symPhaseB}^{}  }}
\newcommand{\SuunitBC}[0]{\textcolor{\complexcolor}{{\symApparentPower}_{ \indexUnit, \symPhaseB\symPhaseC}^{}  }}
\newcommand{\SuunitBN}[0]{\textcolor{\complexcolor}{{\symApparentPower}_{ \indexUnit, \symPhaseB\symPhaseN}^{}  }}

\newcommand{\SuunitCA}[0]{\textcolor{\complexcolor}{{\symApparentPower}_{ \indexUnit, \symPhaseC\symPhaseA}^{}  }}
\newcommand{\SuunitCB}[0]{\textcolor{\complexcolor}{{\symApparentPower}_{ \indexUnit, \symPhaseC\symPhaseB}^{}  }}
\newcommand{\SuunitCC}[0]{\textcolor{\complexcolor}{{\symApparentPower}_{ \indexUnit, \symPhaseC\symPhaseC}^{}  }}
\newcommand{\SuunitCN}[0]{\textcolor{\complexcolor}{{\symApparentPower}_{ \indexUnit, \symPhaseC\symPhaseN}^{}  }}
\newcommand{\SuunitNN}[0]{\textcolor{\complexcolor}{{\symApparentPower}_{ \indexUnit, \symPhaseN\symPhaseN}^{}  }}
\newcommand{\SuunitGG}[0]{\textcolor{\complexcolor}{{\symApparentPower}_{ \indexUnit, \symPhaseG\symPhaseG}^{}  }}

\newcommand{\SuunitNA}[0]{\textcolor{\complexcolor}{{\symApparentPower}_{ \indexUnit, \symPhaseN\symPhaseA}^{}  }}
\newcommand{\SuunitNB}[0]{\textcolor{\complexcolor}{{\symApparentPower}_{ \indexUnit, \symPhaseN\symPhaseB}^{}  }}
\newcommand{\SuunitNC}[0]{\textcolor{\complexcolor}{{\symApparentPower}_{ \indexUnit, \symPhaseN\symPhaseC}^{}  }}

\newcommand{\PuunitNN}[0]{\textcolor{black}{{\symPower}_{ \indexUnit, \symPhaseN\symPhaseN}^{}  }}
\newcommand{\PuunitNA}[0]{\textcolor{black}{{\symPower}_{ \indexUnit, \symPhaseN\symPhaseA}^{}  }}
\newcommand{\PuunitNB}[0]{\textcolor{black}{{\symPower}_{ \indexUnit, \symPhaseN\symPhaseB}^{}  }}
\newcommand{\PuunitNC}[0]{\textcolor{black}{{\symPower}_{ \indexUnit, \symPhaseN\symPhaseC}^{}  }}

\newcommand{\PuunitAA}[0]{\textcolor{black}{{\symPower}_{ \indexUnit, \symPhaseA\symPhaseA}^{}  }}
\newcommand{\PuunitAB}[0]{\textcolor{black}{{\symPower}_{ \indexUnit, \symPhaseA\symPhaseB}^{}  }}
\newcommand{\PuunitAC}[0]{\textcolor{black}{{\symPower}_{ \indexUnit, \symPhaseA\symPhaseC}^{}  }}

\newcommand{\PuunitBA}[0]{\textcolor{black}{{\symPower}_{ \indexUnit, \symPhaseB\symPhaseA}^{}  }}
\newcommand{\PuunitBB}[0]{\textcolor{black}{{\symPower}_{ \indexUnit, \symPhaseB\symPhaseB}^{}  }}
\newcommand{\PuunitBC}[0]{\textcolor{black}{{\symPower}_{ \indexUnit, \symPhaseB\symPhaseC}^{}  }}

\newcommand{\PuunitCA}[0]{\textcolor{black}{{\symPower}_{ \indexUnit, \symPhaseC\symPhaseA}^{}  }}
\newcommand{\PuunitCB}[0]{\textcolor{black}{{\symPower}_{ \indexUnit, \symPhaseC\symPhaseB}^{}  }}
\newcommand{\PuunitCC}[0]{\textcolor{black}{{\symPower}_{ \indexUnit, \symPhaseC\symPhaseC}^{}  }}

\newcommand{\QuunitNN}[0]{\textcolor{black}{{\symReactivePower}_{ \indexUnit, \symPhaseN\symPhaseN}^{}  }}
\newcommand{\QuunitNA}[0]{\textcolor{black}{{\symReactivePower}_{ \indexUnit, \symPhaseN\symPhaseA}^{}  }}
\newcommand{\QuunitNB}[0]{\textcolor{black}{{\symReactivePower}_{ \indexUnit, \symPhaseN\symPhaseB}^{}  }}
\newcommand{\QuunitNC}[0]{\textcolor{black}{{\symReactivePower}_{ \indexUnit, \symPhaseN\symPhaseC}^{}  }}

\newcommand{\QuunitAA}[0]{\textcolor{black}{{\symReactivePower}_{ \indexUnit, \symPhaseA\symPhaseA}^{}  }}
\newcommand{\QuunitAB}[0]{\textcolor{black}{{\symReactivePower}_{ \indexUnit, \symPhaseA\symPhaseB}^{}  }}
\newcommand{\QuunitAC}[0]{\textcolor{black}{{\symReactivePower}_{ \indexUnit, \symPhaseA\symPhaseC}^{}  }}

\newcommand{\QuunitBA}[0]{\textcolor{black}{{\symReactivePower}_{ \indexUnit, \symPhaseB\symPhaseA}^{}  }}
\newcommand{\QuunitBB}[0]{\textcolor{black}{{\symReactivePower}_{ \indexUnit, \symPhaseB\symPhaseB}^{}  }}
\newcommand{\QuunitBC}[0]{\textcolor{black}{{\symReactivePower}_{ \indexUnit, \symPhaseB\symPhaseC}^{}  }}

\newcommand{\QuunitCA}[0]{\textcolor{black}{{\symReactivePower}_{ \indexUnit, \symPhaseC\symPhaseA}^{}  }}
\newcommand{\QuunitCB}[0]{\textcolor{black}{{\symReactivePower}_{ \indexUnit, \symPhaseC\symPhaseB}^{}  }}
\newcommand{\QuunitCC}[0]{\textcolor{black}{{\symReactivePower}_{ \indexUnit, \symPhaseC\symPhaseC}^{}  }}

\newcommand{\QuunitAN}[0]{\textcolor{black}{{\symReactivePower}_{ \indexUnit, \symPhaseA}^{}  }}
\newcommand{\QuunitBN}[0]{\textcolor{black}{{\symReactivePower}_{ \indexUnit, \symPhaseB}^{}  }}
\newcommand{\QuunitCN}[0]{\textcolor{black}{{\symReactivePower}_{ \indexUnit, \symPhaseC}^{}  }}

\newcommand{\SuunitrefSDP}[0]{\textcolor{\complexparamcolor}{\mathbf{\symApparentPower}_{ \indexUnit}^{\refss}  }}
\newcommand{\PuunitrefSDP}[0]{\textcolor{\paramcolor}{\mathbf{\symPower}_{ \indexUnit}^{\refss}  }}
\newcommand{\QuunitrefSDP}[0]{\textcolor{\paramcolor}{\mathbf{\symReactivePower}_{ \indexUnit}^{\refss}  }}

\newcommand{\PuunitrefAN}[0]{\textcolor{\paramcolor}{{\symPower}_{ \indexUnit, \symPhaseA}^{\refss}  }}
\newcommand{\PuunitrefBN}[0]{\textcolor{\paramcolor}{{\symPower}_{ \indexUnit, \symPhaseB}^{\refss}  }}
\newcommand{\PuunitrefCN}[0]{\textcolor{\paramcolor}{{\symPower}_{ \indexUnit, \symPhaseC}^{\refss}  }}

\newcommand{\QuunitrefAN}[0]{\textcolor{\paramcolor}{{\symReactivePower}_{ \indexUnit, \symPhaseA}^{\refss}  }}
\newcommand{\QuunitrefBN}[0]{\textcolor{\paramcolor}{{\symReactivePower}_{ \indexUnit, \symPhaseB}^{\refss}  }}
\newcommand{\QuunitrefCN}[0]{\textcolor{\paramcolor}{{\symReactivePower}_{ \indexUnit, \symPhaseC}^{\refss}  }}

\newcommand{\SbNN}[0]{\textcolor{\complexcolor}{{\symApparentPower}_{ \indexShunt, \symPhaseN\symPhaseN}^{}  }}

\newcommand{\PijkrefAA}[0]{\textcolor{\paramcolor}{\symPower_{\indexGridLines\indexGridNode \indexGridNodeTwo, \symPhaseA\symPhaseA}^{\refss}  }}
\newcommand{\PijkrefBB}[0]{\textcolor{\paramcolor}{\symPower_{\indexGridLines\indexGridNode \indexGridNodeTwo, \symPhaseB\symPhaseB}^{\refss}  }}
\newcommand{\PijkrefCC}[0]{\textcolor{\paramcolor}{\symPower_{\indexGridLines\indexGridNode \indexGridNodeTwo, \symPhaseC\symPhaseC}^{\refss}  }}

\newcommand{\Pijkpp}[0]{\textcolor{black}{\symPower_{\indexGridLines\indexGridNode \indexGridNodeTwo, \indexPhases\indexPhases}^{}  }}
\newcommand{\Qijkpp}[0]{\textcolor{black}{\symReactivePower_{\indexGridLines\indexGridNode \indexGridNodeTwo, \indexPhases\indexPhases}^{}  }}

\newcommand{\PjikAA}[0]{\textcolor{black}{\symPower_{\indexGridLines\indexGridNodeTwo\indexGridNode , \symPhaseA\symPhaseA}^{}  }}
\newcommand{\PijkAAloss}[0]{\textcolor{black}{\symPower_{\indexGridLines\indexGridNode \indexGridNodeTwo, \symPhaseA\symPhaseA}^{\lossss}  }}
\newcommand{\PjikAAloss}[0]{\textcolor{black}{\symPower_{\indexGridLines\indexGridNodeTwo\indexGridNode , \symPhaseA\symPhaseA}^{\lossss}  }}

\newcommand{\SijkAA}[0]{\textcolor{\complexcolor}{\symApparentPower_{\indexGridLines\indexGridNode \indexGridNodeTwo, \symPhaseA\symPhaseA}^{}  }}
\newcommand{\SijkAB}[0]{\textcolor{\complexcolor}{\symApparentPower_{\indexGridLines\indexGridNode \indexGridNodeTwo, \symPhaseA\symPhaseB}^{}  }}
\newcommand{\SijkAC}[0]{\textcolor{\complexcolor}{\symApparentPower_{\indexGridLines\indexGridNode \indexGridNodeTwo, \symPhaseA\symPhaseC}^{}  }}
\newcommand{\SijkAN}[0]{\textcolor{\complexcolor}{\symApparentPower_{\indexGridLines\indexGridNode \indexGridNodeTwo, \symPhaseA\symPhaseN}^{}  }}
\newcommand{\SijkAG}[0]{\textcolor{\complexcolor}{\symApparentPower_{\indexGridLines\indexGridNode \indexGridNodeTwo, \symPhaseA\symPhaseG}^{}  }}
\newcommand{\SijkBA}[0]{\textcolor{\complexcolor}{\symApparentPower_{\indexGridLines\indexGridNode \indexGridNodeTwo, \symPhaseB\symPhaseA}^{}  }}
\newcommand{\SijkBB}[0]{\textcolor{\complexcolor}{\symApparentPower_{\indexGridLines\indexGridNode \indexGridNodeTwo, \symPhaseB\symPhaseB}^{}  }}
\newcommand{\SijkBC}[0]{\textcolor{\complexcolor}{\symApparentPower_{\indexGridLines\indexGridNode \indexGridNodeTwo, \symPhaseB\symPhaseC}^{}  }}
\newcommand{\SijkBN}[0]{\textcolor{\complexcolor}{\symApparentPower_{\indexGridLines\indexGridNode \indexGridNodeTwo, \symPhaseB\symPhaseN}^{}  }}
\newcommand{\SijkBG}[0]{\textcolor{\complexcolor}{\symApparentPower_{\indexGridLines\indexGridNode \indexGridNodeTwo, \symPhaseB\symPhaseG}^{}  }}
\newcommand{\SijkCA}[0]{\textcolor{\complexcolor}{\symApparentPower_{\indexGridLines\indexGridNode \indexGridNodeTwo, \symPhaseC\symPhaseA}^{}  }}
\newcommand{\SijkCB}[0]{\textcolor{\complexcolor}{\symApparentPower_{\indexGridLines\indexGridNode \indexGridNodeTwo, \symPhaseC\symPhaseB}^{}  }}
\newcommand{\SijkCC}[0]{\textcolor{\complexcolor}{\symApparentPower_{\indexGridLines\indexGridNode \indexGridNodeTwo, \symPhaseC\symPhaseC}^{}  }}
\newcommand{\SijkCN}[0]{\textcolor{\complexcolor}{\symApparentPower_{\indexGridLines\indexGridNode \indexGridNodeTwo, \symPhaseC\symPhaseN}^{}  }}
\newcommand{\SijkCG}[0]{\textcolor{\complexcolor}{\symApparentPower_{\indexGridLines\indexGridNode \indexGridNodeTwo, \symPhaseC\symPhaseG}^{}  }}

\newcommand{\SijkNA}[0]{\textcolor{\complexcolor}{\symApparentPower_{\indexGridLines\indexGridNode \indexGridNodeTwo, \symPhaseN\symPhaseA}^{}  }}
\newcommand{\SijkNB}[0]{\textcolor{\complexcolor}{\symApparentPower_{\indexGridLines\indexGridNode \indexGridNodeTwo, \symPhaseN\symPhaseB}^{}  }}
\newcommand{\SijkNC}[0]{\textcolor{\complexcolor}{\symApparentPower_{\indexGridLines\indexGridNode \indexGridNodeTwo, \symPhaseN\symPhaseC}^{}  }}
\newcommand{\SijkNN}[0]{\textcolor{\complexcolor}{\symApparentPower_{\indexGridLines\indexGridNode \indexGridNodeTwo, \symPhaseN\symPhaseN}^{}  }}
\newcommand{\SijkNG}[0]{\textcolor{\complexcolor}{\symApparentPower_{\indexGridLines\indexGridNode \indexGridNodeTwo, \symPhaseN\symPhaseG}^{}  }}
\newcommand{\PijkNG}[0]{\textcolor{black}{\symPower_{\indexGridLines\indexGridNode \indexGridNodeTwo, \symPhaseN\symPhaseG}^{}  }}
\newcommand{\QijkNG}[0]{\textcolor{black}{\symReactivePower_{\indexGridLines\indexGridNode \indexGridNodeTwo, \symPhaseN\symPhaseG}^{}  }}

\newcommand{\SijkGA}[0]{\textcolor{\complexcolor}{\symApparentPower_{\indexGridLines\indexGridNode \indexGridNodeTwo, \symPhaseG\symPhaseA}^{}  }}
\newcommand{\SijkGB}[0]{\textcolor{\complexcolor}{\symApparentPower_{\indexGridLines\indexGridNode \indexGridNodeTwo, \symPhaseG\symPhaseB}^{}  }}
\newcommand{\SijkGC}[0]{\textcolor{\complexcolor}{\symApparentPower_{\indexGridLines\indexGridNode \indexGridNodeTwo, \symPhaseG\symPhaseC}^{}  }}
\newcommand{\SijkGN}[0]{\textcolor{\complexcolor}{\symApparentPower_{\indexGridLines\indexGridNode \indexGridNodeTwo, \symPhaseG\symPhaseN}^{}  }}
\newcommand{\SijkGG}[0]{\textcolor{\complexcolor}{\symApparentPower_{\indexGridLines\indexGridNode \indexGridNodeTwo, \symPhaseG\symPhaseG}^{}  }}

\newcommand{\SjikAA}[0]{\textcolor{\complexcolor}{\symApparentPower_{\indexGridLines\indexGridNodeTwo\indexGridNode , \symPhaseA\symPhaseA}^{}  }}
\newcommand{\SjikAB}[0]{\textcolor{\complexcolor}{\symApparentPower_{\indexGridLines\indexGridNodeTwo\indexGridNode , \symPhaseA\symPhaseB}^{}  }}
\newcommand{\SjikAC}[0]{\textcolor{\complexcolor}{\symApparentPower_{\indexGridLines\indexGridNodeTwo\indexGridNode , \symPhaseA\symPhaseC}^{}  }}
\newcommand{\SjikAN}[0]{\textcolor{\complexcolor}{\symApparentPower_{\indexGridLines\indexGridNodeTwo\indexGridNode , \symPhaseA\symPhaseN}^{}  }}
\newcommand{\SjikBB}[0]{\textcolor{\complexcolor}{\symApparentPower_{\indexGridLines\indexGridNodeTwo\indexGridNode , \symPhaseB\symPhaseB}^{}  }}
\newcommand{\SjikCC}[0]{\textcolor{\complexcolor}{\symApparentPower_{\indexGridLines\indexGridNodeTwo\indexGridNode , \symPhaseC\symPhaseC}^{}  }}
\newcommand{\SjikNA}[0]{\textcolor{\complexcolor}{\symApparentPower_{\indexGridLines\indexGridNodeTwo\indexGridNode , \symPhaseN\symPhaseA}^{}  }}
\newcommand{\SjikNB}[0]{\textcolor{\complexcolor}{\symApparentPower_{\indexGridLines\indexGridNodeTwo\indexGridNode , \symPhaseN\symPhaseB}^{}  }}
\newcommand{\SjikNC}[0]{\textcolor{\complexcolor}{\symApparentPower_{\indexGridLines\indexGridNodeTwo\indexGridNode , \symPhaseN\symPhaseC}^{}  }}
\newcommand{\SjikNN}[0]{\textcolor{\complexcolor}{\symApparentPower_{\indexGridLines\indexGridNodeTwo\indexGridNode , \symPhaseN\symPhaseN}^{}  }}

\newcommand{\SikNG}[0]{\textcolor{\complexcolor}{\symApparentPower_{\indexGridNode, \symPhaseN\symPhaseG}^{}  }}
\newcommand{\PikNG}[0]{\textcolor{black}{\symPower_{\indexGridNode, \symPhaseN\symPhaseG}^{}  }}
\newcommand{\QikNG}[0]{\textcolor{black}{\symReactivePower_{\indexGridNode, \symPhaseN\symPhaseG}^{}  }}

\newcommand{\PijkAA}[0]{\textcolor{black}{\symPower_{\indexGridLines\indexGridNode \indexGridNodeTwo, \symPhaseA\symPhaseA}^{}  }}
\newcommand{\PijkAB}[0]{\textcolor{black}{\symPower_{\indexGridLines\indexGridNode \indexGridNodeTwo, \symPhaseA\symPhaseB}^{}  }}
\newcommand{\PijkAC}[0]{\textcolor{black}{\symPower_{\indexGridLines\indexGridNode \indexGridNodeTwo, \symPhaseA\symPhaseC}^{}  }}
\newcommand{\PijkBA}[0]{\textcolor{black}{\symPower_{\indexGridLines\indexGridNode \indexGridNodeTwo, \symPhaseB\symPhaseA}^{}  }}
\newcommand{\PijkBB}[0]{\textcolor{black}{\symPower_{\indexGridLines\indexGridNode \indexGridNodeTwo, \symPhaseB\symPhaseB}^{}  }}
\newcommand{\PijkBC}[0]{\textcolor{black}{\symPower_{\indexGridLines\indexGridNode \indexGridNodeTwo, \symPhaseB\symPhaseC}^{}  }}
\newcommand{\PijkCA}[0]{\textcolor{black}{\symPower_{\indexGridLines\indexGridNode \indexGridNodeTwo, \symPhaseC\symPhaseA}^{}  }}
\newcommand{\PijkCB}[0]{\textcolor{black}{\symPower_{\indexGridLines\indexGridNode \indexGridNodeTwo, \symPhaseC\symPhaseB}^{}  }}
\newcommand{\PijkCC}[0]{\textcolor{black}{\symPower_{\indexGridLines\indexGridNode \indexGridNodeTwo, \symPhaseC\symPhaseC}^{}  }}
\newcommand{\PijkNN}[0]{\textcolor{black}{\symPower_{\indexGridLines\indexGridNode \indexGridNodeTwo, \symPhaseN\symPhaseN}^{}  }}
\newcommand{\QijkAA}[0]{\textcolor{black}{\symReactivePower_{\indexGridLines\indexGridNode \indexGridNodeTwo, \symPhaseA\symPhaseA}^{}  }}
\newcommand{\QijkAB}[0]{\textcolor{black}{\symReactivePower_{\indexGridLines\indexGridNode \indexGridNodeTwo, \symPhaseA\symPhaseB}^{}  }}
\newcommand{\QijkAC}[0]{\textcolor{black}{\symReactivePower_{\indexGridLines\indexGridNode \indexGridNodeTwo, \symPhaseA\symPhaseC}^{}  }}
\newcommand{\QijkBA}[0]{\textcolor{black}{\symReactivePower_{\indexGridLines\indexGridNode \indexGridNodeTwo, \symPhaseB\symPhaseA}^{}  }}
\newcommand{\QijkBB}[0]{\textcolor{black}{\symReactivePower_{\indexGridLines\indexGridNode \indexGridNodeTwo, \symPhaseB\symPhaseB}^{}  }}
\newcommand{\QijkBC}[0]{\textcolor{black}{\symReactivePower_{\indexGridLines\indexGridNode \indexGridNodeTwo, \symPhaseB\symPhaseC}^{}  }}
\newcommand{\QijkCA}[0]{\textcolor{black}{\symReactivePower_{\indexGridLines\indexGridNode \indexGridNodeTwo, \symPhaseC\symPhaseA}^{}  }}
\newcommand{\QijkCB}[0]{\textcolor{black}{\symReactivePower_{\indexGridLines\indexGridNode \indexGridNodeTwo, \symPhaseC\symPhaseB}^{}  }}
\newcommand{\QijkCC}[0]{\textcolor{black}{\symReactivePower_{\indexGridLines\indexGridNode \indexGridNodeTwo, \symPhaseC\symPhaseC}^{}  }}
\newcommand{\QijkNN}[0]{\textcolor{black}{\symReactivePower_{\indexGridLines\indexGridNode \indexGridNodeTwo, \symPhaseN\symPhaseN}^{}  }}

\newcommand{\PjikAB}[0]{\textcolor{black}{\symPower_{\indexGridLines \indexGridNodeTwo \indexGridNode, \symPhaseA\symPhaseB}^{}  }}
\newcommand{\QjikAB}[0]{\textcolor{black}{\symReactivePower_{\indexGridLines \indexGridNodeTwo \indexGridNode, \symPhaseA\symPhaseB}^{}  }}

\newcommand{\PijksAA}[0]{\textcolor{black}{\symPower_{\indexGridLines\indexGridNode \indexGridNodeTwo, \symPhaseA\symPhaseA }^{\seriesss}  }}
\newcommand{\QijksAA}[0]{\textcolor{black}{\symReactivePower_{\indexGridLines\indexGridNode \indexGridNodeTwo, \symPhaseA\symPhaseA }^{\seriesss}  }}
\newcommand{\SijksAA}[0]{\textcolor{\complexcolor}{\symApparentPower_{\indexGridLines\indexGridNode \indexGridNodeTwo, \symPhaseA\symPhaseA }^{\seriesss}  }}
\newcommand{\PijksAB}[0]{\textcolor{black}{\symPower_{\indexGridLines\indexGridNode \indexGridNodeTwo, \symPhaseA\symPhaseB }^{\seriesss}  }}
\newcommand{\QijksAB}[0]{\textcolor{black}{\symReactivePower_{\indexGridLines\indexGridNode \indexGridNodeTwo, \symPhaseA\symPhaseB }^{\seriesss}  }}
\newcommand{\SijksAB}[0]{\textcolor{\complexcolor}{\symApparentPower_{\indexGridLines\indexGridNode \indexGridNodeTwo, \symPhaseA\symPhaseB }^{\seriesss}  }}
\newcommand{\PijksAC}[0]{\textcolor{black}{\symPower_{\indexGridLines\indexGridNode \indexGridNodeTwo, \symPhaseA\symPhaseC }^{\seriesss}  }}
\newcommand{\QijksAC}[0]{\textcolor{black}{\symReactivePower_{\indexGridLines\indexGridNode \indexGridNodeTwo, \symPhaseA\symPhaseC }^{\seriesss}  }}
\newcommand{\SijksAC}[0]{\textcolor{\complexcolor}{\symApparentPower_{\indexGridLines\indexGridNode \indexGridNodeTwo, \symPhaseA\symPhaseC }^{\seriesss}  }}
\newcommand{\PijksBA}[0]{\textcolor{black}{\symPower_{\indexGridLines\indexGridNode \indexGridNodeTwo, \symPhaseB\symPhaseA}^{\seriesss}  }}
\newcommand{\QijksBA}[0]{\textcolor{black}{\symReactivePower_{\indexGridLines\indexGridNode \indexGridNodeTwo, \symPhaseB\symPhaseA}^{\seriesss}  }}
\newcommand{\SijksBA}[0]{\textcolor{\complexcolor}{\symApparentPower_{\indexGridLines\indexGridNode \indexGridNodeTwo, \symPhaseB\symPhaseA}^{\seriesss}  }}
\newcommand{\PijksBB}[0]{\textcolor{black}{\symPower_{\indexGridLines\indexGridNode \indexGridNodeTwo, \symPhaseB\symPhaseB}^{\seriesss}  }}
\newcommand{\QijksBB}[0]{\textcolor{black}{\symReactivePower_{\indexGridLines\indexGridNode \indexGridNodeTwo, \symPhaseB\symPhaseB}^{\seriesss}  }}
\newcommand{\SijksBB}[0]{\textcolor{\complexcolor}{\symApparentPower_{\indexGridLines\indexGridNode \indexGridNodeTwo, \symPhaseB\symPhaseB}^{\seriesss}  }}
\newcommand{\PijksBC}[0]{\textcolor{black}{\symPower_{\indexGridLines\indexGridNode \indexGridNodeTwo, \symPhaseB\symPhaseC}^{\seriesss}  }}
\newcommand{\QijksBC}[0]{\textcolor{black}{\symReactivePower_{\indexGridLines\indexGridNode \indexGridNodeTwo, \symPhaseB\symPhaseC}^{\seriesss}  }}
\newcommand{\SijksBC}[0]{\textcolor{\complexcolor}{\symApparentPower_{\indexGridLines\indexGridNode \indexGridNodeTwo, \symPhaseB\symPhaseC}^{\seriesss}  }}
\newcommand{\PijksCA}[0]{\textcolor{black}{\symPower_{\indexGridLines\indexGridNode \indexGridNodeTwo, \symPhaseC\symPhaseA}^{\seriesss}  }}
\newcommand{\QijksCA}[0]{\textcolor{black}{\symReactivePower_{\indexGridLines\indexGridNode \indexGridNodeTwo, \symPhaseC\symPhaseA}^{\seriesss}  }}
\newcommand{\SijksCA}[0]{\textcolor{\complexcolor}{\symApparentPower_{\indexGridLines\indexGridNode \indexGridNodeTwo, \symPhaseC\symPhaseA}^{\seriesss}  }}
\newcommand{\PijksCB}[0]{\textcolor{black}{\symPower_{\indexGridLines\indexGridNode \indexGridNodeTwo, \symPhaseC\symPhaseB}^{\seriesss}  }}
\newcommand{\QijksCB}[0]{\textcolor{black}{\symReactivePower_{\indexGridLines\indexGridNode \indexGridNodeTwo, \symPhaseC\symPhaseB}^{\seriesss}  }}
\newcommand{\SijksCB}[0]{\textcolor{\complexcolor}{\symApparentPower_{\indexGridLines\indexGridNode \indexGridNodeTwo, \symPhaseC\symPhaseB}^{\seriesss}  }}
\newcommand{\PijksCC}[0]{\textcolor{black}{\symPower_{\indexGridLines\indexGridNode \indexGridNodeTwo, \symPhaseC\symPhaseC}^{\seriesss}  }}
\newcommand{\QijksCC}[0]{\textcolor{black}{\symReactivePower_{\indexGridLines\indexGridNode \indexGridNodeTwo, \symPhaseC\symPhaseC}^{\seriesss}  }}
\newcommand{\SijksCC}[0]{\textcolor{\complexcolor}{\symApparentPower_{\indexGridLines\indexGridNode \indexGridNodeTwo, \symPhaseC\symPhaseC}^{\seriesss}  }}

\newcommand{\PijksAN}[0]{\textcolor{black}{\symPower_{\indexGridLines\indexGridNode \indexGridNodeTwo, \symPhaseA\symPhaseN }^{\seriesss}  }}
\newcommand{\QijksAN}[0]{\textcolor{black}{\symReactivePower_{\indexGridLines\indexGridNode \indexGridNodeTwo, \symPhaseA\symPhaseN }^{\seriesss}  }}
\newcommand{\SijksAN}[0]{\textcolor{\complexcolor}{\symApparentPower_{\indexGridLines\indexGridNode \indexGridNodeTwo, \symPhaseA\symPhaseN }^{\seriesss}  }}
\newcommand{\PijksBN}[0]{\textcolor{black}{\symPower_{\indexGridLines\indexGridNode \indexGridNodeTwo, \symPhaseB\symPhaseN }^{\seriesss}  }}
\newcommand{\QijksBN}[0]{\textcolor{black}{\symReactivePower_{\indexGridLines\indexGridNode \indexGridNodeTwo, \symPhaseB\symPhaseN }^{\seriesss}  }}
\newcommand{\SijksBN}[0]{\textcolor{\complexcolor}{\symApparentPower_{\indexGridLines\indexGridNode \indexGridNodeTwo, \symPhaseB\symPhaseN }^{\seriesss}  }}
\newcommand{\PijksCN}[0]{\textcolor{black}{\symPower_{\indexGridLines\indexGridNode \indexGridNodeTwo, \symPhaseC\symPhaseN }^{\seriesss}  }}
\newcommand{\QijksCN}[0]{\textcolor{black}{\symReactivePower_{\indexGridLines\indexGridNode \indexGridNodeTwo, \symPhaseC\symPhaseN }^{\seriesss}  }}
\newcommand{\SijksCN}[0]{\textcolor{\complexcolor}{\symApparentPower_{\indexGridLines\indexGridNode \indexGridNodeTwo, \symPhaseC\symPhaseN }^{\seriesss}  }}
\newcommand{\PijksNN}[0]{\textcolor{black}{\symPower_{\indexGridLines\indexGridNode \indexGridNodeTwo, \symPhaseN\symPhaseN }^{\seriesss}  }}
\newcommand{\QijksNN}[0]{\textcolor{black}{\symReactivePower_{\indexGridLines\indexGridNode \indexGridNodeTwo, \symPhaseN\symPhaseN }^{\seriesss}  }}
\newcommand{\SijksNN}[0]{\textcolor{\complexcolor}{\symApparentPower_{\indexGridLines\indexGridNode \indexGridNodeTwo, \symPhaseN\symPhaseN }^{\seriesss}  }}
\newcommand{\PijksNA}[0]{\textcolor{black}{\symPower_{\indexGridLines\indexGridNode \indexGridNodeTwo, \symPhaseN\symPhaseA }^{\seriesss}  }}
\newcommand{\QijksNA}[0]{\textcolor{black}{\symReactivePower_{\indexGridLines\indexGridNode \indexGridNodeTwo, \symPhaseN\symPhaseA }^{\seriesss}  }}
\newcommand{\SijksNA}[0]{\textcolor{\complexcolor}{\symApparentPower_{\indexGridLines\indexGridNode \indexGridNodeTwo, \symPhaseN\symPhaseA }^{\seriesss}  }}
\newcommand{\PijksNB}[0]{\textcolor{black}{\symPower_{\indexGridLines\indexGridNode \indexGridNodeTwo, \symPhaseN\symPhaseB }^{\seriesss}  }}
\newcommand{\QijksNB}[0]{\textcolor{black}{\symReactivePower_{\indexGridLines\indexGridNode \indexGridNodeTwo, \symPhaseN\symPhaseB }^{\seriesss}  }}
\newcommand{\SijksNB}[0]{\textcolor{\complexcolor}{\symApparentPower_{\indexGridLines\indexGridNode \indexGridNodeTwo, \symPhaseN\symPhaseB }^{\seriesss}  }}
\newcommand{\PijksNC}[0]{\textcolor{black}{\symPower_{\indexGridLines\indexGridNode \indexGridNodeTwo, \symPhaseN\symPhaseC }^{\seriesss}  }}
\newcommand{\QijksNC}[0]{\textcolor{black}{\symReactivePower_{\indexGridLines\indexGridNode \indexGridNodeTwo, \symPhaseN\symPhaseC }^{\seriesss}  }}
\newcommand{\SijksNC}[0]{\textcolor{\complexcolor}{\symApparentPower_{\indexGridLines\indexGridNode \indexGridNodeTwo, \symPhaseN\symPhaseC }^{\seriesss}  }}

\newcommand{\PijklossAA}[0]{\textcolor{black}{\symPower_{\indexGridNode \indexGridNodeTwo, \symPhaseA\symPhaseA}^{\lossss}  }}
\newcommand{\PijklossAB}[0]{\textcolor{black}{\symPower_{\indexGridNode \indexGridNodeTwo, \symPhaseA\symPhaseB}^{\lossss}  }}
\newcommand{\PijklossAC}[0]{\textcolor{black}{\symPower_{\indexGridNode \indexGridNodeTwo, \symPhaseA\symPhaseC}^{\lossss}  }}
\newcommand{\PijklossBA}[0]{\textcolor{black}{\symPower_{\indexGridNode \indexGridNodeTwo, \symPhaseB\symPhaseA}^{\lossss}  }}
\newcommand{\PijklossBB}[0]{\textcolor{black}{\symPower_{\indexGridNode \indexGridNodeTwo, \symPhaseB\symPhaseB}^{\lossss}  }}
\newcommand{\PijklossBC}[0]{\textcolor{black}{\symPower_{\indexGridNode \indexGridNodeTwo, \symPhaseB\symPhaseC}^{\lossss}  }}
\newcommand{\PijklossCA}[0]{\textcolor{black}{\symPower_{\indexGridNode \indexGridNodeTwo, \symPhaseC\symPhaseA}^{\lossss}  }}
\newcommand{\PijklossCB}[0]{\textcolor{black}{\symPower_{\indexGridNode \indexGridNodeTwo, \symPhaseC\symPhaseB}^{\lossss}  }}
\newcommand{\PijklossCC}[0]{\textcolor{black}{\symPower_{\indexGridNode \indexGridNodeTwo, \symPhaseC\symPhaseC}^{\lossss}  }}
\newcommand{\QijklossAA}[0]{\textcolor{black}{\symReactivePower_{\indexGridNode \indexGridNodeTwo, \symPhaseA\symPhaseA}^{\lossss}  }}
\newcommand{\QijklossAB}[0]{\textcolor{black}{\symReactivePower_{\indexGridNode \indexGridNodeTwo, \symPhaseA\symPhaseB}^{\lossss}  }}
\newcommand{\QijklossAC}[0]{\textcolor{black}{\symReactivePower_{\indexGridNode \indexGridNodeTwo, \symPhaseA\symPhaseC}^{\lossss}  }}
\newcommand{\QijklossBA}[0]{\textcolor{black}{\symReactivePower_{\indexGridNode \indexGridNodeTwo, \symPhaseB\symPhaseA}^{\lossss}  }}
\newcommand{\QijklossBB}[0]{\textcolor{black}{\symReactivePower_{\indexGridNode \indexGridNodeTwo, \symPhaseB\symPhaseB}^{\lossss}  }}
\newcommand{\QijklossBC}[0]{\textcolor{black}{\symReactivePower_{\indexGridNode \indexGridNodeTwo, \symPhaseB\symPhaseC}^{\lossss}  }}
\newcommand{\QijklossCA}[0]{\textcolor{black}{\symReactivePower_{\indexGridNode \indexGridNodeTwo, \symPhaseC\symPhaseA}^{\lossss}  }}
\newcommand{\QijklossCB}[0]{\textcolor{black}{\symReactivePower_{\indexGridNode \indexGridNodeTwo, \symPhaseC\symPhaseB}^{\lossss}  }}
\newcommand{\QijklossCC}[0]{\textcolor{black}{\symReactivePower_{\indexGridNode \indexGridNodeTwo, \symPhaseC\symPhaseC}^{\lossss}  }}


\newcommand{\Iik}[0]{\textcolor{\complexcolor}{\symCurrent_{\indexGridNode }^{}  }}

\newcommand{\Ijik}[0]{\textcolor{\complexcolor}{\symCurrent_{\indexGridLines\indexGridNodeTwo  \indexGridNode }^{}  }}

\newcommand{\Iijkre}[0]{\textcolor{black}{\symCurrent_{\indexGridLines\indexGridNode \indexGridNodeTwo}^{\realss}  }}
\newcommand{\Iijkim}[0]{\textcolor{black}{\symCurrent_{\indexGridLines\indexGridNode \indexGridNodeTwo}^{\imagss}  }}
\newcommand{\Iijkangle}[0]{\textcolor{black}{\symCurrent_{\indexGridLines\indexGridNode \indexGridNodeTwo}^{\imagangle}  }}
\newcommand{\Iijkanglemin}[0]{\textcolor{\paramcolor}{\symCurrent_{\indexGridLines\indexGridNode \indexGridNodeTwo}^{\imagangle\minss}  }}
\newcommand{\Iijkanglemax}[0]{\textcolor{\paramcolor}{\symCurrent_{\indexGridLines\indexGridNode \indexGridNodeTwo}^{\imagangle\maxss}  }}

\newcommand{\IijkSDPseq}[0]{\textcolor{\complexcolor}{\mathbf{\symCurrent}_{\indexGridLines\indexGridNode \indexGridNodeTwo   }^{\fortescuess}  }}

\newcommand{\IijkSDPreal}[0]{\textcolor{black}{\mathbf{\symCurrent}_{\indexGridLines\indexGridNode \indexGridNodeTwo   }^{\realss}  }}
\newcommand{\IijkSDPimag}[0]{\textcolor{black}{\mathbf{\symCurrent}_{\indexGridLines\indexGridNode \indexGridNodeTwo   }^{\imagss}  }}

\newcommand{\IijkSDP}[0]{\textcolor{\complexcolor}{\mathbf{\symCurrent}_{\indexGridLines\indexGridNode \indexGridNodeTwo   }^{}  }}
\newcommand{\IijkSDPH}[0]{\textcolor{\complexcolor}{(\mathbf{\symCurrent}_{\indexGridLines\indexGridNode \indexGridNodeTwo   })^{\hermitiantranspose}  }}

\newcommand{\IijkSDPref}[0]{\textcolor{\complexparamcolor}{\mathbf{\symCurrent}_{\indexGridLines\indexGridNode \indexGridNodeTwo   }^{\refss}  }}
\newcommand{\IijkSDPrefre}[0]{\textcolor{\paramcolor}{\mathbf{\symCurrent}_{\indexGridLines\indexGridNode \indexGridNodeTwo   }^{\refss, \realss}  }}
\newcommand{\IijkSDPrefim}[0]{\textcolor{\paramcolor}{\mathbf{\symCurrent}_{\indexGridLines\indexGridNode \indexGridNodeTwo   }^{\refss, \imagss}  }}

\newcommand{\IijksSDPref}[0]{\textcolor{\complexparamcolor}{\mathbf{\symCurrent}_{\indexGridLines\indexGridNode \indexGridNodeTwo   }^{\seriesss, \refss}  }}
\newcommand{\IijksSDPrefre}[0]{\textcolor{\paramcolor}{\mathbf{\symCurrent}_{\indexGridLines\indexGridNode \indexGridNodeTwo   }^{\seriesss, \refss, \realss}  }}
\newcommand{\IijksSDPrefim}[0]{\textcolor{\paramcolor}{\mathbf{\symCurrent}_{\indexGridLines\indexGridNode \indexGridNodeTwo   }^{\seriesss, \refss, \imagss}  }}

\newcommand{\IijkshSDPref}[0]{\textcolor{\complexparamcolor}{\mathbf{\symCurrent}_{\indexGridLines\indexGridNode \indexGridNodeTwo   }^{\shuntss, \refss}  }}
\newcommand{\IijkshSDPrefre}[0]{\textcolor{\paramcolor}{\mathbf{\symCurrent}_{\indexGridLines\indexGridNode \indexGridNodeTwo   }^{\shuntss, \refss, \realss}  }}
\newcommand{\IijkshSDPrefim}[0]{\textcolor{\paramcolor}{\mathbf{\symCurrent}_{\indexGridLines\indexGridNode \indexGridNodeTwo   }^{\shuntss, \refss, \imagss}  }}

\newcommand{\IijkSDPdelta}[0]{\textcolor{\complexcolor}{\mathbf{\symCurrent}_{\indexGridLines\indexGridNode \indexGridNodeTwo   }^{\Delta}  }}
\newcommand{\IijkSDPdeltare}[0]{\textcolor{black}{\mathbf{\symCurrent}_{\indexGridLines\indexGridNode \indexGridNodeTwo   }^{\Delta, \realss}  }}
\newcommand{\IijkSDPdeltaim}[0]{\textcolor{black}{\mathbf{\symCurrent}_{\indexGridLines\indexGridNode \indexGridNodeTwo   }^{\Delta, \imagss}  }}

\newcommand{\IijksSDPdelta}[0]{\textcolor{\complexcolor}{\mathbf{\symCurrent}_{\indexGridLines\indexGridNode \indexGridNodeTwo   }^{\seriesss, \Delta}  }}
\newcommand{\IijksSDPdeltare}[0]{\textcolor{black}{\mathbf{\symCurrent}_{\indexGridLines\indexGridNode \indexGridNodeTwo   }^{\seriesss, \Delta, \realss}  }}
\newcommand{\IijksSDPdeltaim}[0]{\textcolor{black}{\mathbf{\symCurrent}_{\indexGridLines\indexGridNode \indexGridNodeTwo   }^{\seriesss, \Delta, \imagss}  }}

\newcommand{\IsqijkSDPdelta}[0]{\textcolor{\complexcolor}{\mathbf{\symCurrentSquared}_{\indexGridLines\indexGridNode \indexGridNodeTwo   }^{\Delta}  }}
\newcommand{\IsqijkSDPdeltare}[0]{\textcolor{black}{\mathbf{\symCurrentSquared}_{\indexGridLines\indexGridNode \indexGridNodeTwo   }^{\Delta, \realss}  }}
\newcommand{\IsqijkSDPdeltaim}[0]{\textcolor{black}{\mathbf{\symCurrentSquared}_{\indexGridLines\indexGridNode \indexGridNodeTwo   }^{\Delta, \imagss}  }}

\newcommand{\IsqijksSDPdelta}[0]{\textcolor{\complexcolor}{\mathbf{\symCurrentSquared}_{\indexGridLines\indexGridNode \indexGridNodeTwo   }^{\seriesss, \Delta}  }}
\newcommand{\IsqijksSDPdeltare}[0]{\textcolor{black}{\mathbf{\symCurrentSquared}_{\indexGridLines\indexGridNode \indexGridNodeTwo   }^{\seriesss, \Delta, \realss}  }}
\newcommand{\IsqijksSDPdeltaim}[0]{\textcolor{black}{\mathbf{\symCurrentSquared}_{\indexGridLines\indexGridNode \indexGridNodeTwo   }^{\seriesss, \Delta, \imagss}  }}

\newcommand{\IjikSDP}[0]{\textcolor{\complexcolor}{\mathbf{\symCurrent}_{\indexGridLines\indexGridNodeTwo\indexGridNode    }^{}  }}

\newcommand{\IijkrefSDP}[0]{\textcolor{\complexparamcolor}{\mathbf{\symCurrent}_{\indexGridLines\indexGridNode \indexGridNodeTwo   }^{\refss}  }}

\newcommand{\IijksrefSDP}[0]{\textcolor{\complexparamcolor}{\mathbf{\symCurrent}_{\indexGridLines\indexGridNode \indexGridNodeTwo, \seriesss   }^{\refss}  }}

\newcommand{\IijskSDPseq}[0]{\textcolor{\complexcolor}{\mathbf{\symCurrent}_{\indexGridLines\indexGridNode \indexGridNodeTwo   }^{\seriesss, \fortescuess}  }}

\newcommand{\IijskSDPH}[0]{\textcolor{\complexcolor}{(\mathbf{\symCurrent}_{\indexGridLines\indexGridNode \indexGridNodeTwo   }^{\seriesss})^{\hermitiantranspose}  }}
\newcommand{\IijskSDP}[0]{\textcolor{\complexcolor}{\mathbf{\symCurrent}_{\indexGridLines\indexGridNode \indexGridNodeTwo   }^{\seriesss}  }}
\newcommand{\IijskSDPreal}[0]{\textcolor{black}{\mathbf{\symCurrent}_{\indexGridLines\indexGridNode \indexGridNodeTwo   }^{\seriesss, \realss}  }}
\newcommand{\IijskSDPimag}[0]{\textcolor{black}{\mathbf{\symCurrent}_{\indexGridLines\indexGridNode \indexGridNodeTwo   }^{\seriesss, \imagss}  }}

\newcommand{\IjiskSDP}[0]{\textcolor{\complexcolor}{\mathbf{\symCurrent}_{\indexGridLines\indexGridNodeTwo\indexGridNode    }^{\seriesss}  }}

\newcommand{\IijshkSDP}[0]{\textcolor{\complexcolor}{\mathbf{\symCurrent}_{\indexGridLines\indexGridNode \indexGridNodeTwo   }^{\shuntss}  }}
\newcommand{\IjishkSDP}[0]{\textcolor{\complexcolor}{\mathbf{\symCurrent}_{ \indexGridLines\indexGridNodeTwo \indexGridNode   }^{\shuntss}  }}

\newcommand{\IijshkSDPreal}[0]{\textcolor{black}{\mathbf{\symCurrent}_{\indexGridLines\indexGridNode \indexGridNodeTwo   }^{\shuntss, \realss}  }}

\newcommand{\IijshkSDPimag}[0]{\textcolor{black}{\mathbf{\symCurrent}_{\indexGridLines\indexGridNode \indexGridNodeTwo   }^{\shuntss, \imagss}  }}

\newcommand{\Iijk}[0]{\textcolor{\complexcolor}{\symCurrent_{\indexGridLines\indexGridNode \indexGridNodeTwo}^{}  }}
\newcommand{\Iijsk}[0]{\textcolor{\complexcolor}{\symCurrent_{\indexGridLines\indexGridNode \indexGridNodeTwo,\seriesss}^{}  }}
\newcommand{\Iijshk}[0]{\textcolor{\complexcolor}{\symCurrent_{\indexGridLines\indexGridNode \indexGridNodeTwo,\shuntss}^{}  }}
\newcommand{\Iijsqk}[0]{\textcolor{black}{\symCurrent_{\indexGridLines\indexGridNode \indexGridNodeTwo}^{2}  }}
\newcommand{\Iijsqsk}[0]{\textcolor{black}{\symCurrent_{\indexGridLines\indexGridLines,\seriesss}^{2}  }}
\newcommand{\Iijsqshk}[0]{\textcolor{black}{\symCurrent_{\indexGridLines\indexGridNode \indexGridNodeTwo,\shuntss}^{2}  }}

\newcommand{\IijkH}[0]{\textcolor{\complexcolor}{\symCurrent_{\indexGridLines\indexGridNode \indexGridNodeTwo}^{\hermitiantranspose}  }}
\newcommand{\Ijisk}[0]{\textcolor{\complexcolor}{\symCurrent_{\indexGridLines\indexGridNodeTwo\indexGridNode ,\seriesss}^{}  }}
\newcommand{\Ijishk}[0]{\textcolor{\complexcolor}{\symCurrent_{\indexGridLines\indexGridNodeTwo\indexGridNode ,\shuntss}^{}  }}

\newcommand{\Iijkmax}[0]{\textcolor{\boundscolor}{\symCurrent_{\indexGridLines\indexGridNode \indexGridNodeTwo}^{\maxss}  }}
\newcommand{\Iijkmin}[0]{\textcolor{\boundscolor}{\symCurrent_{\indexGridLines\indexGridNode \indexGridNodeTwo}^{\minss}  }}
\newcommand{\Iijrated}[0]{\textcolor{\sizingcolor}{\symCurrent_{\indexGridLines\indexGridNode \indexGridNodeTwo}^{\ratedss} }}
\newcommand{\Ijirated}[0]{\textcolor{\sizingcolor}{\symCurrent_{\indexGridLines\indexGridNodeTwo  \indexGridNode }^{\ratedss} }}
\newcommand{\IijratedA}[0]{\textcolor{\sizingcolor}{\symCurrent_{\indexGridLines\indexGridNode \indexGridNodeTwo,\symPhaseA}^{\ratedss} }}
\newcommand{\IjiratedA}[0]{\textcolor{\sizingcolor}{\symCurrent_{\indexGridLines\indexGridNodeTwo  \indexGridNode,\symPhaseA}^{\ratedss} }}
\newcommand{\IijratedB}[0]{\textcolor{\sizingcolor}{\symCurrent_{\indexGridLines\indexGridNode \indexGridNodeTwo,\symPhaseB}^{\ratedss} }}
\newcommand{\IjiratedB}[0]{\textcolor{\sizingcolor}{\symCurrent_{\indexGridLines\indexGridNodeTwo  \indexGridNode,\symPhaseB}^{\ratedss} }}
\newcommand{\IijratedC}[0]{\textcolor{\sizingcolor}{\symCurrent_{\indexGridLines\indexGridNode \indexGridNodeTwo,\symPhaseC}^{\ratedss} }}
\newcommand{\IjiratedC}[0]{\textcolor{\sizingcolor}{\symCurrent_{\indexGridLines\indexGridNodeTwo  \indexGridNode,\symPhaseC}^{\ratedss} }}

\newcommand{\IijratedN}[0]{\textcolor{\sizingcolor}{\symCurrent_{\indexGridLines\indexGridNode \indexGridNodeTwo,\symPhaseN}^{\ratedss} }}

\newcommand{\IijsqskSDPreal}[0]{\textcolor{black}{\mathbf{\symCurrentSOCP}_{\indexGridLines  }^{\seriesss,\realss}  }}
\newcommand{\IijsqskSDPimag}[0]{\textcolor{black}{\mathbf{\symCurrentSOCP}_{\indexGridLines  }^{\seriesss,\imagss}  }}

\newcommand{\IijsqkSDPreal}[0]{\textcolor{black}{\mathbf{\symCurrentSOCP}_{\indexGridLines\indexGridNode \indexGridNodeTwo }^{\realss}  }}
\newcommand{\IijsqkSDPimag}[0]{\textcolor{black}{\mathbf{\symCurrentSOCP}_{\indexGridLines\indexGridNode \indexGridNodeTwo }^{\imagss}  }}

\newcommand{\IijsqkSDPAdot}[0]{\textcolor{\complexcolor}{\mathbf{\symCurrentSOCP}_{\indexGridLines\indexGridNode \indexGridNodeTwo }^{\symPhaseA\!\cdot}  }}
\newcommand{\IijsqkSDPBdot}[0]{\textcolor{\complexcolor}{\mathbf{\symCurrentSOCP}_{\indexGridLines\indexGridNode \indexGridNodeTwo }^{\symPhaseB\!\cdot}  }}
\newcommand{\IijsqkSDPCdot}[0]{\textcolor{\complexcolor}{\mathbf{\symCurrentSOCP}_{\indexGridLines\indexGridNode \indexGridNodeTwo }^{\symPhaseC\!\cdot}  }}

\newcommand{\IjisqkSDP}[0]{\textcolor{\complexcolor}{\mathbf{\symCurrentSOCP}_{\indexGridLines \indexGridNodeTwo \indexGridNode  }^{}  }}
\newcommand{\IijsqkSDP}[0]{\textcolor{\complexcolor}{\mathbf{\symCurrentSOCP}_{\indexGridLines\indexGridNode \indexGridNodeTwo }^{}  }}
\newcommand{\IijsqkmaxSDP}[0]{\textcolor{\complexparamcolor}{\mathbf{\symCurrentSOCP}_{\indexGridLines\indexGridNode \indexGridNodeTwo }^{\maxss}  }}

\newcommand{\IijkmaxSDP}[0]{\textcolor{\paramcolor}{\mathbf{\symCurrent}_{\indexGridLines\indexGridNode \indexGridNodeTwo }^{\maxss}  }}

\newcommand{\IjisqkSDPreal}[0]{\textcolor{black}{\mathbf{\symCurrentSOCP}_{\indexGridLines \indexGridNodeTwo \indexGridNode  }^{\realss}  }}

\newcommand{\IijsqshkSDP}[0]{\textcolor{\complexcolor}{\mathbf{\symCurrentSOCP}_{\indexGridLines\indexGridNode \indexGridNodeTwo }^{\shuntss}  }}
\newcommand{\IjisqshkSDP}[0]{\textcolor{\complexcolor}{\mathbf{\symCurrentSOCP}_{\indexGridLines\indexGridNodeTwo \indexGridNode  }^{\shuntss}  }}

\newcommand{\IusqkSDPAdot}[0]{\textcolor{\complexcolor}{\mathbf{\symCurrentSOCP}_{\indexUnit  }^{\symPhaseA\!\cdot}  }}
\newcommand{\IusqkSDPBdot}[0]{\textcolor{\complexcolor}{\mathbf{\symCurrentSOCP}_{\indexUnit  }^{\symPhaseB\!\cdot}  }}
\newcommand{\IusqkSDPCdot}[0]{\textcolor{\complexcolor}{\mathbf{\symCurrentSOCP}_{\indexUnit  }^{\symPhaseC\!\cdot}  }}
\newcommand{\IusqkSDPNdot}[0]{\textcolor{\complexcolor}{\mathbf{\symCurrentSOCP}_{\indexUnit  }^{\symPhaseN\!\cdot}  }}

\newcommand{\IusqkSDPdotA}[0]{\textcolor{\complexcolor}{\mathbf{\symCurrentSOCP}_{\indexUnit  }^{\!\cdot\!\symPhaseA}  }}
\newcommand{\IusqkSDPdotB}[0]{\textcolor{\complexcolor}{\mathbf{\symCurrentSOCP}_{\indexUnit  }^{\!\cdot\!\symPhaseB}  }}
\newcommand{\IusqkSDPdotC}[0]{\textcolor{\complexcolor}{\mathbf{\symCurrentSOCP}_{\indexUnit  }^{\!\cdot\!\symPhaseC}  }}
\newcommand{\IusqkSDPdotN}[0]{\textcolor{\complexcolor}{\mathbf{\symCurrentSOCP}_{\indexUnit  }^{\!\cdot\!\symPhaseN}  }}

\newcommand{\IijsqskSDPseq}[0]{\textcolor{\complexcolor}{\mathbf{\symCurrentSOCP}_{\indexGridLines  }^{\seriesss, \fortescuess}  }}
\newcommand{\IusqkSDP}[0]{\textcolor{\complexcolor}{\mathbf{\symCurrentSOCP}_{\indexUnit  }  }}
\newcommand{\IusqkSDPreal}[0]{\textcolor{black}{\mathbf{\symCurrentSOCP}_{\indexUnit  }^{\realss}  }}
\newcommand{\IusqkSDPimag}[0]{\textcolor{black}{\mathbf{\symCurrentSOCP}_{\indexUnit  }^{\imagss}  }}

\newcommand{\IusqkSDPdelta}[0]{\textcolor{\complexcolor}{\mathbf{\symCurrentSOCP}_{\indexUnit  }^{\Delta}  }}
\newcommand{\IusqkSDPdeltareal}[0]{\textcolor{black}{\mathbf{\symCurrentSOCP}_{\indexUnit  }^{\Delta\realss}  }}
\newcommand{\IusqkSDPdeltaimag}[0]{\textcolor{black}{\mathbf{\symCurrentSOCP}_{\indexUnit  }^{\Delta\imagss}  }}

\newcommand{\IijsqskSDP}[0]{\textcolor{\complexcolor}{\mathbf{\symCurrentSOCP}_{\indexGridLines  }^{\seriesss}  }}
\newcommand{\IjisqskSDP}[0]{\textcolor{\complexcolor}{\mathbf{\symCurrentSOCP}_{\indexGridLines  }^{\seriesss}  }}
\newcommand{\Iijsqsksocp}[0]{\textcolor{black}{\symCurrentSOCP_{\indexGridLines\indexGridNode \indexGridNodeTwo ,\seriesss}  }}

\newcommand{\Iijsqksocp}[0]{\textcolor{black}{\symCurrentSOCP_{\indexGridLines\indexGridNode \indexGridNodeTwo}  }}

\newcommand{\IijsqsksocpAA}[0]{\textcolor{black}{\symCurrentSOCP_{\indexGridLines\indexGridNode \indexGridNodeTwo,\seriesss, \symPhaseA\symPhaseA}  }}
\newcommand{\IijsqsksocpAB}[0]{\textcolor{\complexcolor}{\symCurrentSOCP_{\indexGridLines\indexGridNode \indexGridNodeTwo,\seriesss, \symPhaseA\symPhaseB}  }}
\newcommand{\IijsqsksocpAC}[0]{\textcolor{\complexcolor}{\symCurrentSOCP_{\indexGridLines\indexGridNode \indexGridNodeTwo,\seriesss, \symPhaseA\symPhaseC}  }}
\newcommand{\IijsqsksocpBA}[0]{\textcolor{\complexcolor}{\symCurrentSOCP_{\indexGridLines\indexGridNode \indexGridNodeTwo,\seriesss, \symPhaseB\symPhaseA}  }}
\newcommand{\IijsqsksocpBB}[0]{\textcolor{black}{\symCurrentSOCP_{\indexGridLines\indexGridNode \indexGridNodeTwo,\seriesss, \symPhaseB\symPhaseB}  }}
\newcommand{\IijsqsksocpBC}[0]{\textcolor{\complexcolor}{\symCurrentSOCP_{\indexGridLines\indexGridNode \indexGridNodeTwo,\seriesss, \symPhaseB\symPhaseC}  }}
\newcommand{\IijsqsksocpCA}[0]{\textcolor{\complexcolor}{\symCurrentSOCP_{\indexGridLines\indexGridNode \indexGridNodeTwo,\seriesss, \symPhaseC\symPhaseA}  }}
\newcommand{\IijsqsksocpCB}[0]{\textcolor{\complexcolor}{\symCurrentSOCP_{\indexGridLines\indexGridNode \indexGridNodeTwo,\seriesss, \symPhaseC\symPhaseB}  }}
\newcommand{\IijsqsksocpCC}[0]{\textcolor{black}{\symCurrentSOCP_{\indexGridLines\indexGridNode \indexGridNodeTwo,\seriesss, \symPhaseC\symPhaseC}  }}

\newcommand{\IijsqksocpAAreal}[0]{\textcolor{black}{\symCurrentSOCP_{\indexGridLines\indexGridNode \indexGridNodeTwo, \symPhaseA\symPhaseA}^{\realss}  }}
\newcommand{\IijsqksocpBBreal}[0]{\textcolor{black}{\symCurrentSOCP_{\indexGridLines\indexGridNode \indexGridNodeTwo, \symPhaseB\symPhaseB}^{\realss}  }}
\newcommand{\IijsqksocpCCreal}[0]{\textcolor{black}{\symCurrentSOCP_{\indexGridLines\indexGridNode \indexGridNodeTwo, \symPhaseC\symPhaseC}^{\realss}  }}

\newcommand{\IijsqksocpNNreal}[0]{\textcolor{black}{\symCurrentSOCP_{\indexGridLines\indexGridNode \indexGridNodeTwo, \symPhaseN\symPhaseN}^{\realss}  }}
\newcommand{\IijsqksocpGGreal}[0]{\textcolor{black}{\symCurrentSOCP_{\indexGridLines\indexGridNode \indexGridNodeTwo, \symPhaseG\symPhaseG}^{\realss}  }}
\newcommand{\kifsqksocpGGreal}[0]{\textcolor{black}{\symCurrentSOCP_{k\indexGridNode f, \symPhaseG\symPhaseG}^{\realss}  }}

\newcommand{\IijingsqksocpGG}[0]{\textcolor{\complexcolor}{\symCurrentSOCP_{\indexGridLines\indexGridNode \indexGridNodeTwo, \symPhaseN\symPhaseG}  }}

\newcommand{\kifingsqksocpGG}[0]{\textcolor{\complexcolor}{\symCurrentSOCP_{k\indexGridNode f, \symPhaseN\symPhaseG}  }}
    \newcommand{\parm}{\mathord{\color{black!33}\bullet}}%

\newcommand{\IijsqsksocpAAreal}[0]{\textcolor{black}{\symCurrentSOCP_{\indexGridLines\indexGridNode \indexGridNodeTwo, \symPhaseA\symPhaseA}^{\seriesss,\realss}  }}
\newcommand{\IijsqsksocpABreal}[0]{\textcolor{black}{\symCurrentSOCP_{\indexGridLines\indexGridNode \indexGridNodeTwo, \symPhaseA\symPhaseB}^{\seriesss,\realss}  }}
\newcommand{\IijsqsksocpACreal}[0]{\textcolor{black}{\symCurrentSOCP_{\indexGridLines\indexGridNode \indexGridNodeTwo, \symPhaseA\symPhaseC}^{\seriesss,\realss}  }}
\newcommand{\IijsqsksocpBBreal}[0]{\textcolor{black}{\symCurrentSOCP_{\indexGridLines\indexGridNode \indexGridNodeTwo, \symPhaseB\symPhaseB}^{\seriesss,\realss}  }}
\newcommand{\IijsqsksocpBCreal}[0]{\textcolor{black}{\symCurrentSOCP_{\indexGridLines\indexGridNode \indexGridNodeTwo, \symPhaseB\symPhaseC}^{\seriesss,\realss}  }}
\newcommand{\IijsqsksocpCCreal}[0]{\textcolor{black}{\symCurrentSOCP_{\indexGridLines\indexGridNode \indexGridNodeTwo, \symPhaseC\symPhaseC}^{\seriesss,\realss}  }}
\newcommand{\IijsqsksocpNNreal}[0]{\textcolor{black}{\symCurrentSOCP_{\indexGridLines\indexGridNode \indexGridNodeTwo, \symPhaseN\symPhaseN}^{\seriesss,\realss}  }}

\newcommand{\IijsqsksocpANreal}[0]{\textcolor{black}{\symCurrentSOCP_{\indexGridLines\indexGridNode \indexGridNodeTwo, \symPhaseA\symPhaseN}^{\seriesss,\realss}  }}
\newcommand{\IijsqsksocpBNreal}[0]{\textcolor{black}{\symCurrentSOCP_{\indexGridLines\indexGridNode \indexGridNodeTwo, \symPhaseB\symPhaseN}^{\seriesss,\realss}  }}
\newcommand{\IijsqsksocpCNreal}[0]{\textcolor{black}{\symCurrentSOCP_{\indexGridLines\indexGridNode \indexGridNodeTwo, \symPhaseC\symPhaseN}^{\seriesss,\realss}  }}

\newcommand{\IijsqsksocpANimag}[0]{\textcolor{black}{\symCurrentSOCP_{\indexGridLines\indexGridNode \indexGridNodeTwo, \symPhaseA\symPhaseN}^{\seriesss,\imagss}  }}
\newcommand{\IijsqsksocpBNimag}[0]{\textcolor{black}{\symCurrentSOCP_{\indexGridLines\indexGridNode \indexGridNodeTwo, \symPhaseB\symPhaseN}^{\seriesss,\imagss}  }}
\newcommand{\IijsqsksocpCNimag}[0]{\textcolor{black}{\symCurrentSOCP_{\indexGridLines\indexGridNode \indexGridNodeTwo, \symPhaseC\symPhaseN}^{\seriesss,\imagss}  }}

\newcommand{\IijsqsksocpGGreal}[0]{\textcolor{black}{\symCurrentSOCP_{\indexGridLines\indexGridNode \indexGridNodeTwo, \symPhaseG\symPhaseG}^{\seriesss,\realss}  }}
\newcommand{\IisqsksocpNGreal}[0]{\textcolor{black}{\symCurrentSOCP_{\indexGridNode , \symPhaseN\symPhaseG}  }}

\newcommand{\IjisqsksocpBBreal}[0]{\textcolor{black}{\symCurrentSOCP_{\indexGridLines\indexGridNodeTwo\indexGridNode , \symPhaseB\symPhaseB}^{\seriesss,\realss}  }}

\newcommand{\IijsqsksocpABimag}[0]{\textcolor{black}{\symCurrentSOCP_{\indexGridLines\indexGridNode \indexGridNodeTwo, \symPhaseA\symPhaseB}^{\seriesss,\imagss}  }}
\newcommand{\IijsqsksocpACimag}[0]{\textcolor{black}{\symCurrentSOCP_{\indexGridLines\indexGridNode \indexGridNodeTwo, \symPhaseA\symPhaseC}^{\seriesss,\imagss}  }}
\newcommand{\IijsqsksocpBCimag}[0]{\textcolor{black}{\symCurrentSOCP_{\indexGridLines\indexGridNode \indexGridNodeTwo, \symPhaseB\symPhaseC}^{\seriesss,\imagss}  }}

\newcommand{\VjsqksocpAAreal}[0]{\textcolor{black}{\symVoltageSOCP_{\indexGridNodeTwo, \symPhaseA\symPhaseA}^{\realss}  }}
\newcommand{\VjsqksocpBBreal}[0]{\textcolor{black}{\symVoltageSOCP_{\indexGridNodeTwo, \symPhaseB\symPhaseB}^{\realss}  }}
\newcommand{\VjsqksocpCCreal}[0]{\textcolor{black}{\symVoltageSOCP_{\indexGridNodeTwo, \symPhaseC\symPhaseC}^{\realss}  }}

\newcommand{\VisqksocpNNreal}[0]{\textcolor{black}{\symVoltageSOCP_{\indexGridNode, \symPhaseN\symPhaseN}^{\realss}  }}
\newcommand{\VisqksocpGGreal}[0]{\textcolor{black}{\symVoltageSOCP_{\indexGridNode, \symPhaseG\symPhaseG}^{\realss}  }}

\newcommand{\VisqksocpAAreal}[0]{\textcolor{black}{\symVoltageSOCP_{\indexGridNode, \symPhaseA\symPhaseA}^{\realss}  }}
\newcommand{\VisqksocpABreal}[0]{\textcolor{black}{\symVoltageSOCP_{\indexGridNode, \symPhaseA\symPhaseB}^{\realss}  }}
\newcommand{\VisqksocpACreal}[0]{\textcolor{black}{\symVoltageSOCP_{\indexGridNode, \symPhaseA\symPhaseC}^{\realss}  }}

\newcommand{\VisqksocpANreal}[0]{\textcolor{black}{\symVoltageSOCP_{\indexGridNode, \symPhaseA\symPhaseN}^{\realss}  }}
\newcommand{\VisqksocpBNreal}[0]{\textcolor{black}{\symVoltageSOCP_{\indexGridNode, \symPhaseB\symPhaseN}^{\realss}  }}
\newcommand{\VisqksocpCNreal}[0]{\textcolor{black}{\symVoltageSOCP_{\indexGridNode, \symPhaseC\symPhaseN}^{\realss}  }}
\newcommand{\VisqksocpANimag}[0]{\textcolor{black}{\symVoltageSOCP_{\indexGridNode, \symPhaseA\symPhaseN}^{\imagss}  }}
\newcommand{\VisqksocpBNimag}[0]{\textcolor{black}{\symVoltageSOCP_{\indexGridNode, \symPhaseB\symPhaseN}^{\imagss}  }}
\newcommand{\VisqksocpCNimag}[0]{\textcolor{black}{\symVoltageSOCP_{\indexGridNode, \symPhaseC\symPhaseN}^{\imagss}  }}

\newcommand{\VisqksocpBBreal}[0]{\textcolor{black}{\symVoltageSOCP_{\indexGridNode, \symPhaseB\symPhaseB}^{\realss}  }}
\newcommand{\VisqksocpBCreal}[0]{\textcolor{black}{\symVoltageSOCP_{\indexGridNode, \symPhaseB\symPhaseC}^{\realss}  }}
\newcommand{\VisqksocpCCreal}[0]{\textcolor{black}{\symVoltageSOCP_{\indexGridNode, \symPhaseC\symPhaseC}^{\realss}  }}

\newcommand{\VisqksocpABimag}[0]{\textcolor{black}{\symVoltageSOCP_{\indexGridNode, \symPhaseA\symPhaseB}^{\imagss}  }}
\newcommand{\VisqksocpACimag}[0]{\textcolor{black}{\symVoltageSOCP_{\indexGridNode, \symPhaseA\symPhaseC}^{\imagss}  }}
\newcommand{\VisqksocpBCimag}[0]{\textcolor{black}{\symVoltageSOCP_{\indexGridNode, \symPhaseB\symPhaseC}^{\imagss}  }}

\newcommand{\VisqksocpAA}[0]{\textcolor{black}{\symVoltageSOCP_{\indexGridNode, \symPhaseA\symPhaseA} }}
\newcommand{\VisqksocpAB}[0]{\textcolor{\complexcolor}{\symVoltageSOCP_{\indexGridNode, \symPhaseA\symPhaseB}}}
\newcommand{\VisqksocpAC}[0]{\textcolor{\complexcolor}{\symVoltageSOCP_{\indexGridNode, \symPhaseA\symPhaseC}  }}
\newcommand{\VisqksocpBA}[0]{\textcolor{\complexcolor}{\symVoltageSOCP_{\indexGridNode, \symPhaseB\symPhaseA} }}
\newcommand{\VisqksocpBB}[0]{\textcolor{black}{\symVoltageSOCP_{\indexGridNode, \symPhaseB\symPhaseB} }}
\newcommand{\VisqksocpBC}[0]{\textcolor{\complexcolor}{\symVoltageSOCP_{\indexGridNode, \symPhaseB\symPhaseC}}}
\newcommand{\VisqksocpCA}[0]{\textcolor{\complexcolor}{\symVoltageSOCP_{\indexGridNode, \symPhaseC\symPhaseA}  }}
\newcommand{\VisqksocpCB}[0]{\textcolor{\complexcolor}{\symVoltageSOCP_{\indexGridNode, \symPhaseC\symPhaseB}  }}
\newcommand{\VisqksocpCC}[0]{\textcolor{black}{\symVoltageSOCP_{\indexGridNode, \symPhaseC\symPhaseC}  }}

\newcommand{\VijsqksocpAAreal}[0]{\textcolor{black}{\symVoltageSOCP_{\indexGridNode\indexGridNodeTwo, \symPhaseA\symPhaseA}^{\realss}  }}
\newcommand{\VijsqksocpABreal}[0]{\textcolor{black}{\symVoltageSOCP_{\indexGridNode\indexGridNodeTwo, \symPhaseA\symPhaseB}^{\realss}  }}
\newcommand{\VijsqksocpACreal}[0]{\textcolor{black}{\symVoltageSOCP_{\indexGridNode\indexGridNodeTwo, \symPhaseA\symPhaseC}^{\realss}  }}

\newcommand{\VijsqksocpBAreal}[0]{\textcolor{black}{\symVoltageSOCP_{\indexGridNode\indexGridNodeTwo, \symPhaseB\symPhaseA}^{\realss}  }}
\newcommand{\VijsqksocpBBreal}[0]{\textcolor{black}{\symVoltageSOCP_{\indexGridNode\indexGridNodeTwo, \symPhaseB\symPhaseB}^{\realss}  }}
\newcommand{\VijsqksocpBCreal}[0]{\textcolor{black}{\symVoltageSOCP_{\indexGridNode\indexGridNodeTwo, \symPhaseB\symPhaseC}^{\realss}  }}

\newcommand{\VijsqksocpCAreal}[0]{\textcolor{black}{\symVoltageSOCP_{\indexGridNode\indexGridNodeTwo, \symPhaseC\symPhaseA}^{\realss}  }}
\newcommand{\VijsqksocpCBreal}[0]{\textcolor{black}{\symVoltageSOCP_{\indexGridNode\indexGridNodeTwo, \symPhaseC\symPhaseB}^{\realss}  }}
\newcommand{\VijsqksocpCCreal}[0]{\textcolor{black}{\symVoltageSOCP_{\indexGridNode\indexGridNodeTwo, \symPhaseC\symPhaseC}^{\realss}  }}

\newcommand{\VjsqksocpABimag}[0]{\textcolor{black}{\symVoltageSOCP_{\indexGridNodeTwo, \symPhaseA\symPhaseB}^{\imagss}  }}
\newcommand{\VjsqksocpACimag}[0]{\textcolor{black}{\symVoltageSOCP_{\indexGridNodeTwo, \symPhaseA\symPhaseC}^{\imagss}  }}
\newcommand{\VjsqksocpBCimag}[0]{\textcolor{black}{\symVoltageSOCP_{\indexGridNodeTwo, \symPhaseB\symPhaseC}^{\imagss}  }}

\newcommand{\VjsqksocpABreal}[0]{\textcolor{black}{\symVoltageSOCP_{\indexGridNodeTwo, \symPhaseA\symPhaseB}^{\realss}  }}
\newcommand{\VjsqksocpACreal}[0]{\textcolor{black}{\symVoltageSOCP_{\indexGridNodeTwo, \symPhaseA\symPhaseC}^{\realss}  }}
\newcommand{\VjsqksocpBCreal}[0]{\textcolor{black}{\symVoltageSOCP_{\indexGridNodeTwo, \symPhaseB\symPhaseC}^{\realss}  }}

\newcommand{\VijsqksocpAAimag}[0]{\textcolor{black}{\symVoltageSOCP_{\indexGridNode\indexGridNodeTwo, \symPhaseA\symPhaseA}^{\imagss}  }}
\newcommand{\VijsqksocpABimag}[0]{\textcolor{black}{\symVoltageSOCP_{\indexGridNode\indexGridNodeTwo, \symPhaseA\symPhaseB}^{\imagss}  }}
\newcommand{\VijsqksocpACimag}[0]{\textcolor{black}{\symVoltageSOCP_{\indexGridNode\indexGridNodeTwo, \symPhaseA\symPhaseC}^{\imagss}  }}

\newcommand{\VijsqksocpBAimag}[0]{\textcolor{black}{\symVoltageSOCP_{\indexGridNode\indexGridNodeTwo, \symPhaseB\symPhaseA}^{\imagss}  }}
\newcommand{\VijsqksocpBBimag}[0]{\textcolor{black}{\symVoltageSOCP_{\indexGridNode\indexGridNodeTwo, \symPhaseB\symPhaseB}^{\imagss}  }}
\newcommand{\VijsqksocpBCimag}[0]{\textcolor{black}{\symVoltageSOCP_{\indexGridNode\indexGridNodeTwo, \symPhaseB\symPhaseC}^{\imagss}  }}

\newcommand{\VijsqksocpCAimag}[0]{\textcolor{black}{\symVoltageSOCP_{\indexGridNode\indexGridNodeTwo, \symPhaseC\symPhaseA}^{\imagss}  }}
\newcommand{\VijsqksocpCBimag}[0]{\textcolor{black}{\symVoltageSOCP_{\indexGridNode\indexGridNodeTwo, \symPhaseC\symPhaseB}^{\imagss}  }}
\newcommand{\VijsqksocpCCimag}[0]{\textcolor{black}{\symVoltageSOCP_{\indexGridNode\indexGridNodeTwo, \symPhaseC\symPhaseC}^{\imagss}  }}

\newcommand{\IijsqskA}[0]{\textcolor{black}{\symCurrent_{\indexGridLines\indexGridNode \indexGridNodeTwo,\seriesss, \symPhaseA}^{2}  }}
\newcommand{\IijsqskB}[0]{\textcolor{black}{\symCurrent_{\indexGridLines\indexGridNode \indexGridNodeTwo,\seriesss, \symPhaseB}^{2}  }}
\newcommand{\IijsqskC}[0]{\textcolor{black}{\symCurrent_{\indexGridLines\indexGridNode \indexGridNodeTwo,\seriesss, \symPhaseC}^{2}  }}

\newcommand{\Sijrated}[0]{\textcolor{\sizingcolor}{\symApparentPower_{\indexGridLines\indexGridNode \indexGridNodeTwo}^{\ratedss} }}
\newcommand{\Sjirated}[0]{\textcolor{\sizingcolor}{\symApparentPower_{\indexGridLines\indexGridNodeTwo \indexGridNode }^{\ratedss} }}
\newcommand{\SijratedA}[0]{\textcolor{\sizingcolor}{\symApparentPower_{\indexGridLines\indexGridNode \indexGridNodeTwo,\symPhaseA}^{\ratedss} }}
\newcommand{\SjiratedA}[0]{\textcolor{\sizingcolor}{\symApparentPower_{\indexGridLines\indexGridNodeTwo \indexGridNode ,\symPhaseA}^{\ratedss} }}
\newcommand{\SijratedB}[0]{\textcolor{\sizingcolor}{\symApparentPower_{\indexGridLines\indexGridNode \indexGridNodeTwo,\symPhaseB}^{\ratedss} }}
\newcommand{\SjiratedB}[0]{\textcolor{\sizingcolor}{\symApparentPower_{\indexGridLines\indexGridNodeTwo \indexGridNode ,\symPhaseB}^{\ratedss} }}
\newcommand{\SijratedC}[0]{\textcolor{\sizingcolor}{\symApparentPower_{\indexGridLines\indexGridNode \indexGridNodeTwo,\symPhaseC}^{\ratedss} }}
\newcommand{\SijratedN}[0]{\textcolor{\sizingcolor}{\symApparentPower_{\indexGridLines\indexGridNode \indexGridNodeTwo,\symPhaseN}^{\ratedss} }}
\newcommand{\SjiratedC}[0]{\textcolor{\sizingcolor}{\symApparentPower_{\indexGridLines\indexGridNodeTwo \indexGridNode ,\symPhaseC}^{\ratedss} }}

\newcommand{\VikSDPmin}[0]{\textcolor{\paramcolor}{\mathbf{\symVoltage}_{\indexGridNode }^{\minss}  }}
\newcommand{\VikSDPmax}[0]{\textcolor{\paramcolor}{\mathbf{\symVoltage}_{\indexGridNode }^{\maxss}  }}

\newcommand{\VikSDPwye}[0]{\textcolor{\complexcolor}{\mathbf{\symVoltage}_{\indexGridNode }^{\text{pn}}  }}

\newcommand{\VikSDPminwye}[0]{\textcolor{\paramcolor}{\mathbf{\symVoltage}_{\indexGridNode }^{\text{pn},\minss}  }}
\newcommand{\VikSDPmaxwye}[0]{\textcolor{\paramcolor}{\mathbf{\symVoltage}_{\indexGridNode }^{\text{pn},\maxss}  }}

\newcommand{\VikSDPmindelta}[0]{\textcolor{\paramcolor}{\mathbf{\symVoltage}_{\indexGridNode }^{\Delta\minss}  }}
\newcommand{\VikSDPmaxdelta}[0]{\textcolor{\paramcolor}{\mathbf{\symVoltage}_{\indexGridNode }^{\Delta\maxss}  }}

\newcommand{\VjkSDPmin}[0]{\textcolor{\paramcolor}{\mathbf{\symVoltage}_{\indexGridNodeTwo }^{\minss}  }}
\newcommand{\VjkSDPmax}[0]{\textcolor{\paramcolor}{\mathbf{\symVoltage}_{\indexGridNodeTwo }^{\maxss}  }}

\newcommand{\VjkSDPseq}[0]{\textcolor{\complexcolor}{\mathbf{\symVoltage}_{\indexGridNodeTwo }^{\fortescuess}  }}

\newcommand{\VikSDPdelta}[0]{\textcolor{\complexcolor}{\mathbf{\symVoltage}_{\indexGridNode }^{\Delta}  }}

\newcommand{\VikSDPseq}[0]{\textcolor{\complexcolor}{\mathbf{\symVoltage}_{\indexGridNode }^{\fortescuess}  }}
\newcommand{\VikSDP}[0]{\textcolor{\complexcolor}{\mathbf{\symVoltage}_{\indexGridNode }^{}  }}
\newcommand{\VikSDPreal}[0]{\textcolor{black}{\mathbf{\symVoltage}_{\indexGridNode }^{\realss}  }}
\newcommand{\VikSDPimag}[0]{\textcolor{black}{\mathbf{\symVoltage}_{\indexGridNode }^{\imagss}  }}
\newcommand{\VjkSDPreal}[0]{\textcolor{black}{\mathbf{\symVoltage}_{\indexGridNodeTwo }^{\realss}  }}
\newcommand{\VjkSDPimag}[0]{\textcolor{black}{\mathbf{\symVoltage}_{\indexGridNodeTwo }^{\imagss}  }}
\newcommand{\VikSDPH}[0]{\textcolor{\complexcolor}{\mathbf{\symVoltage}_{\indexGridNode }^{\hermitiantranspose}  }}
\newcommand{\VjkSDP}[0]{\textcolor{\complexcolor}{\mathbf{\symVoltage}_{\indexGridNodeTwo }^{}  }}
\newcommand{\VzkSDP}[0]{\textcolor{\complexcolor}{\mathbf{\symVoltage}_{z }^{}  }}
\newcommand{\VjkSDPH}[0]{\textcolor{\complexcolor}{\mathbf{\symVoltage}_{\indexGridNodeTwo }^{\hermitiantranspose}  }}
\newcommand{\Vik}[0]{\textcolor{black}{\symVoltage_{\indexGridNode }^{}  }}
\newcommand{\Vjk}[0]{\textcolor{black}{\symVoltage_{\indexGridNodeTwo }^{}  }}
\newcommand{\Vikdelta}[0]{\textcolor{black}{\Delta\symVoltage_{\indexGridNode }^{}  }}
\newcommand{\Vjkdelta}[0]{\textcolor{black}{\Delta\symVoltage_{\indexGridNodeTwo }^{}  }}
\newcommand{\Viacck}[0]{\textcolor{black}{\symVoltage_{\indexGridNode }^{'}  }}
\newcommand{\Vjacck}[0]{\textcolor{black}{\symVoltage_{\indexGridNodeTwo }^{'}  }}
\newcommand{\Viacckl}[0]{\textcolor{black}{\symVoltage_{\indexGridNode }^{\star}  }}
\newcommand{\Vjacckl}[0]{\textcolor{black}{\symVoltage_{\indexGridNodeTwo }^{\star}  }}

\newcommand{\VukrefSDP}[0]{\textcolor{\complexparamcolor}{\mathbf{\symVoltage}_{\indexUnit }^{\refss}  }}

\newcommand{\IurefSDP}[0]{\textcolor{\complexparamcolor}{\mathbf{\symCurrent}_{\indexGridNode }^{\refss}  }}

\newcommand{\VikrefSDP}[0]{\textcolor{\complexparamcolor}{\mathbf{\symVoltage}_{\indexGridNode }^{\refss}  }}
\newcommand{\VikrefSDPreal}[0]{\textcolor{\paramcolor}{\mathbf{\symVoltage}_{\indexGridNode }^{\refss,\realss}  }}
\newcommand{\VikrefSDPimag}[0]{\textcolor{\paramcolor}{\mathbf{\symVoltage}_{\indexGridNode }^{\refss,\imagss}  }}
\newcommand{\VikdeltaSDP}[0]{\textcolor{\complexcolor}{\mathbf{\symVoltage}_{\indexGridNode }^{\Delta}  }}
\newcommand{\VikdeltaSDPreal}[0]{\textcolor{black}{\mathbf{\symVoltage}_{\indexGridNode }^{\Delta,\realss}  }}
\newcommand{\VikdeltaSDPimag}[0]{\textcolor{black}{\mathbf{\symVoltage}_{\indexGridNode }^{\Delta,\imagss}  }}

\newcommand{\VsqikdeltaSDP}[0]{\textcolor{\complexcolor}{\mathbf{\symVoltageSOCP}_{\indexGridNode }^{\Delta}  }}

\newcommand{\VikaccSDP}[0]{\textcolor{\complexcolor}{\mathbf{\symVoltage}_{\indexGridNode }^{'}  }}
\newcommand{\VjkaccSDP}[0]{\textcolor{\complexcolor}{\mathbf{\symVoltage}_{\indexGridNodeTwo }^{'}  }}
\newcommand{\Cdelta}[0]{\textcolor{\complexcolor}{\mathbf{\symConnectivity}_{\Delta }^{}  }}
\newcommand{\Czigzag}[0]{\textcolor{\complexcolor}{\mathbf{\symConnectivity}_{Z }^{}  }}
\newcommand{\Cwye}[0]{\textcolor{\complexcolor}{\mathbf{\symConnectivity}_{Y }^{}  }}

\newcommand{\CijkSDP}[0]{\textcolor{\complexcolor}{\mathbf{\symConnectivity}_{\indexGridNode\indexGridNodeTwo }^{}  }}
\newcommand{\CjikSDP}[0]{\textcolor{\complexcolor}{\mathbf{\symConnectivity}_{\indexGridNodeTwo\indexGridNode }^{}  }}

\newcommand{\Visqk}[0]{\textcolor{black}{\symVoltage_{\indexGridNode }^{2}  }}
\newcommand{\Vjsqk}[0]{\textcolor{black}{\symVoltage_{\indexGridNodeTwo }^{2}  }}

\newcommand{\VikAN}[0]{\textcolor{\complexcolor}{\symVoltage_{\indexGridNode,\symPhaseA\symPhaseN }^{}  }}
\newcommand{\VikBN}[0]{\textcolor{\complexcolor}{\symVoltage_{\indexGridNode,\symPhaseB \symPhaseN}^{}  }}
\newcommand{\VikCN}[0]{\textcolor{\complexcolor}{\symVoltage_{\indexGridNode,\symPhaseC\symPhaseN }^{}  }}

\newcommand{\VikAB}[0]{\textcolor{\complexcolor}{\symVoltage_{\indexGridNode,\symPhaseA\symPhaseB }^{\Delta}  }}
\newcommand{\VikBC}[0]{\textcolor{\complexcolor}{\symVoltage_{\indexGridNode,\symPhaseB\symPhaseC }^{\Delta}  }}
\newcommand{\VikCA}[0]{\textcolor{\complexcolor}{\symVoltage_{\indexGridNode,\symPhaseC\symPhaseA }^{\Delta}  }}

\newcommand{\Vikzero}[0]{\textcolor{\complexcolor}{\symVoltage_{\indexGridNode,0 }^{}  }}
\newcommand{\Vikone}[0]{\textcolor{\complexcolor}{\symVoltage_{\indexGridNode,1 }^{}  }}
\newcommand{\Viktwo}[0]{\textcolor{\complexcolor}{\symVoltage_{\indexGridNode,2 }^{}  }}

\newcommand{\VikANreal}[0]{\textcolor{black}{\symVoltage_{\indexGridNode,\symPhaseA }^{\realss}  }}
\newcommand{\VikBNreal}[0]{\textcolor{black}{\symVoltage_{\indexGridNode,\symPhaseB }^{\realss}  }}
\newcommand{\VikCNreal}[0]{\textcolor{black}{\symVoltage_{\indexGridNode,\symPhaseC }^{\realss}  }}
\newcommand{\VikNreal}[0]{\textcolor{black}{\symVoltage_{\indexGridNode,\symPhaseN }^{\realss}  }}

\newcommand{\VikANimag}[0]{\textcolor{black}{\symVoltage_{\indexGridNode,\symPhaseA }^{\imagss}  }}
\newcommand{\VikBNimag}[0]{\textcolor{black}{\symVoltage_{\indexGridNode,\symPhaseB }^{\imagss}  }}
\newcommand{\VikCNimag}[0]{\textcolor{black}{\symVoltage_{\indexGridNode,\symPhaseC }^{\imagss}  }}
\newcommand{\VikNimag}[0]{\textcolor{black}{\symVoltage_{\indexGridNode,\symPhaseN }^{\imagss}  }}

\newcommand{\VikpNreal}[0]{\textcolor{black}{\symVoltage_{\indexGridNode,\indexPhases }^{\realss}  }}
\newcommand{\VikpNimag}[0]{\textcolor{black}{\symVoltage_{\indexGridNode,\indexPhases }^{\imagss}  }}
\newcommand{\VikppNreal}[0]{\textcolor{black}{\symVoltage_{\indexGridNode,\indexPhasesTwo }^{\realss}  }}
\newcommand{\VikppNimag}[0]{\textcolor{black}{\symVoltage_{\indexGridNode,\indexPhasesTwo }^{\imagss}  }}
\newcommand{\VjkppNreal}[0]{\textcolor{black}{\symVoltage_{\indexGridNodeTwo,\indexPhasesTwo }^{\realss}  }}
\newcommand{\VjkppNimag}[0]{\textcolor{black}{\symVoltage_{\indexGridNodeTwo,\indexPhasesTwo }^{\imagss}  }}

\newcommand{\VikpNmag}[0]{\textcolor{black}{\symVoltage_{\indexGridNode,\indexPhases }^{\text{mag}}  }}
\newcommand{\VikppNmag}[0]{\textcolor{black}{\symVoltage_{\indexGridNode,\indexPhasesTwo }^{\text{mag}}  }}
\newcommand{\VjkppNmag}[0]{\textcolor{black}{\symVoltage_{\indexGridNodeTwo,\indexPhasesTwo }^{\text{mag}}  }}

\newcommand{\VikANmag}[0]{\textcolor{black}{\symVoltage_{\indexGridNode,\symPhaseA }^{\text{mag}}  }}
\newcommand{\VikBNmag}[0]{\textcolor{black}{\symVoltage_{\indexGridNode,\symPhaseB }^{\text{mag}}  }}
\newcommand{\VikCNmag}[0]{\textcolor{black}{\symVoltage_{\indexGridNode,\symPhaseC }^{\text{mag}}  }}
\newcommand{\Vikmag}[0]{\textcolor{black}{\symVoltage_{\indexGridNode }^{\text{mag}}  }}

\newcommand{\VjkANmag}[0]{\textcolor{black}{\symVoltage_{\indexGridNodeTwo,\symPhaseA }^{\text{mag}}  }}
\newcommand{\VjkBNmag}[0]{\textcolor{black}{\symVoltage_{\indexGridNodeTwo,\symPhaseB }^{\text{mag}}  }}
\newcommand{\VjkCNmag}[0]{\textcolor{black}{\symVoltage_{\indexGridNodeTwo,\symPhaseC }^{\text{mag}}  }}

\newcommand{\VikANangle}[0]{\textcolor{black}{\symVoltageAngle_{\indexGridNode,\symPhaseA }^{}  }}
\newcommand{\VikBNangle}[0]{\textcolor{black}{\symVoltageAngle_{\indexGridNode,\symPhaseB }^{}  }}
\newcommand{\VikCNangle}[0]{\textcolor{black}{\symVoltageAngle_{\indexGridNode,\symPhaseC }^{}  }}
\newcommand{\VikNangle}[0]{\textcolor{black}{\symVoltageAngle_{\indexGridNode,\symPhaseN }^{}  }}

\newcommand{\VjkANangle}[0]{\textcolor{black}{\symVoltageAngle_{\indexGridNodeTwo,\symPhaseA }^{}  }}
\newcommand{\VjkBNangle}[0]{\textcolor{black}{\symVoltageAngle_{\indexGridNodeTwo,\symPhaseB }^{}  }}
\newcommand{\VjkCNangle}[0]{\textcolor{black}{\symVoltageAngle_{\indexGridNodeTwo,\symPhaseC }^{}  }}
\newcommand{\VjkNangle}[0]{\textcolor{black}{\symVoltageAngle_{\indexGridNodeTwo,\symPhaseN }^{}  }}

\newcommand{\VikpNangle}[0]{\textcolor{black}{\symVoltageAngle_{\indexGridNode,\indexPhases }^{}  }}
\newcommand{\VikhNangle}[0]{\textcolor{black}{\symVoltageAngle_{\indexGridNode,\indexPhasesTwo }^{}  }}
\newcommand{\VjkpNangle}[0]{\textcolor{black}{\symVoltageAngle_{\indexGridNodeTwo,\indexPhases }^{}  }}
\newcommand{\VjkhNangle}[0]{\textcolor{black}{\symVoltageAngle_{\indexGridNodeTwo,\indexPhasesTwo }^{}  }}

\newcommand{\VikrefAN}[0]{\textcolor{\complexparamcolor}{\symVoltage_{\indexGridNode,\symPhaseA }^{\refss}  }}
\newcommand{\VikrefBN}[0]{\textcolor{\complexparamcolor}{\symVoltage_{\indexGridNode,\symPhaseB }^{\refss}  }}
\newcommand{\VikrefCN}[0]{\textcolor{\complexparamcolor}{\symVoltage_{\indexGridNode,\symPhaseC  }^{\refss}  }}

\newcommand{\VikANref}[0]{\textcolor{\paramcolor}{\symVoltage_{\indexGridNode,\symPhaseA }^{\refss}  }}
\newcommand{\VikBNref}[0]{\textcolor{\paramcolor}{\symVoltage_{\indexGridNode,\symPhaseB }^{\refss}  }}
\newcommand{\VikCNref}[0]{\textcolor{\paramcolor}{\symVoltage_{\indexGridNode,\symPhaseC }^{\refss}  }}
\newcommand{\VikNref}[0]{\textcolor{\paramcolor}{\symVoltage_{\indexGridNode,n }^{\refss}  }}

\newcommand{\thetaikANref}[0]{\textcolor{\paramcolor}{\symVoltageAngle_{\indexGridNode,\symPhaseA }^{\refss}  }}
\newcommand{\thetaikBNref}[0]{\textcolor{\paramcolor}{\symVoltageAngle_{\indexGridNode,\symPhaseB }^{\refss}  }}
\newcommand{\thetaikCNref}[0]{\textcolor{\paramcolor}{\symVoltageAngle_{\indexGridNode,\symPhaseC }^{\refss}  }}
\newcommand{\thetaikNref}[0]{\textcolor{\paramcolor}{\symVoltageAngle_{\indexGridNode,n }^{\refss}  }}

\newcommand{\VikAmin}[0]{\textcolor{\paramcolor}{\symVoltage_{\indexGridNode,\symPhaseA }^{\minss}  }}
\newcommand{\VikBmin}[0]{\textcolor{\paramcolor}{\symVoltage_{\indexGridNode,\symPhaseB }^{\minss}  }}
\newcommand{\VikCmin}[0]{\textcolor{\paramcolor}{\symVoltage_{\indexGridNode,\symPhaseC }^{\minss}  }}
\newcommand{\VikAmax}[0]{\textcolor{\paramcolor}{\symVoltage_{\indexGridNode,\symPhaseA }^{\maxss}  }}
\newcommand{\VikBmax}[0]{\textcolor{\paramcolor}{\symVoltage_{\indexGridNode,\symPhaseB }^{\maxss}  }}
\newcommand{\VikCmax}[0]{\textcolor{\paramcolor}{\symVoltage_{\indexGridNode,\symPhaseC }^{\maxss}  }}

\newcommand{\VikANmin}[0]{\textcolor{\paramcolor}{\symVoltage_{\indexGridNode,\symPhaseA\symPhaseN }^{\minss}  }}
\newcommand{\VikBNmin}[0]{\textcolor{\paramcolor}{\symVoltage_{\indexGridNode,\symPhaseB\symPhaseN }^{\minss}  }}
\newcommand{\VikCNmin}[0]{\textcolor{\paramcolor}{\symVoltage_{\indexGridNode,\symPhaseC\symPhaseN }^{\minss}  }}
\newcommand{\VikNmin}[0]{\textcolor{\paramcolor}{\symVoltage_{\indexGridNode,\symPhaseN }^{\minss}  }}
\newcommand{\VikANmax}[0]{\textcolor{\paramcolor}{\symVoltage_{\indexGridNode,\symPhaseA\symPhaseN }^{\maxss}  }}
\newcommand{\VikBNmax}[0]{\textcolor{\paramcolor}{\symVoltage_{\indexGridNode,\symPhaseB\symPhaseN }^{\maxss}  }}
\newcommand{\VikCNmax}[0]{\textcolor{\paramcolor}{\symVoltage_{\indexGridNode,\symPhaseC\symPhaseN }^{\maxss}  }}
\newcommand{\VikNmax}[0]{\textcolor{\paramcolor}{\symVoltage_{\indexGridNode,\symPhaseN }^{\maxss}  }}

\newcommand{\VikABmax}[0]{\textcolor{\paramcolor}{\symVoltage_{\indexGridNode,\symPhaseA\symPhaseB }^{\maxss}  }}
\newcommand{\VikBCmax}[0]{\textcolor{\paramcolor}{\symVoltage_{\indexGridNode,\symPhaseB\symPhaseC }^{\maxss}  }}
\newcommand{\VikCAmax}[0]{\textcolor{\paramcolor}{\symVoltage_{\indexGridNode,\symPhaseC\symPhaseA }^{\maxss}  }}
\newcommand{\VikABmin}[0]{\textcolor{\paramcolor}{\symVoltage_{\indexGridNode,\symPhaseA\symPhaseB }^{\minss}  }}
\newcommand{\VikBCmin}[0]{\textcolor{\paramcolor}{\symVoltage_{\indexGridNode,\symPhaseB\symPhaseC }^{\minss}  }}
\newcommand{\VikCAmin}[0]{\textcolor{\paramcolor}{\symVoltage_{\indexGridNode,\symPhaseC\symPhaseA }^{\minss}  }}

\newcommand{\VjkAN}[0]{\textcolor{\complexcolor}{\symVoltage_{\indexGridNodeTwo,\symPhaseA }^{}  }}
\newcommand{\VjkBN}[0]{\textcolor{\complexcolor}{\symVoltage_{\indexGridNodeTwo,\symPhaseB }^{}  }}
\newcommand{\VjkCN}[0]{\textcolor{\complexcolor}{\symVoltage_{\indexGridNodeTwo,\symPhaseC }^{}  }}
\newcommand{\VjkANmax}[0]{\textcolor{\paramcolor}{\symVoltage_{\indexGridNodeTwo,\symPhaseA }^{\maxss}  }}

\newcommand{\VjkA}[0]{\textcolor{\complexcolor}{\symVoltage_{\indexGridNodeTwo,\symPhaseA }^{}  }}
\newcommand{\VjkB}[0]{\textcolor{\complexcolor}{\symVoltage_{\indexGridNodeTwo,\symPhaseB }^{}  }}
\newcommand{\VjkC}[0]{\textcolor{\complexcolor}{\symVoltage_{\indexGridNodeTwo,\symPhaseC }^{}  }}
\newcommand{\VjkN}[0]{\textcolor{\complexcolor}{\symVoltage_{\indexGridNodeTwo,\symPhaseN }^{}  }}
\newcommand{\VjkG}[0]{\textcolor{\complexcolor}{\symVoltage_{\indexGridNodeTwo,\symPhaseG }^{}  }}

\newcommand{\VikA}[0]{\textcolor{\complexcolor}{\symVoltage_{\indexGridNode,\symPhaseA }^{}  }}
\newcommand{\VikB}[0]{\textcolor{\complexcolor}{\symVoltage_{\indexGridNode,\symPhaseB }^{}  }}
\newcommand{\VikC}[0]{\textcolor{\complexcolor}{\symVoltage_{\indexGridNode,\symPhaseC }^{}  }}
\newcommand{\VikN}[0]{\textcolor{\complexcolor}{\symVoltage_{\indexGridNode,\symPhaseN }^{}  }}
\newcommand{\VikG}[0]{\textcolor{\complexcolor}{\symVoltage_{\indexGridNode,\symPhaseG }^{}  }}

\newcommand{\Vikmin}[0]{\textcolor{\boundscolor}{\symVoltage^{\minss}_{\indexGridNode }  }}
\newcommand{\Vjkmin}[0]{\textcolor{\boundscolor}{\symVoltage^{\minss}_{\indexGridNodeTwo }  }}
\newcommand{\Vikmax}[0]{\textcolor{\boundscolor}{\symVoltage^{\maxss}_{\indexGridNode }  }}
\newcommand{\Vjkmax}[0]{\textcolor{\boundscolor}{\symVoltage^{\maxss}_{\indexGridNodeTwo }  }}
\newcommand{\bigMikmax}[0]{\textcolor{\boundscolor}{M^{}_{\indexGridNode }  }}
\newcommand{\bigMjkmax}[0]{\textcolor{\boundscolor}{M^{}_{\indexGridNodeTwo }  }}
\newcommand{\Virated}[0]{\textcolor{\sizingcolor}{\symVoltage^{\ratedss}_{\indexGridNode } }}
\newcommand{\Vjrated}[0]{\textcolor{\sizingcolor}{\symVoltage^{\ratedss}_{\indexGridNodeTwo } }}
\newcommand{\Vijrated}[0]{\textcolor{\sizingcolor}{\symVoltage^{\ratedss}_{\indexGridNode \indexGridNodeTwo} }}
\newcommand{\Vjirated}[0]{\textcolor{\sizingcolor}{\symVoltage^{\ratedss}_{\indexGridNodeTwo \indexGridNode } }}

\newcommand{\ViratedAN}[0]{\textcolor{\sizingcolor}{\symVoltage^{\ratedss}_{\indexGridNode,\symPhaseA\symPhaseN } }}
\newcommand{\ViratedBN}[0]{\textcolor{\sizingcolor}{\symVoltage^{\ratedss}_{\indexGridNode,\symPhaseB\symPhaseN } }}
\newcommand{\ViratedCN}[0]{\textcolor{\sizingcolor}{\symVoltage^{\ratedss}_{\indexGridNode,\symPhaseC\symPhaseN } }}

\newcommand{\VijratedAN}[0]{\textcolor{\sizingcolor}{\symVoltage^{\ratedss}_{\indexGridNode \indexGridNodeTwo,\symPhaseA\symPhaseN} }}
\newcommand{\VjiratedAN}[0]{\textcolor{\sizingcolor}{\symVoltage^{\ratedss}_{\indexGridNodeTwo \indexGridNode,\symPhaseA\symPhaseN } }}
\newcommand{\VijratedBN}[0]{\textcolor{\sizingcolor}{\symVoltage^{\ratedss}_{\indexGridNode \indexGridNodeTwo,\symPhaseB\symPhaseN} }}
\newcommand{\VjiratedBN}[0]{\textcolor{\sizingcolor}{\symVoltage^{\ratedss}_{\indexGridNodeTwo \indexGridNode ,\symPhaseB\symPhaseN} }}
\newcommand{\VijratedCN}[0]{\textcolor{\sizingcolor}{\symVoltage^{\ratedss}_{\indexGridNode \indexGridNodeTwo,\symPhaseC\symPhaseN} }}
\newcommand{\VjiratedCN}[0]{\textcolor{\sizingcolor}{\symVoltage^{\ratedss}_{\indexGridNodeTwo \indexGridNode,\symPhaseC\symPhaseN } }}

\newcommand{\Visqksocplin}[0]{\textcolor{black}{\symVoltage_{\indexGridNode }^{\linss}  }}

\newcommand{\VisqkSDPseq}[0]{\textcolor{\complexcolor}{\mathbf{\symVoltageSOCP}_{\indexGridNode }^{\fortescuess} }}

\newcommand{\VusqkSDP}[0]{\textcolor{\complexcolor}{\mathbf{\symVoltageSOCP}_{\indexUnit }  }}

\newcommand{\VisqkSDPdelta}[0]{\textcolor{\complexcolor}{\mathbf{\symVoltageSOCP}_{\indexGridNode }^{\Delta}  }}
\newcommand{\VisqkSDPdeltareal}[0]{\textcolor{black}{\mathbf{\symVoltageSOCP}_{\indexGridNode }^{\Delta \realss}  }}
\newcommand{\VisqkSDPdeltaimag}[0]{\textcolor{black}{\mathbf{\symVoltageSOCP}_{\indexGridNode }^{\Delta \imagss}  }}

\newcommand{\VisqkSDP}[0]{\textcolor{\complexcolor}{\mathbf{\symVoltageSOCP}_{\indexGridNode }  }}
\newcommand{\VjsqkSDP}[0]{\textcolor{\complexcolor}{\mathbf{\symVoltageSOCP}_{\indexGridNodeTwo }  }}
\newcommand{\VisqkSDPreal}[0]{\textcolor{black}{\mathbf{\symVoltageSOCP}_{\indexGridNode }^{\realss}  }}
\newcommand{\VisqkSDPimag}[0]{\textcolor{black}{\mathbf{\symVoltageSOCP}_{\indexGridNode }^{\imagss}  }}
\newcommand{\VjsqkSDPreal}[0]{\textcolor{black}{\mathbf{\symVoltageSOCP}_{\indexGridNodeTwo }^{\realss}  }}
\newcommand{\VjsqkSDPimag}[0]{\textcolor{black}{\mathbf{\symVoltageSOCP}_{\indexGridNodeTwo }^{\imagss}  }}
\newcommand{\VisqkaccSDP}[0]{\textcolor{\complexcolor}{\mathbf{\symVoltageSOCP}_{\indexGridNode }^{'}  }}

\newcommand{\VzsqkSDP}[0]{\textcolor{\complexcolor}{\mathbf{\symVoltageSOCP}_{z }  }}
\newcommand{\VsqkSDP}[0]{\textcolor{\complexcolor}{\mathbf{M}  }}
\newcommand{\VsqkSDPtworeal}[0]{\textcolor{black}{\mathbf{M}^{2\realss}  }}
\newcommand{\VsqkSDPreal}[0]{\textcolor{black}{\mathbf{M}^{\realss}  }}
\newcommand{\VsqkSDPimag}[0]{\textcolor{black}{\mathbf{M}^{\imagss}  }}
\newcommand{\MisqkSDP}[0]{\textcolor{\complexcolor}{\mathbf{M}_i  }}
\newcommand{\MisqkSDPtworeal}[0]{\textcolor{black}{\mathbf{M}^{2\realss}_{\indexGridNode}  }}

\newcommand{\VbpsqkSDP}[0]{\textcolor{\complexcolor}{\mathbf{M}_{\indexGridNode\indexGridNodeTwo}  }}
\newcommand{\VbpsqkSDPtworeal}[0]{\textcolor{black}{\mathbf{M}^{2\realss}_{\indexGridNode\indexGridNodeTwo}  }}
\newcommand{\VbpsqkSDPreal}[0]{\textcolor{black}{\mathbf{M}_{\indexGridNode\indexGridNodeTwo} ^{\realss}  }}
\newcommand{\VbpsqkSDPimag}[0]{\textcolor{black}{\mathbf{M}_{\indexGridNode\indexGridNodeTwo} ^{\imagss}  }}

\newcommand{\VlinesqkSDP}[0]{\textcolor{\complexcolor}{\mathbf{M}_{\indexGridLines\indexGridNode\indexGridNodeTwo}  }}
\newcommand{\VlinesqkSDPreal}[0]{\textcolor{black}{\mathbf{M}_{\indexGridLines\indexGridNode\indexGridNodeTwo}^{\realss}  }}
\newcommand{\VlinesqkSDPrealtwo}[0]{\textcolor{black}{\mathbf{M}_{\indexGridLines\indexGridNode\indexGridNodeTwo}^{2\realss}  }}

\newcommand{\VlinetosqkSDP}[0]{\textcolor{\complexcolor}{\mathbf{M}_{\indexGridLines\indexGridNodeTwo\indexGridNode}  }}

\newcommand{\VunitsqkSDP}[0]{\textcolor{\complexcolor}{\mathbf{M}_{\indexGridNode\indexUnit}  }}
\newcommand{\VunitsqkSDPvoltage}[0]{\textcolor{\complexcolor}{\mathbf{M}^{\symVoltage}_{\indexGridNode\indexUnit}  }}
\newcommand{\VunitsqkSDPcurrent}[0]{\textcolor{\complexcolor}{\mathbf{M}_{\indexGridNode\indexUnit}  }}
\newcommand{\VunitsqkSDPcurrenttworeal}[0]{\textcolor{black}{\mathbf{M}^{2\realss}_{\indexGridNode\indexUnit}  }}

\newcommand{\VunitsqkSDPcurrentdelta}[0]{\textcolor{\complexcolor}{\mathbf{M}^{\Delta}_{\indexGridNode\indexUnit}  }}
\newcommand{\VunitsqkSDPcurrentdeltap}[0]{\textcolor{\complexcolor}{\mathbf{M}^{\Delta'}_{\indexGridNode\indexUnit}  }}

\newcommand{\VjzsqkSDP}[0]{\textcolor{\complexcolor}{\mathbf{\symVoltageSOCP}_{jz }  }}
\newcommand{\VizsqkSDP}[0]{\textcolor{\complexcolor}{\mathbf{\symVoltageSOCP}_{iz }  }}
\newcommand{\VzjsqkSDP}[0]{\textcolor{\complexcolor}{\mathbf{\symVoltageSOCP}_{zj }  }}
\newcommand{\VzisqkSDP}[0]{\textcolor{\complexcolor}{\mathbf{\symVoltageSOCP}_{zi }  }}

\newcommand{\ViusqkSDP}[0]{\textcolor{\complexcolor}{\mathbf{\symVoltageSOCP}_{\indexGridNode \indexUnit}  }}
\newcommand{\VuisqkSDP}[0]{\textcolor{\complexcolor}{\mathbf{\symVoltageSOCP}_{ \indexUnit \indexGridNode}  }}
\newcommand{\ViusqkSDPreal}[0]{\textcolor{black}{\mathbf{\symVoltageSOCP}_{\indexGridNode \indexUnit}^{\realss}  }}
\newcommand{\VuisqkSDPreal}[0]{\textcolor{black}{\mathbf{\symVoltageSOCP}_{ \indexUnit \indexGridNode}^{\realss}  }}
\newcommand{\ViusqkSDPimag}[0]{\textcolor{black}{\mathbf{\symVoltageSOCP}_{\indexGridNode \indexUnit}^{\imagss}   }}
\newcommand{\VuisqkSDPimag}[0]{\textcolor{black}{\mathbf{\symVoltageSOCP}_{ \indexUnit \indexGridNode}^{\imagss}   }}

\newcommand{\VijsqkSDP}[0]{\textcolor{\complexcolor}{\mathbf{\symVoltageSOCP}_{\indexGridNode \indexGridNodeTwo}  }}
\newcommand{\VjisqkSDP}[0]{\textcolor{\complexcolor}{\mathbf{\symVoltageSOCP}_{\indexGridNodeTwo \indexGridNode}  }}
\newcommand{\VijsqkSDPreal}[0]{\textcolor{black}{\mathbf{\symVoltageSOCP}_{\indexGridNode \indexGridNodeTwo}^{\realss}  }}
\newcommand{\VijsqkSDPimag}[0]{\textcolor{black}{\mathbf{\symVoltageSOCP}_{\indexGridNode \indexGridNodeTwo}^{\imagss}  }}

\newcommand{\ratioij}[0]{\textcolor{black}{\mathbf{\symRatio}_{\indexGridLines }  }}
\newcommand{\ratiosqij}[0]{\textcolor{black}{\mathbf{H}_{\indexGridLines }  }}
\newcommand{\ratioijA}[0]{\textcolor{\complexcolor}{{\symRatio}_{\indexGridLines , \symPhaseA}  }}
\newcommand{\ratioijB}[0]{\textcolor{\complexcolor}{{\symRatio}_{\indexGridLines , \symPhaseB}  }}
\newcommand{\ratioijC}[0]{\textcolor{\complexcolor}{{\symRatio}_{\indexGridLines , \symPhaseC}  }}

\newcommand{\ratiomijA}[0]{\textcolor{black}{{\symRatio}^{\text{mag}}_{\indexGridLines , \symPhaseA}  }}
\newcommand{\ratiomijB}[0]{\textcolor{black}{{\symRatio}^{\text{mag}}_{\indexGridLines , \symPhaseB}  }}
\newcommand{\ratiomijC}[0]{\textcolor{black}{{\symRatio}^{\text{mag}}_{\indexGridLines , \symPhaseC}  }}

\newcommand{\ratiomminijA}[0]{\textcolor{\paramcolor}{{\symRatio}^{\text{min}}_{\indexGridLines , \symPhaseA}  }}
\newcommand{\ratiomminijB}[0]{\textcolor{\paramcolor}{{\symRatio}^{\text{min}}_{\indexGridLines , \symPhaseB}  }}
\newcommand{\ratiomminijC}[0]{\textcolor{\paramcolor}{{\symRatio}^{\text{min}}_{\indexGridLines , \symPhaseC}  }}

\newcommand{\ratiommaxijA}[0]{\textcolor{\paramcolor}{{\symRatio}^{\text{max}}_{\indexGridLines , \symPhaseA}  }}
\newcommand{\ratiommaxijB}[0]{\textcolor{\paramcolor}{{\symRatio}^{\text{max}}_{\indexGridLines , \symPhaseB}  }}
\newcommand{\ratiommaxijC}[0]{\textcolor{\paramcolor}{{\symRatio}^{\text{max}}_{\indexGridLines , \symPhaseC}  }}

\newcommand{\ratioaijA}[0]{\textcolor{black}{{\symRatio}^{\angle}_{\indexGridLines\indexGridNode \indexGridNodeTwo, \symPhaseA}  }}
\newcommand{\ratioaijB}[0]{\textcolor{black}{{\symRatio}^{\angle}_{\indexGridLines\indexGridNode \indexGridNodeTwo, \symPhaseB}  }}
\newcommand{\ratioaijC}[0]{\textcolor{black}{{\symRatio}^{\angle}_{\indexGridLines\indexGridNode \indexGridNodeTwo, \symPhaseC}  }}

\newcommand{\ratioaijref}[0]{\textcolor{\paramcolor}{{\symRatio}^{\angle \text{ref}}_{\indexGridLines\indexGridNode \indexGridNodeTwo}  }}

\newcommand{\ratiomaxij}[0]{\textcolor{\paramcolor}{\mathbf{\symRatio}^{\text{max}}_{\indexGridLines }  }}
\newcommand{\ratiominij}[0]{\textcolor{\paramcolor}{\mathbf{\symRatio}^{\text{min}}_{\indexGridLines }  }}

\newcommand{\Visqksocp}[0]{\textcolor{black}{\symVoltageSOCP_{\indexGridNode }  }}
\newcommand{\Vjsqksocp}[0]{\textcolor{black}{\symVoltageSOCP_{\indexGridNodeTwo } }}
\newcommand{\Visqksocpacc}[0]{\textcolor{black}{\symVoltageSOCP_{\indexGridNode }^{'}  }}
\newcommand{\Vjsqksocpacc}[0]{\textcolor{black}{\symVoltageSOCP_{\indexGridNodeTwo }^{'} }}
\newcommand{\Visqksocpaccl}[0]{\textcolor{black}{\symVoltageSOCP_{\indexGridNode }^{\star}  }}
\newcommand{\Vjsqksocpaccl}[0]{\textcolor{black}{\symVoltageSOCP_{\indexGridNodeTwo }^{\star} }}
\newcommand{\Visqksocpjabr}[0]{\textcolor{black}{\symVoltageSOCP^{\star}_{\indexGridNode }  }}
\newcommand{\Vjsqksocpjabr}[0]{\textcolor{black}{\symVoltageSOCP^{\star}_{\indexGridNodeTwo } }}
\newcommand{\Visqksocpjabrl}[0]{\textcolor{black}{\symVoltageSOCP^{\star\indexGridLines}_{\indexGridNode }  }}
\newcommand{\Vjsqksocpjabrl}[0]{\textcolor{black}{\symVoltageSOCP^{\star\indexGridLines}_{\indexGridNodeTwo } }}

\newcommand{\Ajabr}[0]{\textcolor{\paramcolor}{A^{'}_{\indexGridNode \indexGridNodeTwo } }}
\newcommand{\Bjabr}[0]{\textcolor{\paramcolor}{B^{'}_{\indexGridNode \indexGridNodeTwo } }}
\newcommand{\Cjabr}[0]{\textcolor{\paramcolor}{C^{'}_{\indexGridNode \indexGridNodeTwo } }}
\newcommand{\Djabr}[0]{\textcolor{\paramcolor}{D^{'}_{\indexGridNode \indexGridNodeTwo } }}

\newcommand{\bigM}[0]{\textcolor{\paramcolor}{M^{} }}

\newcommand{\thetaikref}[0]{\textcolor{\paramcolor}{\symVoltageAngle_{\indexGridNode }^{\refss}  }}
\newcommand{\thetaikrefone}[0]{\textcolor{\paramcolor}{\symVoltageAngle_{\indexGridNode,1 }^{\refss}  }}
\newcommand{\Vikref}[0]{\textcolor{\paramcolor}{\symVoltage_{\indexGridNode }^{\refss}  }}
\newcommand{\Viref}[0]{\textcolor{\paramcolor}{\symVoltage_{\indexGridNode }^{\refss} }}

\newcommand{\Vjkref}[0]{\textcolor{\paramcolor}{\symVoltage_{\indexGridNodeTwo }^{\refss}  }}
\newcommand{\Vjref}[0]{\textcolor{\paramcolor}{\symVoltage_{\indexGridNodeTwo }^{\refss} }}

\newcommand{\thetaik}[0]{\textcolor{black}{\symVoltageAngle_{\indexGridNode }^{}  }}
\newcommand{\thetajk}[0]{\textcolor{black}{\symVoltageAngle_{\indexGridNodeTwo }^{}  }}
\newcommand{\thetaijk}[0]{\textcolor{black}{\symVoltageAngle_{\indexGridNode\indexGridNodeTwo }^{}  }}
\newcommand{\thetaijacckmin}[0]{\textcolor{\paramcolor}{\symVoltageAngle_{\indexGridNode\indexGridNodeTwo }^{'\minss}  }}
\newcommand{\thetaijacckmax}[0]{\textcolor{\paramcolor}{\symVoltageAngle_{\indexGridNode\indexGridNodeTwo }^{'\maxss}  }}
\newcommand{\thetaikmin}[0]{\textcolor{\boundscolor}{\symVoltageAngle_{\indexGridNode }^{\minss}  }}
\newcommand{\thetaikmax}[0]{\textcolor{\boundscolor}{\symVoltageAngle_{\indexGridNode }^{\maxss}  }}

\newcommand{\thetaijkmax}[0]{\textcolor{\boundscolor}{\symVoltageAngle_{\indexGridNode \indexGridNodeTwo}^{\maxss}  }}
\newcommand{\thetaijkmin}[0]{\textcolor{\boundscolor}{\symVoltageAngle_{\indexGridNode \indexGridNodeTwo}^{\minss}  }}

\newcommand{\thetaijkabs}[0]{\textcolor{\boundscolor}{\symVoltageAngle_{\indexGridNode \indexGridNodeTwo}^{\text{absmax}}  }}

\newcommand{\thetaijkmaxacc}[0]{\textcolor{\boundscolor}{\symVoltageAngle_{\indexGridNode \indexGridNodeTwo}^{'\maxss}  }}
\newcommand{\thetaijkminacc}[0]{\textcolor{\boundscolor}{\symVoltageAngle_{\indexGridNode \indexGridNodeTwo}^{'\minss}  }}

\newcommand{\thetaijkabsacc}[0]{\textcolor{\boundscolor}{\symVoltageAngle_{\indexGridNode \indexGridNodeTwo}^{'\text{absmax}}  }}

\newcommand{\phiijk}[0]{\textcolor{black}{\symPhaseDifference_{\indexGridNode\indexGridNodeTwo }^{}  }}
\newcommand{\phiikmax}[0]{\textcolor{\boundscolor}{\symPhaseDifference_{\indexGridNode }^{\maxss}  }}
\newcommand{\phiikmin}[0]{\textcolor{\boundscolor}{\symPhaseDifference_{\indexGridNode }^{\minss}  }}
\newcommand{\phiijkmin}[0]{\textcolor{\boundscolor}{\symPhaseDifference_{\indexGridNode   \indexGridNodeTwo }^{\minss}  }}
\newcommand{\phiijkmax}[0]{\textcolor{\boundscolor}{\symPhaseDifference_{\indexGridNode   \indexGridNodeTwo }^{\maxss}  }}

\newcommand{\phijik}[0]{\textcolor{black}{\symPhaseDifference_{\indexGridNodeTwo   \indexGridNode }^{}  }}
\newcommand{\phijikmin}[0]{\textcolor{\boundscolor}{\symPhaseDifference_{\indexGridNodeTwo   \indexGridNode    }^{\minss}  }}
\newcommand{\phijikmax}[0]{\textcolor{\boundscolor}{\symPhaseDifference_{\indexGridNodeTwo  \indexGridNode    }^{\maxss}  }}

\newcommand{\thetaiacck}[0]{\textcolor{black}{\symVoltageAngle_{\indexGridNode }^{'}  }}
\newcommand{\thetajacck}[0]{\textcolor{black}{\symVoltageAngle_{\indexGridNodeTwo }^{'}  }}
\newcommand{\thetaijacck}[0]{\textcolor{black}{\symVoltageAngle_{\indexGridNode  \indexGridNodeTwo }^{'}  }}

\newcommand{\thetaiacckl}[0]{\textcolor{black}{\symVoltageAngle_{\indexGridNode }^{\star}  }}
\newcommand{\thetajacckl}[0]{\textcolor{black}{\symVoltageAngle_{\indexGridNodeTwo }^{\star}  }}

\newcommand{\yijkshSDPseq}[0]{\textcolor{\complexparamcolor}{\mathbf{\symAdmittance}_{\indexGridLines \indexGridNode   \indexGridNodeTwo}^{\shuntss, \fortescuess}    }}

\newcommand{\yijkshSDPAdot}[0]{\textcolor{\complexparamcolor}{\mathbf{\symAdmittance}_{\indexGridLines \indexGridNode   \indexGridNodeTwo}^{\shuntss, \symPhaseA\!\cdot}    }}
\newcommand{\yijkshSDPBdot}[0]{\textcolor{\complexparamcolor}{\mathbf{\symAdmittance}_{\indexGridLines \indexGridNode   \indexGridNodeTwo}^{\shuntss, \symPhaseB\!\cdot}    }}
\newcommand{\yijkshSDPCdot}[0]{\textcolor{\complexparamcolor}{\mathbf{\symAdmittance}_{\indexGridLines \indexGridNode   \indexGridNodeTwo}^{\shuntss, \symPhaseC\!\cdot}    }}
\newcommand{\yijkshSDPNdot}[0]{\textcolor{\complexparamcolor}{\mathbf{\symAdmittance}_{\indexGridLines \indexGridNode   \indexGridNodeTwo}^{\shuntss, \symPhaseN\!\cdot}    }}

\newcommand{\yijkshSDPdotN}[0]{\textcolor{\complexparamcolor}{\mathbf{\symAdmittance}_{\indexGridLines \indexGridNode   \indexGridNodeTwo}^{\shuntss, \cdot \! \symPhaseN}    }}

\newcommand{\yijkshSDP}[0]{\textcolor{\complexparamcolor}{\mathbf{\symAdmittance}_{\indexGridLines \indexGridNode   \indexGridNodeTwo}^{\shuntss}    }}
\newcommand{\gijkshSDP}[0]{\textcolor{\paramcolor}{\mathbf{\symConductance}_{ \indexGridLines\indexGridNode   \indexGridNodeTwo}^{\shuntss}    }}
\newcommand{\bijkshSDP}[0]{\textcolor{\paramcolor}{\mathbf{\symSusceptance}_{ \indexGridLines\indexGridNode   \indexGridNodeTwo}^{\shuntss}    }}
\newcommand{\zijkshSDP}[0]{\textcolor{\complexparamcolor}{\mathbf{\symImpedance}_{\indexGridLines \indexGridNode   \indexGridNodeTwo}^{\shuntss}    }}

\newcommand{\yuunitSDP}[0]{\textcolor{\complexparamcolor}{\mathbf{\symAdmittance}_{    \indexUnit}  }}
\newcommand{\yuunitSDPDelta}[0]{\textcolor{\complexparamcolor}{\mathbf{\symAdmittance}_{    \indexUnit}^{\Delta}  }}

\newcommand{\ybSDP}[0]{\textcolor{\complexparamcolor}{\mathbf{\symAdmittance}_{    \indexShunt}  }}
\newcommand{\gbSDP}[0]{\textcolor{\paramcolor}{\mathbf{\symConductance}_{    \indexShunt}  }}
\newcommand{\bbSDP}[0]{\textcolor{\paramcolor}{\mathbf{\symSusceptance}_{    \indexShunt}  }}

\newcommand{\yjikshSDP}[0]{\textcolor{\complexparamcolor}{\mathbf{\symAdmittance}_{ \indexGridLines\indexGridNodeTwo   \indexGridNode}^{\shuntss}    }}
\newcommand{\gjikshSDP}[0]{\textcolor{\paramcolor}{\mathbf{\symConductance}_{ \indexGridLines\indexGridNodeTwo   \indexGridNode}^{\shuntss}    }}
\newcommand{\bjikshSDP}[0]{\textcolor{\paramcolor}{\mathbf{\symSusceptance}_{ \indexGridLines\indexGridNodeTwo   \indexGridNode}^{\shuntss}    }}
\newcommand{\zjikshSDP}[0]{\textcolor{\complexparamcolor}{\mathbf{\symImpedance}_{\indexGridLines \indexGridNodeTwo   \indexGridNode}^{\shuntss}    }}

\newcommand{\yikSDP}[0]{\textcolor{\complexparamcolor}{\mathbf{\symAdmittance}_{ \indexGridNode   }^{}    }}
\newcommand{\gikSDP}[0]{\textcolor{\paramcolor}{\mathbf{\symConductance}_{ \indexGridNode   }^{}    }}
\newcommand{\bikSDP}[0]{\textcolor{\paramcolor}{\mathbf{\symSusceptance}_{ \indexGridNode   }^{}    }}

\newcommand{\zijksSDPkron}[0]{\textcolor{\complexparamcolor}{\mathbf{{\symImpedance}}_{ \indexGridLines}^{\seriesss, \text{Kron}}    }}
\newcommand{\zijksSDPphase}[0]{\textcolor{\complexparamcolor}{\mathbf{{\symImpedance}}_{ \indexGridLines}^{\seriesss, \text{phase}}    }}
\newcommand{\zijksSDPprim}[0]{\textcolor{\complexparamcolor}{\mathbf{\hat{\symImpedance}}_{ \indexGridLines}^{\seriesss, \text{circ}}    }}
\newcommand{\zijksSDPvec}[0]{\textcolor{\complexparamcolor}{\mathbf{\hat{\symImpedance}}_{ \indexGridLines}^{\seriesss, \text{pn}}    }}
\newcommand{\zijksSDPvect}[0]{\textcolor{\complexparamcolor}{\mathbf{\hat{\symImpedance}}_{ \indexGridLines}^{\seriesss,  \text{np}}    }}
\newcommand{\zijksSDPhat}[0]{\textcolor{\complexparamcolor}{\mathbf{\hat{\symImpedance}}_{ \indexGridLines}^{\seriesss,  \text{pp}}    }}

\newcommand{\zijksSDPcond}[0]{\textcolor{\complexparamcolor}{\mathbf{\bar{\symImpedance}}_{ \indexGridLines}^{\seriesss, \text{cond}}    }}
\newcommand{\zijksSDPseq}[0]{\textcolor{\complexparamcolor}{\mathbf{\symImpedance}_{ \indexGridLines}^{\seriesss, \fortescuess}    }}

\newcommand{\zijksSDPseqdiag}[0]{\textcolor{\complexparamcolor}{\mathbf{\symImpedance}_{ \indexGridLines}^{\seriesss, \fortescuess, \text{diag}}    }}

\newcommand{\seq}[0]{\textcolor{\complexparamcolor}{\mathbf{A}}}
\newcommand{\seqinv}[0]{\textcolor{\complexparamcolor}{\mathbf{A}^{-1}}}
\newcommand{\seqH}[0]{\textcolor{\complexparamcolor}{\mathbf{A}^{\hermitiantranspose}}}
\newcommand{\seqconj}[0]{\textcolor{\complexparamcolor}{\mathbf{A}^{*}}}
\newcommand{\seqT}[0]{\textcolor{\complexparamcolor}{\mathbf{A}^{\transpose}}}

\newcommand{\seqvec}[0]{\textcolor{\complexparamcolor}{\mathbf{a}}}

\newcommand{\zijksSDP}[0]{\textcolor{\complexparamcolor}{\mathbf{\symImpedance}_{ \indexGridLines}^{\seriesss}    }}
\newcommand{\rijksSDP}[0]{\textcolor{\paramcolor}{\mathbf{\symResistance}_{  \indexGridLines   }^{\seriesss}    }}
\newcommand{\xijksSDP}[0]{\textcolor{\paramcolor}{\mathbf{\symReactance}_{  \indexGridLines   }^{\seriesss}    }}

\newcommand{\zjiksSDP}[0]{\textcolor{\complexparamcolor}{\mathbf{\symImpedance}_{ \indexGridLines   }^{\seriesss}    }}
\newcommand{\rjiksSDP}[0]{\textcolor{\paramcolor}{\mathbf{\symResistance}_{ \indexGridLines   }^{\seriesss}    }}
\newcommand{\xjiksSDP}[0]{\textcolor{\paramcolor}{\mathbf{\symReactance}_{ \indexGridLines   }^{\seriesss}    }}

\newcommand{\zijks}[0]{\textcolor{\complexparamcolor}{\symImpedance_{ \indexGridLines,\seriesss}    }}
\newcommand{\rijks}[0]{\textcolor{\paramcolor}{\symResistance_{  \indexGridLines,\seriesss}    }}
\newcommand{\xijks}[0]{\textcolor{\paramcolor}{\symReactance_{  \indexGridLines,\seriesss}    }}
\newcommand{\yijks}[0]{\textcolor{\complexparamcolor}{\symAdmittance_{  \indexGridLines}^{\seriesss}    }}
\newcommand{\bijks}[0]{\textcolor{\paramcolor}{\symSusceptance_{  \indexGridLines}^{\seriesss}    }}
\newcommand{\gijks}[0]{\textcolor{\paramcolor}{\symConductance_{  \indexGridLines}^{\seriesss}    }}

\newcommand{\zijksSDPH}[0]{\textcolor{\complexparamcolor}{(\mathbf{\symImpedance}_{ \indexGridLines}^{\seriesss})^{\hermitiantranspose}    }}
\newcommand{\rijksSDPH}[0]{\textcolor{\paramcolor}{(\mathbf{\symResistance}_{ \indexGridLines,\seriesss})^{\transpose}    }}
\newcommand{\xijksSDPH}[0]{\textcolor{\paramcolor}{(\mathbf{\symReactance}_{ \indexGridLines,\seriesss})^{\transpose}    }}

\newcommand{\yukSDP}[0]{\textcolor{\complexparamcolor}{\mathbf{\symAdmittance}_{ \indexUnit}    }}
\newcommand{\gukSDP}[0]{\textcolor{\paramcolor}{\mathbf{\symConductance}_{ \indexUnit}    }}
\newcommand{\bukSDP}[0]{\textcolor{\paramcolor}{\mathbf{\symSusceptance}_{ \indexUnit}    }}

\newcommand{\zukSDP}[0]{\textcolor{\complexparamcolor}{\mathbf{\symImpedance}_{ \indexUnit}    }}
\newcommand{\rukSDP}[0]{\textcolor{\paramcolor}{\mathbf{\symResistance}_{ \indexUnit}    }}
\newcommand{\xukSDP}[0]{\textcolor{\paramcolor}{\mathbf{\symReactance}_{ \indexUnit}    }}

\newcommand{\yijksSDP}[0]{\textcolor{\complexparamcolor}{\mathbf{\symAdmittance}_{ \indexGridLines}^{\seriesss}    }}
\newcommand{\gijksSDP}[0]{\textcolor{\paramcolor}{\mathbf{\symConductance}_{ \indexGridLines}^{\seriesss}    }}
\newcommand{\bijksSDP}[0]{\textcolor{\paramcolor}{\mathbf{\symSusceptance}_{ \indexGridLines}^{\seriesss}    }}

\newcommand{\yjiksSDP}[0]{\textcolor{\complexparamcolor}{\mathbf{\symAdmittance}_{  \indexGridLines\indexGridNodeTwo \indexGridNode }^{\seriesss}    }}
\newcommand{\gjiksSDP}[0]{\textcolor{\paramcolor}{\mathbf{\symConductance}_{ \indexGridLines \indexGridNodeTwo \indexGridNode }^{\seriesss}    }}
\newcommand{\bjiksSDP}[0]{\textcolor{\paramcolor}{\mathbf{\symSusceptance}_{  \indexGridLines\indexGridNodeTwo \indexGridNode }^{\seriesss}    }}

\newcommand{\zijksh}[0]{\textcolor{\complexparamcolor}{\symImpedance_{ \indexGridNode \indexGridNodeTwo,\shuntss}    }}
\newcommand{\rijksh}[0]{\textcolor{\paramcolor}{\symResistance_{ \indexGridNode \indexGridNodeTwo,\shuntss}    }}
\newcommand{\xijksh}[0]{\textcolor{\paramcolor}{\symReactance_{ \indexGridNode \indexGridNodeTwo,\shuntss}    }}
\newcommand{\yijksh}[0]{\textcolor{\complexparamcolor}{\symAdmittance_{ \indexGridNode \indexGridNodeTwo,\shuntss}    }}
\newcommand{\gijksh}[0]{\textcolor{\paramcolor}{\symConductance_{ \indexGridNode \indexGridNodeTwo,\shuntss}    }}
\newcommand{\bijksh}[0]{\textcolor{\paramcolor}{\symSusceptance_{ \indexGridNode \indexGridNodeTwo,\shuntss}    }}

\newcommand{\zjiksh}[0]{\textcolor{\complexparamcolor}{\symImpedance_{ \indexGridNodeTwo \indexGridNode ,\shuntss}    }}
\newcommand{\rjiksh}[0]{\textcolor{\paramcolor}{\symResistance_{\indexGridNodeTwo \indexGridNode ,\shuntss}    }}
\newcommand{\xjiksh}[0]{\textcolor{\paramcolor}{\symReactance_{ \indexGridNodeTwo\indexGridNode ,\shuntss}    }}
\newcommand{\yjiksh}[0]{\textcolor{\complexparamcolor}{\symAdmittance_{\indexGridNodeTwo \indexGridNode ,\shuntss}    }}
\newcommand{\gjiksh}[0]{\textcolor{\paramcolor}{\symConductance_{\indexGridNodeTwo \indexGridNode ,\shuntss}    }}
\newcommand{\bjiksh}[0]{\textcolor{\paramcolor}{\symSusceptance_{ \indexGridNodeTwo \indexGridNode ,\shuntss}    }}

\newcommand{\zlksh}[0]{\textcolor{\complexparamcolor}{\symImpedance_{  \indexGridLines,\shuntss}    }}
\newcommand{\rlksh}[0]{\textcolor{\paramcolor}{\symResistance_{  \indexGridLines,\shuntss}    }}
\newcommand{\xlksh}[0]{\textcolor{\paramcolor}{\symReactance_{  \indexGridLines,\shuntss}    }}
\newcommand{\ylksh}[0]{\textcolor{\complexparamcolor}{\symAdmittance_{  \indexGridLines,\shuntss}    }}
\newcommand{\glksh}[0]{\textcolor{\paramcolor}{\symConductance_{  \indexGridLines,\shuntss}    }}
\newcommand{\blksh}[0]{\textcolor{\paramcolor}{\symSusceptance_{  \indexGridLines,\shuntss}    }}

\newcommand{\zlkshspec}[0]{\textcolor{\complexparamcolor}{\dot{\symImpedance}_{  \indexGridLines,\shuntss}    }}
\newcommand{\rlkshspec}[0]{\textcolor{\paramcolor}{\dot{\symResistance}_{  \indexGridLines,\shuntss}    }}
\newcommand{\xlkshv}[0]{\textcolor{\paramcolor}{\dot{\symReactance}_{  \indexGridLines,\shuntss}    }}
\newcommand{\ylkshspec}[0]{\textcolor{\complexparamcolor}{\dot{\symAdmittance}_{  \indexGridLines,\shuntss}    }}
\newcommand{\glkshspec}[0]{\textcolor{\paramcolor}{\dot{\symConductance}_{  \indexGridLines,\shuntss}    }}
\newcommand{\blkshspec}[0]{\textcolor{\paramcolor}{\dot{\symSusceptance}_{  \indexGridLines,\shuntss}    }}

\newcommand{\zijksspec}[0]{\textcolor{\complexparamcolor}{\dot{\symImpedance}_{  \indexGridLines,\seriesss}    }}
\newcommand{\rijksspec}[0]{\textcolor{\paramcolor}{\dot{\symResistance}_{  \indexGridLines,\seriesss}    }}
\newcommand{\xijksspec}[0]{\textcolor{\paramcolor}{\dot{\symReactance}_{  \indexGridLines,\seriesss}    }}
\newcommand{\yijkshspec}[0]{\textcolor{\complexparamcolor}{\dot{\symAdmittance}_{  \indexGridNode \indexGridNodeTwo,\shuntss}    }}
\newcommand{\bijkshspec}[0]{\textcolor{\paramcolor}{\dot{\symSusceptance}_{  \indexGridNode \indexGridNodeTwo,\shuntss}    }}
\newcommand{\gijkshspec}[0]{\textcolor{\paramcolor}{\dot{\symConductance}_{  \indexGridNode \indexGridNodeTwo,\shuntss}    }}

\newcommand{\lijk}[0]{\textcolor{\paramcolor}{{\symLength}_{  \indexGridLines}    }}
\newcommand{\vprimvsecijk}[0]{\textcolor{black}{{\symRatio}_{ \indexGridNode \indexGridNodeTwo}    }}
\newcommand{\vprimvsecijkmin}[0]{\textcolor{\boundscolor}{{\symRatio}_{ \indexGridNode \indexGridNodeTwo}^{\minss}    }}
\newcommand{\vprimvsecijkmax}[0]{\textcolor{\boundscolor}{{\symRatio}_{ \indexGridNode \indexGridNodeTwo}^{\maxss}    }}

\newcommand{\vsecvprimijk}[0]{\textcolor{black}{{\symRatio}_{\indexGridNodeTwo \indexGridNode }    }}
\newcommand{\vsecvprimijkmin}[0]{\textcolor{\boundscolor}{{\symRatio}_{ \indexGridNodeTwo\indexGridNode }^{\minss}    }}
\newcommand{\vsecvprimijkmax}[0]{\textcolor{\boundscolor}{{\symRatio}_{\indexGridNodeTwo \indexGridNode }^{\maxss}    }}

\newcommand{\YPFconvexification}[0]{\textcolor{black}{{\symPenalty}^{\PFconvexss}   }}
\newcommand{\YGIconvexification}[0]{\textcolor{black}{{\symPenalty}^{\GIconvexss}   }}
\newcommand{\YPgridloss}[0]{\textcolor{black}{{\symPenalty}^{\symPower,\lossss}   }}
\newcommand{\YSgridloss}[0]{\textcolor{black}{{\symPenalty}^{\symApparentPower,\lossss}   }}
\newcommand{\Ytotcurrent}[0]{\textcolor{black}{{\symPenalty}^{\text{\symCurrent}}   }}

\newcommand{\Yobj}[0]{\textcolor{black}{{\symPenalty}^{\text{cost}}   }}
\newcommand{\Kobj}[0]{\textcolor{black}{{\symAnnualCost}^{\text{cost}}   }}

\newcommand{\NPFconvexification}[0]{\textcolor{\paramcolor}{{\symPenaltyWeight}_{\PFconvexss}   }}
\newcommand{\NGIconvexification}[0]{\textcolor{\paramcolor}{{\symPenaltyWeight}_{\GIconvexss}   }}

\newcommand{\wPFconvexification}[0]{\textcolor{\paramcolor}{{\symObjectiveWeight}_{\PFconvexss}   }}
\newcommand{\wPFconvexificationopt}[0]{\textcolor{\paramcolor}{{\symObjectiveWeight}_{\PFconvexss}^{\text{opt}}   }}

\newcommand{\zijkszero}[0]{\textcolor{\complexparamcolor}{\symImpedance_{ \indexGridLines , 0}^{\seriesss}    }}
\newcommand{\zijksone}[0]{\textcolor{\complexparamcolor}{\symImpedance_{ \indexGridLines , 1}^{\seriesss}    }}

\newcommand{\zijkssymmzerozero}[0]{\textcolor{\complexparamcolor}{\symImpedance_{ \indexGridLines , 0 0}^{\seriesss}    }}
\newcommand{\zijkssymmzeroone}[0]{\textcolor{\complexparamcolor}{\symImpedance_{ \indexGridLines , 0 1}^{\seriesss}    }}
\newcommand{\zijkssymmzerotwo}[0]{\textcolor{\complexparamcolor}{\symImpedance_{ \indexGridLines , 0 2}^{\seriesss}    }}

\newcommand{\zijkssymmonezero}[0]{\textcolor{\complexparamcolor}{\symImpedance_{ \indexGridLines , 1 0}^{\seriesss}    }}
\newcommand{\zijkssymmoneone}[0]{\textcolor{\complexparamcolor}{\symImpedance_{ \indexGridLines , 1 1}^{\seriesss}    }}
\newcommand{\zijkssymmonetwo}[0]{\textcolor{\complexparamcolor}{\symImpedance_{ \indexGridLines , 1 2}^{\seriesss}    }}

\newcommand{\zijkssymmtwozero}[0]{\textcolor{\complexparamcolor}{\symImpedance_{ \indexGridLines , 2 0}^{\seriesss}    }}
\newcommand{\zijkssymmtwoone}[0]{\textcolor{\complexparamcolor}{\symImpedance_{ \indexGridLines , 2 1}^{\seriesss}    }}
\newcommand{\zijkssymmtwotwo}[0]{\textcolor{\complexparamcolor}{\symImpedance_{ \indexGridLines , 2 2}^{\seriesss}    }}

\newcommand{\zjNG}[0]{\textcolor{\complexparamcolor}{\symImpedance_{ \indexGridNodeTwo , \symPhaseN \symPhaseG}   }}
\newcommand{\ziNG}[0]{\textcolor{\complexparamcolor}{\symImpedance_{ \indexGridNode , \symPhaseN \symPhaseG}   }}
\newcommand{\riNG}[0]{\textcolor{\paramcolor}{\symResistance_{ \indexGridNode , \symPhaseN \symPhaseG}   }}
\newcommand{\xiNG}[0]{\textcolor{\paramcolor}{\symReactance_{ \indexGridNode , \symPhaseN \symPhaseG}   }}

\newcommand{\zijksAA}[0]{\textcolor{\complexparamcolor}{\symImpedance_{ \indexGridLines , \symPhaseA \symPhaseA}^{\seriesss}    }}
\newcommand{\zijksAB}[0]{\textcolor{\complexparamcolor}{\symImpedance_{ \indexGridLines , \symPhaseA \symPhaseB}^{\seriesss}    }}
\newcommand{\zijksAC}[0]{\textcolor{\complexparamcolor}{\symImpedance_{ \indexGridLines , \symPhaseA \symPhaseC}^{\seriesss}    }}
\newcommand{\zijksBA}[0]{\textcolor{\complexparamcolor}{\symImpedance_{ \indexGridLines , \symPhaseB \symPhaseA}^{\seriesss}    }}
\newcommand{\zijksBB}[0]{\textcolor{\complexparamcolor}{\symImpedance_{ \indexGridLines , \symPhaseB \symPhaseB}^{\seriesss}    }}
\newcommand{\zijksBC}[0]{\textcolor{\complexparamcolor}{\symImpedance_{ \indexGridLines , \symPhaseB \symPhaseC}^{\seriesss}    }}
\newcommand{\zijksCA}[0]{\textcolor{\complexparamcolor}{\symImpedance_{ \indexGridLines , \symPhaseC \symPhaseA}^{\seriesss}    }}
\newcommand{\zijksCB}[0]{\textcolor{\complexparamcolor}{\symImpedance_{ \indexGridLines , \symPhaseC \symPhaseB}^{\seriesss}    }}
\newcommand{\zijksCC}[0]{\textcolor{\complexparamcolor}{\symImpedance_{ \indexGridLines , \symPhaseC \symPhaseC}^{\seriesss}    }}

\newcommand{\zijksAN}[0]{\textcolor{\complexparamcolor}{\symImpedance_{ \indexGridLines , \symPhaseA \symPhaseN}^{\seriesss}    }}
\newcommand{\zijksBN}[0]{\textcolor{\complexparamcolor}{\symImpedance_{ \indexGridLines , \symPhaseB \symPhaseN}^{\seriesss}    }}
\newcommand{\zijksCN}[0]{\textcolor{\complexparamcolor}{\symImpedance_{ \indexGridLines , \symPhaseC \symPhaseN}^{\seriesss}    }}
\newcommand{\zijksNA}[0]{\textcolor{\complexparamcolor}{\symImpedance_{ \indexGridLines , \symPhaseN \symPhaseA }^{\seriesss}    }}
\newcommand{\zijksNB}[0]{\textcolor{\complexparamcolor}{\symImpedance_{ \indexGridLines ,\symPhaseN \symPhaseB }^{\seriesss}    }}
\newcommand{\zijksNC}[0]{\textcolor{\complexparamcolor}{\symImpedance_{ \indexGridLines , \symPhaseN\symPhaseC }^{\seriesss}    }}

\newcommand{\zijksNN}[0]{\textcolor{\complexparamcolor}{\symImpedance_{ \indexGridLines , \symPhaseN \symPhaseN}^{\seriesss}    }}

\newcommand{\zijksAG}[0]{\textcolor{\complexparamcolor}{\symImpedance_{ \indexGridLines , \symPhaseA \symPhaseG}^{\seriesss}    }}
\newcommand{\zijksBG}[0]{\textcolor{\complexparamcolor}{\symImpedance_{ \indexGridLines , \symPhaseB \symPhaseG}^{\seriesss}    }}
\newcommand{\zijksCG}[0]{\textcolor{\complexparamcolor}{\symImpedance_{ \indexGridLines , \symPhaseC \symPhaseG}^{\seriesss}    }}
\newcommand{\zijksNG}[0]{\textcolor{\complexparamcolor}{\symImpedance_{ \indexGridLines , \symPhaseN \symPhaseG}^{\seriesss}    }}
\newcommand{\zijksGA}[0]{\textcolor{\complexparamcolor}{\symImpedance_{ \indexGridLines , \symPhaseG \symPhaseA }^{\seriesss}    }}
\newcommand{\zijksGB}[0]{\textcolor{\complexparamcolor}{\symImpedance_{ \indexGridLines ,\symPhaseG \symPhaseB }^{\seriesss}    }}
\newcommand{\zijksGC}[0]{\textcolor{\complexparamcolor}{\symImpedance_{ \indexGridLines , \symPhaseG\symPhaseC }^{\seriesss}    }}
\newcommand{\zijksGN}[0]{\textcolor{\complexparamcolor}{\symImpedance_{ \indexGridLines , \symPhaseG \symPhaseN}^{\seriesss}    }}

\newcommand{\zijksGG}[0]{\textcolor{\complexparamcolor}{\symImpedance_{ \indexGridLines , \symPhaseG \symPhaseG}^{\seriesss}    }}

\newcommand{\zijksNNhat}[0]{\textcolor{\complexparamcolor}{\hat{\symImpedance}_{ \indexGridLines , \symPhaseN \symPhaseN}^{\seriesss}    }}
\newcommand{\zijksANhat}[0]{\textcolor{\complexparamcolor}{\hat{\symImpedance}_{ \indexGridLines , \symPhaseA \symPhaseN}^{\seriesss}    }}
\newcommand{\zijksBNhat}[0]{\textcolor{\complexparamcolor}{\hat{\symImpedance}_{ \indexGridLines , \symPhaseB \symPhaseN}^{\seriesss}    }}
\newcommand{\zijksCNhat}[0]{\textcolor{\complexparamcolor}{\hat{\symImpedance}_{ \indexGridLines , \symPhaseC \symPhaseN}^{\seriesss}    }}

\newcommand{\zijksAAhat}[0]{\textcolor{\complexparamcolor}{\hat{\symImpedance}_{ \indexGridLines , \symPhaseA \symPhaseA}^{\seriesss}    }}
\newcommand{\zijksBAhat}[0]{\textcolor{\complexparamcolor}{\hat{\symImpedance}_{ \indexGridLines , \symPhaseB \symPhaseA}^{\seriesss}    }}
\newcommand{\zijksCAhat}[0]{\textcolor{\complexparamcolor}{\hat{\symImpedance}_{ \indexGridLines , \symPhaseC \symPhaseA}^{\seriesss}    }}

\newcommand{\zijksABhat}[0]{\textcolor{\complexparamcolor}{\hat{\symImpedance}_{ \indexGridLines , \symPhaseA \symPhaseB}^{\seriesss}    }}
\newcommand{\zijksBBhat}[0]{\textcolor{\complexparamcolor}{\hat{\symImpedance}_{ \indexGridLines , \symPhaseB \symPhaseB}^{\seriesss}    }}
\newcommand{\zijksCBhat}[0]{\textcolor{\complexparamcolor}{\hat{\symImpedance}_{ \indexGridLines , \symPhaseC \symPhaseB}^{\seriesss}    }}

\newcommand{\zijksAChat}[0]{\textcolor{\complexparamcolor}{\hat{\symImpedance}_{ \indexGridLines , \symPhaseA \symPhaseC}^{\seriesss}    }}
\newcommand{\zijksBChat}[0]{\textcolor{\complexparamcolor}{\hat{\symImpedance}_{ \indexGridLines , \symPhaseB \symPhaseC}^{\seriesss}    }}
\newcommand{\zijksCChat}[0]{\textcolor{\complexparamcolor}{\hat{\symImpedance}_{ \indexGridLines , \symPhaseC \symPhaseC}^{\seriesss}    }}

\newcommand{\zijksNAhat}[0]{\textcolor{\complexparamcolor}{\hat{\symImpedance}_{ \indexGridLines , \symPhaseN \symPhaseA }^{\seriesss}    }}
\newcommand{\zijksNBhat}[0]{\textcolor{\complexparamcolor}{\hat{\symImpedance}_{ \indexGridLines , \symPhaseN \symPhaseB }^{\seriesss}    }}
\newcommand{\zijksNChat}[0]{\textcolor{\complexparamcolor}{\hat{\symImpedance}_{ \indexGridLines , \symPhaseN \symPhaseC }^{\seriesss}    }}

\newcommand{\yijksAA}[0]{\textcolor{\complexparamcolor}{\symAdmittance_{ \indexGridLines , \symPhaseA \symPhaseA}^{\seriesss}    }}
\newcommand{\yijksAB}[0]{\textcolor{\complexparamcolor}{\symAdmittance_{ \indexGridLines , \symPhaseA \symPhaseB}^{\seriesss}    }}
\newcommand{\yijksAC}[0]{\textcolor{\complexparamcolor}{\symAdmittance_{ \indexGridLines , \symPhaseA \symPhaseC}^{\seriesss}    }}
\newcommand{\yijksBA}[0]{\textcolor{\complexparamcolor}{\symAdmittance_{ \indexGridLines , \symPhaseB \symPhaseA}^{\seriesss}    }}
\newcommand{\yijksBB}[0]{\textcolor{\complexparamcolor}{\symAdmittance_{ \indexGridLines , \symPhaseB \symPhaseB}^{\seriesss}    }}
\newcommand{\yijksBC}[0]{\textcolor{\complexparamcolor}{\symAdmittance_{ \indexGridLines , \symPhaseB \symPhaseC}^{\seriesss}    }}
\newcommand{\yijksCA}[0]{\textcolor{\complexparamcolor}{\symAdmittance_{ \indexGridLines , \symPhaseC \symPhaseA}^{\seriesss}    }}
\newcommand{\yijksCB}[0]{\textcolor{\complexparamcolor}{\symAdmittance_{ \indexGridLines , \symPhaseC \symPhaseB}^{\seriesss}    }}
\newcommand{\yijksCC}[0]{\textcolor{\complexparamcolor}{\symAdmittance_{ \indexGridLines , \symPhaseC \symPhaseC}^{\seriesss}    }}

\newcommand{\gijksAA}[0]{\textcolor{\paramcolor}{\symConductance_{ \indexGridLines , \symPhaseA \symPhaseA}^{\seriesss}    }}
\newcommand{\gijksAB}[0]{\textcolor{\paramcolor}{\symConductance_{ \indexGridLines , \symPhaseA \symPhaseB}^{\seriesss}    }}
\newcommand{\gijksAC}[0]{\textcolor{\paramcolor}{\symConductance_{ \indexGridLines , \symPhaseA \symPhaseC}^{\seriesss}    }}
\newcommand{\gijksBA}[0]{\textcolor{\paramcolor}{\symConductance_{ \indexGridLines , \symPhaseB \symPhaseA}^{\seriesss}    }}
\newcommand{\gijksBB}[0]{\textcolor{\paramcolor}{\symConductance_{ \indexGridLines , \symPhaseB \symPhaseB}^{\seriesss}    }}
\newcommand{\gijksBC}[0]{\textcolor{\paramcolor}{\symConductance_{ \indexGridLines , \symPhaseB \symPhaseC}^{\seriesss}    }}
\newcommand{\gijksCA}[0]{\textcolor{\paramcolor}{\symConductance_{ \indexGridLines , \symPhaseC \symPhaseA}^{\seriesss}    }}
\newcommand{\gijksCB}[0]{\textcolor{\paramcolor}{\symConductance_{ \indexGridLines , \symPhaseC \symPhaseB}^{\seriesss}    }}
\newcommand{\gijksCC}[0]{\textcolor{\paramcolor}{\symConductance_{ \indexGridLines , \symPhaseC \symPhaseC}^{\seriesss}    }}

\newcommand{\bijksAA}[0]{\textcolor{\paramcolor}{\symSusceptance_{ \indexGridLines , \symPhaseA \symPhaseA}^{\seriesss}    }}
\newcommand{\bijksAB}[0]{\textcolor{\paramcolor}{\symSusceptance_{ \indexGridLines , \symPhaseA \symPhaseB}^{\seriesss}    }}
\newcommand{\bijksAC}[0]{\textcolor{\paramcolor}{\symSusceptance_{ \indexGridLines , \symPhaseA \symPhaseC}^{\seriesss}    }}
\newcommand{\bijksBA}[0]{\textcolor{\paramcolor}{\symSusceptance_{ \indexGridLines , \symPhaseB \symPhaseA}^{\seriesss}    }}
\newcommand{\bijksBB}[0]{\textcolor{\paramcolor}{\symSusceptance_{ \indexGridLines , \symPhaseB \symPhaseB}^{\seriesss}    }}
\newcommand{\bijksBC}[0]{\textcolor{\paramcolor}{\symSusceptance_{ \indexGridLines , \symPhaseB \symPhaseC}^{\seriesss}    }}
\newcommand{\bijksCA}[0]{\textcolor{\paramcolor}{\symSusceptance_{ \indexGridLines , \symPhaseC \symPhaseA}^{\seriesss}    }}
\newcommand{\bijksCB}[0]{\textcolor{\paramcolor}{\symSusceptance_{ \indexGridLines , \symPhaseC \symPhaseB}^{\seriesss}    }}
\newcommand{\bijksCC}[0]{\textcolor{\paramcolor}{\symSusceptance_{ \indexGridLines , \symPhaseC \symPhaseC}^{\seriesss}    }}

\newcommand{\IijkA}[0]{\textcolor{\complexcolor}{\symCurrent_{\indexGridLines\indexGridNode \indexGridNodeTwo, \symPhaseA}^{}  }}
\newcommand{\IijkB}[0]{\textcolor{\complexcolor}{\symCurrent_{\indexGridLines\indexGridNode \indexGridNodeTwo, \symPhaseB}^{}  }}
\newcommand{\IijkC}[0]{\textcolor{\complexcolor}{\symCurrent_{\indexGridLines\indexGridNode \indexGridNodeTwo, \symPhaseC}^{}  }}
\newcommand{\IijkN}[0]{\textcolor{\complexcolor}{\symCurrent_{\indexGridLines\indexGridNode \indexGridNodeTwo, \symPhaseN}^{}  }}
\newcommand{\IijkG}[0]{\textcolor{\complexcolor}{\symCurrent_{\indexGridLines\indexGridNode \indexGridNodeTwo, \symPhaseG}^{}  }}

\newcommand{\kifkG}[0]{\textcolor{\complexcolor}{\symCurrent_{k\indexGridNode f, \symPhaseG}^{}  }}

\newcommand{\IijkNreal}[0]{\textcolor{black}{\symCurrent_{\indexGridLines\indexGridNode \indexGridNodeTwo, \symPhaseN}^{\text{re}} }}
\newcommand{\IijkNimag}[0]{\textcolor{black}{\symCurrent_{\indexGridLines\indexGridNode \indexGridNodeTwo, \symPhaseN}^{\text{im}} }}

\newcommand{\IikNG}[0]{\textcolor{\complexcolor}{\symCurrent_{\indexGridNode , \symPhaseN\symPhaseG}}}
\newcommand{\IikNGreal}[0]{\textcolor{black}{\symCurrent_{\indexGridNode , \symPhaseN\symPhaseG}^{\text{re}} }}
\newcommand{\IikNGimag}[0]{\textcolor{black}{\symCurrent_{\indexGridNode , \symPhaseN\symPhaseG}^{\text{im}}}}

\newcommand{\Iijkzero}[0]{\textcolor{\complexcolor}{\symCurrent_{\indexGridLines\indexGridNode \indexGridNodeTwo, 0}^{}  }}
\newcommand{\Iijkone}[0]{\textcolor{\complexcolor}{\symCurrent_{\indexGridLines\indexGridNode \indexGridNodeTwo, 1}^{}  }}
\newcommand{\Iijktwo}[0]{\textcolor{\complexcolor}{\symCurrent_{\indexGridLines\indexGridNode \indexGridNodeTwo, 2}^{}  }}

\newcommand{\IjikA}[0]{\textcolor{\complexcolor}{\symCurrent_{\indexGridLines\indexGridNodeTwo\indexGridNode , \symPhaseA}^{}  }}
\newcommand{\IjikB}[0]{\textcolor{\complexcolor}{\symCurrent_{\indexGridLines\indexGridNodeTwo\indexGridNode  , \symPhaseB}^{}  }}
\newcommand{\IjikC}[0]{\textcolor{\complexcolor}{\symCurrent_{\indexGridLines\indexGridNodeTwo\indexGridNode , \symPhaseC}^{}  }}
\newcommand{\IjikN}[0]{\textcolor{\complexcolor}{\symCurrent_{\indexGridLines\indexGridNodeTwo\indexGridNode , \symPhaseN}^{}  }}
\newcommand{\IjikG}[0]{\textcolor{\complexcolor}{\symCurrent_{\indexGridLines\indexGridNodeTwo\indexGridNode , \symPhaseG}^{}  }}

\newcommand{\IijkrefA}[0]{\textcolor{\complexparamcolor}{\symCurrent_{\indexGridLines\indexGridNode \indexGridNodeTwo,  \symPhaseA}^{\refss}  }}
\newcommand{\IijkrefB}[0]{\textcolor{\complexparamcolor}{\symCurrent_{\indexGridLines\indexGridNode \indexGridNodeTwo, \symPhaseB}^{\refss}  }}
\newcommand{\IijkrefC}[0]{\textcolor{\complexparamcolor}{\symCurrent_{\indexGridLines\indexGridNode \indexGridNodeTwo, \symPhaseC}^{\refss}  }}
\newcommand{\IijksrefA}[0]{\textcolor{\complexparamcolor}{\symCurrent_{\indexGridLines\indexGridNode \indexGridNodeTwo, \seriesss, \symPhaseA}^{\refss}  }}
\newcommand{\IijksrefB}[0]{\textcolor{\complexparamcolor}{\symCurrent_{\indexGridLines\indexGridNode \indexGridNodeTwo, \seriesss, \symPhaseB}^{\refss}  }}
\newcommand{\IijksrefC}[0]{\textcolor{\complexparamcolor}{\symCurrent_{\indexGridLines\indexGridNode \indexGridNodeTwo, \seriesss, \symPhaseC}^{\refss}  }}

\newcommand{\IijksA}[0]{\textcolor{\complexcolor}{\symCurrent_{\indexGridLines\indexGridNode \indexGridNodeTwo, \symPhaseA}^{ \seriesss}  }}
\newcommand{\IijksB}[0]{\textcolor{\complexcolor}{\symCurrent_{\indexGridLines\indexGridNode \indexGridNodeTwo, \symPhaseB}^{\seriesss}  }}
\newcommand{\IijksC}[0]{\textcolor{\complexcolor}{\symCurrent_{\indexGridLines\indexGridNode \indexGridNodeTwo, \symPhaseC}^{\seriesss}  }}
\newcommand{\IijksN}[0]{\textcolor{\complexcolor}{\symCurrent_{\indexGridLines\indexGridNode \indexGridNodeTwo, \symPhaseN}^{\seriesss}  }}
\newcommand{\IijksG}[0]{\textcolor{\complexcolor}{\symCurrent_{\indexGridLines\indexGridNode \indexGridNodeTwo, \symPhaseG}^{\seriesss}  }}

\newcommand{\IjiksA}[0]{\textcolor{\complexcolor}{\symCurrent_{\indexGridLines\indexGridNodeTwo\indexGridNode , \symPhaseA}^{ \seriesss}  }}
\newcommand{\IjiksB}[0]{\textcolor{\complexcolor}{\symCurrent_{\indexGridLines\indexGridNodeTwo\indexGridNode , \symPhaseB}^{\seriesss}  }}
\newcommand{\IjiksC}[0]{\textcolor{\complexcolor}{\symCurrent_{\indexGridLines\indexGridNodeTwo\indexGridNode , \symPhaseC}^{\seriesss}  }}
\newcommand{\IjiksN}[0]{\textcolor{\complexcolor}{\symCurrent_{\indexGridLines\indexGridNodeTwo\indexGridNode , \symPhaseN}^{\seriesss}  }}

\newcommand{\IijkshA}[0]{\textcolor{\complexcolor}{\symCurrent_{\indexGridLines\indexGridNode \indexGridNodeTwo, \symPhaseA}^{ \shuntss}  }}
\newcommand{\IijkshB}[0]{\textcolor{\complexcolor}{\symCurrent_{\indexGridLines\indexGridNode \indexGridNodeTwo, \symPhaseB}^{\shuntss}  }}
\newcommand{\IijkshC}[0]{\textcolor{\complexcolor}{\symCurrent_{\indexGridLines\indexGridNode \indexGridNodeTwo, \symPhaseC}^{\shuntss}  }}
\newcommand{\IijkshN}[0]{\textcolor{\complexcolor}{\symCurrent_{\indexGridLines\indexGridNode \indexGridNodeTwo, \symPhaseN}^{\shuntss}  }}
\newcommand{\IijkshG}[0]{\textcolor{\complexcolor}{\symCurrent_{\indexGridLines\indexGridNode \indexGridNodeTwo, \symPhaseG}^{\shuntss}  }}

\newcommand{\IjikshN}[0]{\textcolor{\complexcolor}{\symCurrent_{\indexGridLines \indexGridNodeTwo \indexGridNode , \symPhaseN}^{\shuntss}  }}

\newcommand{\IikA}[0]{\textcolor{\complexcolor}{\symCurrent_{\indexGridNode , \symPhaseA}^{}  }}
\newcommand{\IikB}[0]{\textcolor{\complexcolor}{\symCurrent_{\indexGridNode , \symPhaseB}^{}  }}
\newcommand{\IikC}[0]{\textcolor{\complexcolor}{\symCurrent_{\indexGridNode , \symPhaseC}^{}  }}

\newcommand{\yijkshAA}[0]{\textcolor{\complexparamcolor}{\symAdmittance_{\indexGridLines \indexGridNode \indexGridNodeTwo, \symPhaseA \symPhaseA}^{\shuntss}    }}
\newcommand{\yijkshAB}[0]{\textcolor{\complexparamcolor}{\symAdmittance_{\indexGridLines \indexGridNode \indexGridNodeTwo, \symPhaseA \symPhaseB}^{\shuntss}    }}
\newcommand{\yijkshAC}[0]{\textcolor{\complexparamcolor}{\symAdmittance_{\indexGridLines \indexGridNode \indexGridNodeTwo, \symPhaseA \symPhaseC}^{\shuntss}    }}
\newcommand{\yijkshAN}[0]{\textcolor{\complexparamcolor}{\symAdmittance_{\indexGridLines \indexGridNode \indexGridNodeTwo, \symPhaseA \symPhaseN}^{\shuntss}    }}
\newcommand{\yijkshAG}[0]{\textcolor{\complexparamcolor}{\symAdmittance_{\indexGridLines \indexGridNode \indexGridNodeTwo, \symPhaseA \symPhaseG}^{\shuntss}    }}

\newcommand{\yijkshBA}[0]{\textcolor{\complexparamcolor}{\symAdmittance_{\indexGridLines \indexGridNode \indexGridNodeTwo, \symPhaseB \symPhaseA}^{\shuntss}    }}
\newcommand{\yijkshBB}[0]{\textcolor{\complexparamcolor}{\symAdmittance_{\indexGridLines \indexGridNode \indexGridNodeTwo, \symPhaseB \symPhaseB}^{\shuntss}    }}
\newcommand{\yijkshBC}[0]{\textcolor{\complexparamcolor}{\symAdmittance_{\indexGridLines \indexGridNode \indexGridNodeTwo, \symPhaseB \symPhaseC}^{\shuntss}    }}
\newcommand{\yijkshBN}[0]{\textcolor{\complexparamcolor}{\symAdmittance_{\indexGridLines \indexGridNode \indexGridNodeTwo, \symPhaseB \symPhaseN}^{\shuntss}    }}
\newcommand{\yijkshBG}[0]{\textcolor{\complexparamcolor}{\symAdmittance_{\indexGridLines \indexGridNode \indexGridNodeTwo, \symPhaseB \symPhaseG}^{\shuntss}    }}

\newcommand{\yijkshCA}[0]{\textcolor{\complexparamcolor}{\symAdmittance_{\indexGridLines \indexGridNode \indexGridNodeTwo, \symPhaseC \symPhaseA}^{\shuntss}    }}
\newcommand{\yijkshCB}[0]{\textcolor{\complexparamcolor}{\symAdmittance_{\indexGridLines \indexGridNode \indexGridNodeTwo, \symPhaseC \symPhaseB}^{\shuntss}    }}
\newcommand{\yijkshCC}[0]{\textcolor{\complexparamcolor}{\symAdmittance_{\indexGridLines \indexGridNode \indexGridNodeTwo, \symPhaseC \symPhaseC}^{\shuntss}    }}\newcommand{\yijkshCN}[0]{\textcolor{\complexparamcolor}{\symAdmittance_{\indexGridLines \indexGridNode \indexGridNodeTwo, \symPhaseC \symPhaseN}^{\shuntss}    }}\newcommand{\yijkshCG}[0]{\textcolor{\complexparamcolor}{\symAdmittance_{\indexGridLines \indexGridNode \indexGridNodeTwo, \symPhaseC \symPhaseG}^{\shuntss}    }}

\newcommand{\yijkshNA}[0]{\textcolor{\complexparamcolor}{\symAdmittance_{\indexGridLines \indexGridNode \indexGridNodeTwo, \symPhaseN \symPhaseA}^{\shuntss}    }}
\newcommand{\yijkshNB}[0]{\textcolor{\complexparamcolor}{\symAdmittance_{\indexGridLines \indexGridNode \indexGridNodeTwo, \symPhaseN \symPhaseB}^{\shuntss}    }}
\newcommand{\yijkshNC}[0]{\textcolor{\complexparamcolor}{\symAdmittance_{\indexGridLines \indexGridNode \indexGridNodeTwo, \symPhaseN \symPhaseC}^{\shuntss}    }}\newcommand{\yijkshNN}[0]{\textcolor{\complexparamcolor}{\symAdmittance_{\indexGridLines \indexGridNode \indexGridNodeTwo, \symPhaseN \symPhaseN}^{\shuntss}    }}\newcommand{\yijkshNG}[0]{\textcolor{\complexparamcolor}{\symAdmittance_{\indexGridLines \indexGridNode \indexGridNodeTwo, \symPhaseN \symPhaseG}^{\shuntss}    }}

\newcommand{\yijkshGA}[0]{\textcolor{\complexparamcolor}{\symAdmittance_{\indexGridLines \indexGridNode \indexGridNodeTwo, \symPhaseG \symPhaseA}^{\shuntss}    }}
\newcommand{\yijkshGB}[0]{\textcolor{\complexparamcolor}{\symAdmittance_{\indexGridLines \indexGridNode \indexGridNodeTwo, \symPhaseG \symPhaseB}^{\shuntss}    }}
\newcommand{\yijkshGC}[0]{\textcolor{\complexparamcolor}{\symAdmittance_{\indexGridLines \indexGridNode \indexGridNodeTwo, \symPhaseG \symPhaseC}^{\shuntss}    }}\newcommand{\yijkshGN}[0]{\textcolor{\complexparamcolor}{\symAdmittance_{\indexGridLines \indexGridNode \indexGridNodeTwo, \symPhaseG \symPhaseN}^{\shuntss}    }}\newcommand{\yijkshGG}[0]{\textcolor{\complexparamcolor}{\symAdmittance_{\indexGridLines \indexGridNode \indexGridNodeTwo, \symPhaseG \symPhaseG}^{\shuntss}    }}

\newcommand{\gijkshAA}[0]{\textcolor{\paramcolor}{\symConductance_{ \indexGridLines\indexGridNode \indexGridNodeTwo, \symPhaseA \symPhaseA}^{\shuntss}    }}
\newcommand{\gijkshAB}[0]{\textcolor{\paramcolor}{\symConductance_{\indexGridLines \indexGridNode \indexGridNodeTwo, \symPhaseA \symPhaseB}^{\shuntss}    }}
\newcommand{\gijkshAC}[0]{\textcolor{\paramcolor}{\symConductance_{\indexGridLines \indexGridNode \indexGridNodeTwo, \symPhaseA \symPhaseC}^{\shuntss}    }}

\newcommand{\gijkshBA}[0]{\textcolor{\paramcolor}{\symConductance_{ \indexGridLines\indexGridNode \indexGridNodeTwo, \symPhaseB \symPhaseA}^{\shuntss}    }}
\newcommand{\gijkshBB}[0]{\textcolor{\paramcolor}{\symConductance_{\indexGridLines \indexGridNode \indexGridNodeTwo, \symPhaseB \symPhaseB}^{\shuntss}    }}
\newcommand{\gijkshBC}[0]{\textcolor{\paramcolor}{\symConductance_{\indexGridLines \indexGridNode \indexGridNodeTwo, \symPhaseB \symPhaseC}^{\shuntss}    }}

\newcommand{\gijkshCA}[0]{\textcolor{\paramcolor}{\symConductance_{ \indexGridLines\indexGridNode \indexGridNodeTwo, \symPhaseC \symPhaseA}^{\shuntss}    }}
\newcommand{\gijkshCB}[0]{\textcolor{\paramcolor}{\symConductance_{ \indexGridLines\indexGridNode \indexGridNodeTwo, \symPhaseC \symPhaseB}^{\shuntss}    }}
\newcommand{\gijkshCC}[0]{\textcolor{\paramcolor}{\symConductance_{\indexGridLines \indexGridNode \indexGridNodeTwo, \symPhaseC \symPhaseC}^{\shuntss}    }}

\newcommand{\bijkshAA}[0]{\textcolor{\paramcolor}{\symSusceptance_{ \indexGridLines\indexGridNode \indexGridNodeTwo, \symPhaseA \symPhaseA}^{\shuntss}    }}
\newcommand{\bijkshAB}[0]{\textcolor{\paramcolor}{\symSusceptance_{ \indexGridLines\indexGridNode \indexGridNodeTwo, \symPhaseA \symPhaseB}^{\shuntss}    }}
\newcommand{\bijkshAC}[0]{\textcolor{\paramcolor}{\symSusceptance_{ \indexGridLines\indexGridNode \indexGridNodeTwo, \symPhaseA \symPhaseC}^{\shuntss}    }}

\newcommand{\bijkshBA}[0]{\textcolor{\paramcolor}{\symSusceptance_{ \indexGridLines\indexGridNode \indexGridNodeTwo, \symPhaseB \symPhaseA}^{\shuntss}    }}
\newcommand{\bijkshBB}[0]{\textcolor{\paramcolor}{\symSusceptance_{ \indexGridLines\indexGridNode \indexGridNodeTwo, \symPhaseB \symPhaseB}^{\shuntss}    }}
\newcommand{\bijkshBC}[0]{\textcolor{\paramcolor}{\symSusceptance_{ \indexGridLines\indexGridNode \indexGridNodeTwo, \symPhaseB \symPhaseC}^{\shuntss}    }}

\newcommand{\bijkshCA}[0]{\textcolor{\paramcolor}{\symSusceptance_{ \indexGridLines\indexGridNode \indexGridNodeTwo, \symPhaseC \symPhaseA}^{\shuntss}    }}
\newcommand{\bijkshCB}[0]{\textcolor{\paramcolor}{\symSusceptance_{\indexGridLines \indexGridNode \indexGridNodeTwo, \symPhaseC \symPhaseB}^{\shuntss}    }}
\newcommand{\bijkshCC}[0]{\textcolor{\paramcolor}{\symSusceptance_{ \indexGridLines\indexGridNode \indexGridNodeTwo, \symPhaseC \symPhaseC}^{\shuntss}    }}

\newcommand{\yjikshAA}[0]{\textcolor{\complexparamcolor}{\symAdmittance_{\indexGridLines\indexGridNodeTwo \indexGridNode , \symPhaseA \symPhaseA}^{\shuntss}    }}
\newcommand{\yjikshAB}[0]{\textcolor{\complexparamcolor}{\symAdmittance_{\indexGridLines\indexGridNodeTwo \indexGridNode , \symPhaseA \symPhaseB}^{\shuntss}    }}
\newcommand{\yjikshAC}[0]{\textcolor{\complexparamcolor}{\symAdmittance_{\indexGridLines \indexGridNodeTwo\indexGridNode , \symPhaseA \symPhaseC}^{\shuntss}    }}
\newcommand{\yjikshAN}[0]{\textcolor{\complexparamcolor}{\symAdmittance_{\indexGridLines \indexGridNodeTwo\indexGridNode , \symPhaseA \symPhaseN}^{\shuntss}    }}
\newcommand{\yjikshAG}[0]{\textcolor{\complexparamcolor}{\symAdmittance_{\indexGridLines \indexGridNodeTwo\indexGridNode , \symPhaseA \symPhaseG}^{\shuntss}    }}

\newcommand{\yjikshBA}[0]{\textcolor{\complexparamcolor}{\symAdmittance_{\indexGridLines\indexGridNodeTwo \indexGridNode , \symPhaseB \symPhaseA}^{\shuntss}    }}
\newcommand{\yjikshBB}[0]{\textcolor{\complexparamcolor}{\symAdmittance_{\indexGridLines\indexGridNodeTwo \indexGridNode , \symPhaseB \symPhaseB}^{\shuntss}    }}
\newcommand{\yjikshBC}[0]{\textcolor{\complexparamcolor}{\symAdmittance_{\indexGridLines \indexGridNodeTwo\indexGridNode , \symPhaseB \symPhaseC}^{\shuntss}    }}
\newcommand{\yjikshBN}[0]{\textcolor{\complexparamcolor}{\symAdmittance_{\indexGridLines \indexGridNodeTwo\indexGridNode , \symPhaseB \symPhaseN}^{\shuntss}    }}
\newcommand{\yjikshBG}[0]{\textcolor{\complexparamcolor}{\symAdmittance_{\indexGridLines \indexGridNodeTwo\indexGridNode , \symPhaseB \symPhaseG}^{\shuntss}    }}

\newcommand{\yjikshCA}[0]{\textcolor{\complexparamcolor}{\symAdmittance_{\indexGridLines\indexGridNodeTwo \indexGridNode , \symPhaseC \symPhaseA}^{\shuntss}    }}
\newcommand{\yjikshCB}[0]{\textcolor{\complexparamcolor}{\symAdmittance_{\indexGridLines\indexGridNodeTwo \indexGridNode , \symPhaseC \symPhaseB}^{\shuntss}    }}
\newcommand{\yjikshCC}[0]{\textcolor{\complexparamcolor}{\symAdmittance_{\indexGridLines \indexGridNodeTwo\indexGridNode , \symPhaseC \symPhaseC}^{\shuntss}    }}
\newcommand{\yjikshCN}[0]{\textcolor{\complexparamcolor}{\symAdmittance_{\indexGridLines \indexGridNodeTwo\indexGridNode , \symPhaseC \symPhaseN}^{\shuntss}    }}
\newcommand{\yjikshCG}[0]{\textcolor{\complexparamcolor}{\symAdmittance_{\indexGridLines \indexGridNodeTwo\indexGridNode , \symPhaseC \symPhaseG}^{\shuntss}    }}

\newcommand{\yjikshNA}[0]{\textcolor{\complexparamcolor}{\symAdmittance_{\indexGridLines\indexGridNodeTwo \indexGridNode , \symPhaseN \symPhaseA}^{\shuntss}    }}
\newcommand{\yjikshNB}[0]{\textcolor{\complexparamcolor}{\symAdmittance_{\indexGridLines\indexGridNodeTwo \indexGridNode , \symPhaseN \symPhaseB}^{\shuntss}    }}
\newcommand{\yjikshNC}[0]{\textcolor{\complexparamcolor}{\symAdmittance_{\indexGridLines \indexGridNodeTwo\indexGridNode , \symPhaseN \symPhaseC}^{\shuntss}    }}
\newcommand{\yjikshNN}[0]{\textcolor{\complexparamcolor}{\symAdmittance_{\indexGridLines \indexGridNodeTwo\indexGridNode , \symPhaseN \symPhaseN}^{\shuntss}    }}
\newcommand{\yjikshNG}[0]{\textcolor{\complexparamcolor}{\symAdmittance_{\indexGridLines \indexGridNodeTwo\indexGridNode , \symPhaseN \symPhaseG}^{\shuntss}    }}

\newcommand{\yjikshGA}[0]{\textcolor{\complexparamcolor}{\symAdmittance_{\indexGridLines\indexGridNodeTwo \indexGridNode , \symPhaseG \symPhaseA}^{\shuntss}    }}
\newcommand{\yjikshGB}[0]{\textcolor{\complexparamcolor}{\symAdmittance_{\indexGridLines\indexGridNodeTwo \indexGridNode , \symPhaseG \symPhaseB}^{\shuntss}    }}
\newcommand{\yjikshGC}[0]{\textcolor{\complexparamcolor}{\symAdmittance_{\indexGridLines \indexGridNodeTwo\indexGridNode , \symPhaseG \symPhaseC}^{\shuntss}    }}
\newcommand{\yjikshGN}[0]{\textcolor{\complexparamcolor}{\symAdmittance_{\indexGridLines \indexGridNodeTwo\indexGridNode , \symPhaseG \symPhaseN}^{\shuntss}    }}
\newcommand{\yjikshGG}[0]{\textcolor{\complexparamcolor}{\symAdmittance_{\indexGridLines \indexGridNodeTwo\indexGridNode , \symPhaseG \symPhaseG}^{\shuntss}    }}

\newcommand{\gijkshpp}[0]{\textcolor{\paramcolor}{\symConductance_{\indexGridLines \indexGridNode \indexGridNodeTwo, \indexPhases\indexPhases}^{\shuntss}    }}
\newcommand{\gijkshph}[0]{\textcolor{\paramcolor}{\symConductance_{ \indexGridLines\indexGridNode \indexGridNodeTwo, \indexPhases\indexPhasesTwo}^{\shuntss}    }}

\newcommand{\gijkspp}[0]{\textcolor{\paramcolor}{\symConductance_{ \indexGridLines  , \indexPhases\indexPhases}^{\seriesss}    }}
\newcommand{\gijksph}[0]{\textcolor{\paramcolor}{\symConductance_{ \indexGridLines  , \indexPhases\indexPhasesTwo}^{\seriesss}    }}

\newcommand{\bijkshpp}[0]{\textcolor{\paramcolor}{\symSusceptance_{\indexGridLines \indexGridNode \indexGridNodeTwo, \indexPhases\indexPhases}^{\shuntss}    }}
\newcommand{\bijkspp}[0]{\textcolor{\paramcolor}{\symSusceptance_{ \indexGridLines, \indexPhases\indexPhases}^{\seriesss}    }}
\newcommand{\bijksph}[0]{\textcolor{\paramcolor}{\symSusceptance_{ \indexGridLines, \indexPhases\indexPhasesTwo}^{\seriesss}    }}

\newcommand{\bijkshph}[0]{\textcolor{\paramcolor}{\symSusceptance_{ \indexGridNode \indexGridNodeTwo, \indexPhases\indexPhasesTwo}^{\shuntss}    }}

\newcommand{\zbental}[0]{{z}}
\newcommand{\ybental}[0]{{y}}
\newcommand{\xbental}[0]{{x}}
\newcommand{\xonebental}[0]{\textcolor{black}{\xbental_{1} }}
\newcommand{\xtwobental}[0]{\textcolor{black}{\xbental_{2} }}
\newcommand{\xthreebental}[0]{\textcolor{black}{\xbental_{3} }}
\newcommand{\xfourbental}[0]{\textcolor{black}{\xbental_{4} }}
\newcommand{\xfivebental}[0]{\textcolor{black}{\xbental_{5} }}
\newcommand{\xsixbental}[0]{\textcolor{black}{\xbental_{6} }}
\newcommand{\xsevenbental}[0]{\textcolor{black}{\xbental_{7} }}
\newcommand{\xibental}[0]{\textcolor{black}{\xbental_{\counterbental} }}
\newcommand{\upperbental}[0]{{U}}
\newcommand{\lowerbental}[0]{{L}}
\newcommand{\xilowerbental}[0]{\textcolor{\paramcolor}{\xibental^{ \lowerbental } }}
\newcommand{\xiupperbental}[0]{\textcolor{\paramcolor}{\xibental^{ \upperbental } }}
\newcommand{\xlowerbental}[1]{\textcolor{\boundscolor}{\xbental_{#1}^{ \lowerbental } }}
\newcommand{\xupperbental}[1]{\textcolor{\boundscolor}{\xbental_{#1}^{ \upperbental } }}

\newcommand{\sigmabental}[1]{\textcolor{black}{\sigma_{(#1)}}}
\newcommand{\etabental}[1]{\textcolor{black}{\eta_{(#1)}}}
\newcommand{\counterbental}[0]{\textcolor{\paramcolor}{i}}
\newcommand{\countbental}[0]{\textcolor{\paramcolor}{\nu}}
\newcommand{\errorbental}[0]{\textcolor{\paramcolor}{\epsilon}}
\newcommand{\countpoly}[0]{\textcolor{\paramcolor}{\mu}}

\newcommand{\circumscribedopoly}[0]{\textcolor{\paramcolor}{\symSetting_{\circumss}}}
\newcommand{\rotatedopoly}[0]{\textcolor{\paramcolor}{\symSetting_{\rotatedss}}}
\newcommand{\genericbinary}[1]{\textcolor{black}{b_{#1}}}


\newcommand{\pipediameterijk}[0]{ \textcolor{\paramcolor}{  \symFluidPipeDiameter_{\indexGridLines}     }}
\newcommand{\pipeareaijk}[0]{ \textcolor{\paramcolor}{  \symFluidPipeArea_{\indexGridLines}     }}

\newcommand{\fluidpressure}[1]{\textcolor{black}{\symPressure_{#1}}}
\newcommand{\fluidpressureik}[0]{\textcolor{black}{\symPressure_{\indexGridNode }    }}
\newcommand{\fluidpressurejk}[0]{\textcolor{black}{\symPressure_{\indexGridNodeTwo }    }}

\newcommand{\fluidpressureikmin}[0]{\textcolor{\paramcolor}{\symPressure_{\indexGridNode }^{\minss}    }}
\newcommand{\fluidpressureikmax}[0]{\textcolor{\paramcolor}{\symPressure_{\indexGridNode }^{\maxss}    }}
\newcommand{\fluidpressurejkmin}[0]{\textcolor{\paramcolor}{\symPressure_{\indexGridNodeTwo }^{\minss}    }}
\newcommand{\fluidpressurejkmax}[0]{\textcolor{\paramcolor}{\symPressure_{\indexGridNodeTwo }^{\maxss}    }}

\newcommand{\fluidpressuresqik}[0]{\textcolor{black}{\symPressureSquared_{\indexGridNode }    }}
\newcommand{\fluidpressuresqjk}[0]{\textcolor{black}{\symPressureSquared_{\indexGridNodeTwo }    }}

\newcommand{\fluidpressuresqikmin}[0]{\textcolor{\paramcolor}{\symPressureSquared_{\indexGridNode }^{\minss}    }}
\newcommand{\fluidpressuresqjkmin}[0]{\textcolor{\paramcolor}{\symPressureSquared_{\indexGridNodeTwo }^{\minss}    }}
\newcommand{\fluidpressuresqikmax}[0]{\textcolor{\paramcolor}{\symPressureSquared_{\indexGridNode }^{\maxss}    }}
\newcommand{\fluidpressuresqjkmax}[0]{\textcolor{\paramcolor}{\symPressureSquared_{\indexGridNodeTwo }^{\maxss}    }}

\newcommand{\fluidflow}[1]{\textcolor{black}{\symVolumeFlow_{#1}}}
\newcommand{\fluidflowijk}[0]{\textcolor{black}{\fluidflow{\indexGridNode \indexGridNodeTwo}   }}

\newcommand{\fluidmassflow}[1]{\textcolor{black}{\symMassFlow_{#1}}}
\newcommand{\fluidmassflowijk}[0]{\textcolor{black}{\fluidmassflow{\indexGridNode \indexGridNodeTwo} }}

\newcommand{\fluidmassflowdirectionone}[0]{\textcolor{black}{y_{ij}^{+}}}
\newcommand{\fluidmassflowdirectiontwo}[0]{\textcolor{black}{y_{ij}^{-}}}

\newcommand{\fluiddensity}[1]{{\symDensity_{#1}  }}
\newcommand{\fluiddensityijk}[0]{\textcolor{\paramcolor}{\fluiddensity{\indexGridLines}   }}

\newcommand{\fluidviscosity}[1]{{\symViscosity_{#1}  }}
\newcommand{\fluidviscosityijk}[0]{\textcolor{\paramcolor}{\fluidviscosity{\indexGridLines}   }}

\newcommand{\heatflow}[1]{\textcolor{black}{\symHeatFlow_{#1}}}
\newcommand{\heatflowijk}[0]{\textcolor{black}{\symHeatFlow_{\indexGridNode \indexGridNodeTwo} }}

\newcommand{\heatresistance}[1]{\textcolor{\paramcolor}{\symHeatResistance_{#1}}}
\newcommand{\heatresistanceijk}[0]{ \textcolor{\paramcolor}{ \heatresistance{\indexGridLines}  }}   

\newcommand{\heattempijk}[0]{\textcolor{black}{\symTemperature_{\indexGridNode \indexGridNodeTwo}   }}
\newcommand{\heattempik}[0]{\textcolor{black}{\symTemperature_{\indexGridNode }   }}
\newcommand{\heattempjk}[0]{\textcolor{black}{\symTemperature_{ \indexGridNodeTwo}   }}

\newcommand{\fluidhagenresistance}[0]{ \textcolor{\paramcolor}{ w^{\hagenposeuilless}_{\indexGridLines}  }    }
\newcommand{\fluiddarcyresistance}[0]{ \textcolor{\paramcolor}{ w^{\darcyweisbachss}_{\indexGridLines}   }    }

\newcommand{\heatcapacity}[0]{ \textcolor{\paramcolor}{ \symHeatCapacity^{\text{\symPressure}} }}   


\newcommand{\symFlowVariable}[0]{   F }
\newcommand{\symFlowVariableSquare}[0]{   f }
\newcommand{\symPotentialVariable}[0]{ U }
\newcommand{\symPotentialAngleVariable}[0]{ \theta }
\newcommand{\symPotentialVariableSquare}[0]{ u }
\newcommand{\setFlowVariable}[0]{   \mathcal{F} }
\newcommand{\setPotentialVariable}[0]{ \mathcal{U} }
\newcommand{\indexFlowVariable}[0]{  f }
\newcommand{\indexPotentialVariable}[0]{ u }
\newcommand{\indexDomains}[0]{ d }
\newcommand{\setDomains}[0]{\mathcal{D}  }
\newcommand{\setDomainsExt}[0]{\mathcal{D}_0  }
\newcommand{\potentiali}[0]{ \symPotentialVariable_{\indexGridNode} }
\newcommand{\potentialz}[0]{ \symPotentialVariable_{\indexZIP} }
\newcommand{\potentialzlin}[0]{ \symPotentialVariable_{\indexZIP}^{\linss} }
\newcommand{\potentialilin}[0]{ \symPotentialVariable_{\indexGridNode}^{\linss} }
\newcommand{\potentialj}[0]{ \symPotentialVariable_{\indexGridNodeTwo} }
\newcommand{\potentialanglei}[0]{ \symPotentialAngleVariable_{\indexGridNode} }
\newcommand{\potentialanglej}[0]{ \symPotentialAngleVariable_{\indexGridNodeTwo} }
\newcommand{\potentialref}[0]{ \textcolor{\paramcolor}{ \symPotentialVariable^{\refss} }}
\newcommand{\potentialiref}[0]{ \textcolor{\paramcolor}{ \symPotentialVariable_{\indexGridNode}^{\refss} }}
\newcommand{\potentialzref}[0]{ \textcolor{\paramcolor}{ \symPotentialVariable_{\indexZIP}^{\refss} }}
\newcommand{\flowij}[0]{ \symFlowVariable_{\indexGridNode \indexGridNodeTwo} }
\newcommand{\flowi}[0]{ \symFlowVariable_{\indexGridNode } }
\newcommand{\flowz}[0]{ \symFlowVariable_{\indexZIP } }
\newcommand{\flowcomplexij}[0]{ \textcolor{\complexcolor}{\symFlowVariable_{\indexGridNode \indexGridNodeTwo} }}
\newcommand{\flowcomplexji}[0]{ \textcolor{\complexcolor}{\symFlowVariable_{\indexGridNodeTwo\indexGridNode } }}
\newcommand{\flowangleij}[0]{ {}{\symFlowVariable_{\indexGridNode \indexGridNodeTwo}^{\imagangle} }}
\newcommand{\flowji}[0]{ \symFlowVariable_{\indexGridNodeTwo \indexGridNode } }
\newcommand{\flowresistancetij}[0]{\textcolor{\paramcolor}{ r_{\indexGridLines } }}
\newcommand{\flowimpedancetij}[0]{\textcolor{\complexcolor}{ z_{\indexGridLines } }}
\newcommand{\flowreactancetij}[0]{\textcolor{\paramcolor}{ x_{\indexGridLines } }}
\newcommand{\potentialsqi}[0]{ \symPotentialVariableSquare_{\indexGridNode} }
\newcommand{\potentialsqz}[0]{ \symPotentialVariableSquare_{\indexZIP} }
\newcommand{\potentialsqj}[0]{ \symPotentialVariableSquare_{\indexGridNodeTwo} }
\newcommand{\flowsqij}[0]{ \symFlowVariableSquare_{\indexGridNode \indexGridNodeTwo} }
\newcommand{\potentialimin}[0]{ \textcolor{\paramcolor}{\symPotentialVariable_{\indexGridNode}^{\minss} }}
\newcommand{\potentialimax}[0]{ \textcolor{\paramcolor}{\symPotentialVariable_{\indexGridNode}^{\maxss} }}
\newcommand{\potentialirated}[0]{ \textcolor{\sizingcolor}{\symPotentialVariable_{\indexGridNode}^{\ratedss} }}
\newcommand{\potentialjmin}[0]{ \textcolor{\paramcolor}{\symPotentialVariable_{\indexGridNodeTwo}^{\minss} }}
\newcommand{\potentialjmax}[0]{ \textcolor{\paramcolor}{\symPotentialVariable_{\indexGridNodeTwo}^{\maxss} }}
\newcommand{\potentialjrated}[0]{ \textcolor{\sizingcolor}{\symPotentialVariable_{\indexGridNodeTwo}^{\ratedss} }}
\newcommand{\potentialijrated}[0]{ \textcolor{\sizingcolor}{\symPotentialVariable_{\indexGridNode\indexGridNodeTwo}^{\ratedss} }}
\newcommand{\potentialjirated}[0]{ \textcolor{\sizingcolor}{\symPotentialVariable_{\indexGridNodeTwo \indexGridNode}^{\ratedss} }}
\newcommand{\flowijmin}[0]{ \textcolor{\paramcolor}{\symFlowVariable_{\indexGridNode \indexGridNodeTwo}^{\minss} }}
\newcommand{\flowijmax}[0]{\textcolor{\paramcolor}{ \symFlowVariable_{\indexGridNode \indexGridNodeTwo}^{\maxss} }}
\newcommand{\flowijrated}[0]{\textcolor{\sizingcolor}{ \symFlowVariable_{\indexGridNode \indexGridNodeTwo}^{\ratedss} }}
\newcommand{\flowjirated}[0]{\textcolor{\sizingcolor}{ \symFlowVariable_{ \indexGridNodeTwo \indexGridNode }^{\ratedss} }}

\newcommand{\powerbalancefactorij}[0]{\textcolor{\paramcolor}{c_{\indexGridLines } }}

\newcommand{\enthalpycombustion}{\textcolor{\paramcolor}{\Delta \symEnthalpy^{\circ}_{\combustionss}}}
\newcommand{\enthalpyformation}{\textcolor{\paramcolor}{\Delta \symEnthalpy^{\circ}_{\formationss}}}

\newcommand{\slackconvexificationij}[0]{{\symSlack_{\indexGridNode \indexGridNodeTwo }^{\PFconvexss} }}
\newcommand{\slackconvexification}[0]{{\symSlack^{\PFconvexss} }}
\newcommand{\slackmaxconvexification}[0]{{\symSlack_{  }^{\maxss} }}
\newcommand{\slackminconvexification}[0]{{\symSlack_{  }^{\minss} }}
\newcommand{\slackabsconvexification}[0]{{\symSlack_{  }^{\absss} }}
\newcommand{\slackcomplexmagnitudeconvexification}[0]{{\symSlack_{  }^{\absss} }}
\newcommand{\slackbental}[0]{{\symSlack_{  }^{\bentallinearization} }}


\newcommand{\mapping}[0]{   \to }
\newcommand{\injection}[0]{   \rightarrowtail }
\newcommand{\surjection}[0]{   \twoheadrightarrow }
\newcommand{\bijection}[0]{   \leftrightarrow }

\newcommand{\symbollistlabelwith}{$\PPNeleclossk$}
\newcommand{\imagnumber}[0]{\textcolor{\paramcolor}{j}}
\newcommand{\imagangle}[0]{\textcolor{\paramcolor}{\angle}}

\begin{abstract}
This paper develops a novel second order cone relaxation of the semidefinite programming formulation of optimal power flow, that does not imply the `angle relaxation'. 
We build on a technique developed by Kim et al., extend it for complex matrices, and apply it to $3 \times 3$ positive semidefinite matrices to generate novel second-order cone constraints that augment upon the well-known $2\times2$ principal-minor based second-order cone constraints. Finally, we apply it to optimal power flow in meshed networks and provide numerical illustrations. 
\end{abstract}

\begin{IEEEkeywords}
Optimal Power Flow, Mathematical Optimization
\end{IEEEkeywords}


 \section{Introduction and Background} \label{intro}
Convex relaxation of nonlinear optimization models has been successfully applied to  optimal power flow  in the context of certifying global optimality or infeasibility \cite{Gopinath2020}, as well as in enabling mixed-integer nonlinear optimization. 
The tightest known relaxations are based on semidefinite programming (SDP) forms, either by construction, or generated as part of the Lasserre hierarchy \cite{Josz2018}. 
Since the scalability of SDP is limited in practice, so there remains a need for tight but  scalable relaxations. 
Second-order cone (SOC) programming is considered more scalable, and mixed-integer (MI) SOC solvers are commercially available. 

It is noted that the physical interpretation of the difference between the SDP and the SOC relaxation is the `angle relaxation' \cite{Low2014}, that is, Kirchhoff's voltage law in loops is \emph{not} represented beyond the implications of Ohm's law on voltage magnitude. 
Therefore, Kocuk et al. developed valid inequalities to include the physics of 3 and 4-cycles in a SOC formulation \cite{Kocuk2015}, by iteratively adding linear inequalities based on 3 and 4-cycle SDP subproblems. 
Alternatively, Coffrin et al. intersect the conventional SDP/SOC BIM power/lifted-voltage formulation with a  convex hull power/polar-voltage formulation \cite{Duben2015}, which therefore also improves upon the angle-related physics representation. 
Hijazi et al. \cite{Coffrin2015a} develop polynomial expressions for $3\times 3$ minors and propose an interpretation w.r.t. Kirchhoff's voltage law.

Alternatively, in this paper, we propose a method to generate tighter SOC relaxations, based on a technique of Kim et al. \cite{Kim2003}, and apply it to SDP optimization models for meshed  power networks. It will be shown that the novel SOC formulation improves upon the conventional, $2\times2$ principal minor based SOC relaxation (abbreviated PM SOC here) of the bus injection model (BIM) optimal power flow (OPF) formulation \cite{Low2014}. Finally, we provide an intuitive interpretation of the novel SOC constraints for 3-cycles by comparing the  expressions with those published by Hijazi et al.  \cite{Coffrin2015a} and Kocuk et al. \cite{Kocuk2017}.

 \section{Kim et al. SDP to SOC Relaxation} \label{soc_kim}
 Kim et al. \cite{Kim2003} develop a SDP to SOC relaxation procedure in the real numbers, whereas this section develops the generalization to  complex matrices; staying in the complex domain enables tighter relaxations. 
 Specifically, we adapt the derivation \S3 `An extension of the SOC relaxation of type 1', Lemma 3.2 \cite{Kim2003}  to  Hermitian matrices (a.k.a. complex symmetric). 
 
 We partition an Hermitian, positive semidefinite (PSD) variable matrix $\textcolor{\complexcolor}{\mathbf{M}}$ as,
  \begin{IEEEeqnarray}{C?s}
\textcolor{\complexcolor}{\mathbf{M}} =  \begin{bmatrix} 
\textcolor{black}{ \alpha } & \textcolor{\complexcolor}{\mathbf{a}}^{\hermitiantranspose} \\
 \textcolor{\complexcolor}{\mathbf{a}} & \textcolor{\complexcolor}{\mathbf{A}} \\ 
\end{bmatrix} \succeq 0 . \label{eq_partition}
 \end{IEEEeqnarray}
 where superscript $\text{H}$ is the conjugate transpose. 
 Note that in general we could position the scalar $\alpha$ anywhere on the diagonal through a permutation of the rows and columns.
 This results in the following matrix variable sizes and structures:
 \begin{IEEEeqnarray}{C?s}
 \textcolor{\complexcolor}{\mathbf{M}} \in \setHermitianMatrix{n+1},  \textcolor{\complexcolor}{\mathbf{A}} =  \mathbf{A}^{\realss}  + \imagnumber  \mathbf{A}^{\imagss}  \in \setHermitianMatrix{n}, \setPSDCone, \nonumber \\
\textcolor{\complexcolor}{\mathbf{a}}= \mathbf{a}^{\realss} + \imagnumber \mathbf{a}^{\imagss}  \in \setComplexMatrix{n}{1}, \textcolor{black}{\alpha} \in \setRealMatrix{1}{1}, \textcolor{black}{\alpha} \geq 0. \nonumber
 \end{IEEEeqnarray}
 We choose a complex matrix parameter $\textcolor{\complexparamcolor}{\mathbf{c}}$ and derive $\textcolor{\complexparamcolor}{\mathbf{C}}$,
 \begin{IEEEeqnarray}{C?s}
\textcolor{\complexparamcolor}{\mathbf{c}} = \textcolor{\paramcolor}{\mathbf{c}}^{\realss} + \imagnumber \textcolor{\paramcolor}{\mathbf{c}}^{\imagss}\in \setComplexMatrix{n}{q} , \quad\nonumber
\textcolor{\complexparamcolor}{\mathbf{C}} = \textcolor{\paramcolor}{\mathbf{C}}^{\realss} + \imagnumber \textcolor{\paramcolor}{\mathbf{C}}^{\imagss} = \textcolor{\complexparamcolor}{\mathbf{c}}\textcolor{\complexparamcolor}{\mathbf{c}}^{\hermitiantranspose}  , \nonumber
 \end{IEEEeqnarray}
  or vice-versa we choose $\textcolor{\complexparamcolor}{\mathbf{C}}  \in \setHermitianMatrix{n},\setPSDCone$ and factor it into $\textcolor{\complexparamcolor}{\mathbf{c}}$.
 We use $\bullet$ to represent the Frobenius inner product of two matrices, 
 \begin{IEEEeqnarray}{C?s}
\textcolor{\complexparamcolor}{\mathbf{C}} \bullet \textcolor{\complexcolor}{\mathbf{A}} = \sum_{ij} ([\textcolor{\complexparamcolor}{\mathbf{C}}]_{ij})^{*} [ \textcolor{\complexcolor}{\mathbf{A}} ]_{ij},
  \end{IEEEeqnarray}
  where $[\mathbf{X}]_{ij}$ means extracting the scalar at index $(i,j)$ in matrix $\mathbf{X}$, and superscript $*$ is the complex conjugate.
For PSD Hermitian matrices this product is nonnegative,
 \begin{IEEEeqnarray}{C?s}
\textcolor{\complexparamcolor}{\mathbf{C}} \bullet \textcolor{\complexcolor}{\mathbf{A}}  \geq 0 \label{eq_frob_def}.
 \end{IEEEeqnarray}
We now illustrate Kim et al.'s main result, a valid convex quadratic inequality for a sampled $\textcolor{\complexparamcolor}{\mathbf{C}}  = \textcolor{\complexparamcolor}{\mathbf{c}}\textcolor{\complexparamcolor}{\mathbf{c}}^{\hermitiantranspose} $, for any Hermitian and PSD matrix $\textcolor{\complexcolor}{\mathbf{M}}$ partitioned with a scalar $\alpha$ on the diagonal:
 \begin{IEEEeqnarray}{C?s}
(\textcolor{\complexparamcolor}{\mathbf{c}}^{\hermitiantranspose} \textcolor{\complexcolor}{\mathbf{a}} )^{\hermitiantranspose} (\textcolor{\complexparamcolor}{\mathbf{c}}^{\hermitiantranspose} \textcolor{\complexcolor}{\mathbf{a}} ) \leq \textcolor{black}{ \alpha } \cdot (\textcolor{\complexparamcolor}{\mathbf{C}} \bullet \textcolor{\complexcolor}{\mathbf{A}} ) . \label{eq_kimkojima_result}
 \end{IEEEeqnarray}
 Implementation can use identities w.r.t. the reals,
 \begin{IEEEeqnarray}{C?s}
 \textcolor{\complexparamcolor}{\mathbf{C}} \bullet \textcolor{\complexcolor}{\mathbf{A}} = 
\mathbf{C}^{\realss} \bullet \mathbf{A}^{\realss}
  + \mathbf{C}^{\imagss} \bullet \mathbf{A}^{\imagss}, \nonumber \\
  \real( \textcolor{\complexparamcolor}{\mathbf{c}}^{\hermitiantranspose} \textcolor{\complexcolor}{\mathbf{a}} ) = (\textcolor{\paramcolor}{\mathbf{c}}^{\realss})^{\transpose} \mathbf{a}^{\realss} + (\textcolor{\paramcolor}{\mathbf{c}}^{\imagss})^{\transpose} \mathbf{a}^{\imagss} , \nonumber \\
  \imag( \textcolor{\complexparamcolor}{\mathbf{c}}^{\hermitiantranspose} \textcolor{\complexcolor}{\mathbf{a}} ) = (\textcolor{\paramcolor}{\mathbf{c}}^{\realss})^{\transpose} \mathbf{a}^{\imagss} - (\textcolor{\paramcolor}{\mathbf{c}}^{\imagss})^{\transpose} \mathbf{a}^{\realss}  \nonumber.
 \end{IEEEeqnarray}
 Note that  we  directly derive (real-value) SOCs for complex-value PSD constraints, \emph{without} having to invoke the isomorphism between the complex numbers and $2\times 2$ matrices (in its  block matrix form) on $\textcolor{\complexcolor}{\mathbf{A}}$ or $\textcolor{\complexcolor}{\mathbf{M}}$.
 
 \section{Application to a $3\times 3$ PSD Voltage Matrix}
 We focus on a 3-bus fully connected power network in this section, and discuss generalization later.
Using the BIM OPF \cite{Low2014} model for a 3-bus  system, with bus voltages $\textcolor{\complexcolor}{U_i} = U^{\realss}_{i} + \imagnumber U^{\imagss}_{i}$, the system `$\textcolor{\complexcolor}{\mathbf{M}}$' matrix in the `W' variable space corresponds to:
  \begin{IEEEeqnarray}{L?s}
       \begin{bmatrix} 
U^{\realss}_{1} + \imagnumber U^{\imagss}_{1}   \\
U^{\realss}_{2} + \imagnumber U^{\imagss}_{2}  \\
 U^{\realss}_{3} + \imagnumber U^{\imagss}_{3}    \\
\end{bmatrix} \begin{bmatrix} 
U^{\realss}_{1} + \imagnumber U^{\imagss}_{1}   \\
U^{\realss}_{2} + \imagnumber U^{\imagss}_{2}  \\
 U^{\realss}_{3} + \imagnumber U^{\imagss}_{3}    \\
\end{bmatrix}^{\hermitiantranspose} \\
=
\begin{bmatrix} 
W^{\realss}_{11}& W^{\realss}_{12} + \imagnumber W^{\imagss}_{12}   & W^{\realss}_{13} + \imagnumber W^{\imagss}_{13}   \\
W^{\realss}_{12} - \imagnumber W^{\imagss}_{12}   & W^{\realss}_{22} &  W^{\realss}_{23} + \imagnumber W^{\imagss}_{23} \\
 W^{\realss}_{13} - \imagnumber W^{\imagss}_{13}  & W^{\realss}_{23} - \imagnumber W^{\imagss}_{23}  & W^{\realss}_{33}   \\
\end{bmatrix} = \textcolor{\complexcolor}{\mathbf{M}} \succeq 0 \nonumber
  \end{IEEEeqnarray}
The partitions \eqref{eq_partition} then are, 
  \begin{IEEEeqnarray}{C?s}
\textcolor{\complexcolor}{\mathbf{A}} \!:\!
 \begin{bmatrix} 
W^{\realss}_{22}\!\! &\! \! W^{\realss}_{23}\! +\! \imagnumber W^{\imagss}_{23}  \\
W^{\realss}_{23}\! -\! \imagnumber W^{\imagss}_{23} \! \!&\!\! W^{\realss}_{33}   \\
\end{bmatrix}\nonumber, 
 \textcolor{\complexcolor}{\mathbf{a}} \!: \!
 \begin{bmatrix}
W^{\realss}_{12} \!-\! \imagnumber W^{\imagss}_{12} \\
 W^{\realss}_{13} \!-\! \imagnumber W^{\imagss}_{13} 
   \end{bmatrix}\nonumber, 
      \alpha \! : \!W^{\realss}_{11}.\nonumber 
     \end{IEEEeqnarray}
\subsection{$2\times 2$ Principal Minor SOC Constraints}

  For this partition, substituting $\textcolor{\complexparamcolor}{\mathbf{c}}=\begin{bmatrix} 1 & 0  \end{bmatrix}^{\text{T}} $ or $\begin{bmatrix} 0 & 1  \end{bmatrix}^{\text{T}}$ into \eqref{eq_kimkojima_result} generates  known  $2\times 2$ PM SOCs, i.e.
     \begin{IEEEeqnarray}{C?s}
      (W^{\realss}_{12})^2 + (W^{\imagss}_{12})^2 \leq W^{\realss}_{11} W^{\realss}_{22}  . \IEEEyesnumber \label{eq_pm_aabb_leq}\\
 (W^{\realss}_{13})^2 + (W^{\imagss}_{13})^2 \leq W^{\realss}_{11} W^{\realss}_{33}  , \IEEEyesnumber\IEEEyessubnumber \label{eq_pm_ccaa_leq}
 \end{IEEEeqnarray}
   A partition with $\alpha = W^{\realss}_{22}$ or $\alpha = W^{\realss}_{33}$ and a one-hot $ \textcolor{\complexparamcolor}{\mathbf{c}}$ generates the third unique PM SOC, i.e.,
        \begin{IEEEeqnarray}{C?s}
   (W^{\realss}_{23})^2 + (W^{\imagss}_{23})^2 \leq W^{\realss}_{22} W^{\realss}_{33}  . \IEEEyessubnumber \label{eq_pm_ccbb_leq}
 \end{IEEEeqnarray}
  
\subsection{Kim et al. SOC Constraints}
We sample from the space of vectors $   \textcolor{\complexparamcolor}{\mathbf{c}} \in \setComplexMatrix{2}{1}$ with the first element normalized, and  parameterize  through $\textcolor{\paramcolor}{r}, \textcolor{\paramcolor}{\theta}$, 
   \begin{IEEEeqnarray}{L}
   \textcolor{\complexparamcolor}{\mathbf{c}} = 
   \begin{bmatrix} 
   1 \\ \textcolor{\paramcolor}{r}e^{\imagnumber \textcolor{\paramcolor}{\theta}}
    \end{bmatrix} = 
     \begin{bmatrix} 
   1 \\ \textcolor{\paramcolor}{r}\cos(\textcolor{\paramcolor}{\theta})
    \end{bmatrix}
    + \imagnumber  
       \begin{bmatrix} 
   0 \\ \textcolor{\paramcolor}{r}\sin( \textcolor{\paramcolor}{\theta})
    \end{bmatrix},\nonumber
     \end{IEEEeqnarray}
Next we obtain a parameterized family of valid SOCs:
   \begin{IEEEeqnarray}{L}
  (W^{\realss}_{12}+ \textcolor{\paramcolor}{r}\cos(\textcolor{\paramcolor}{\theta}) W^{\realss}_{13} -  \textcolor{\paramcolor}{r}\sin(\textcolor{\paramcolor}{\theta}) W^{\imagss}_{13} )^2 \nonumber\\
  + (  W^{\imagss}_{12}  +\textcolor{\paramcolor}{r}\cos(\textcolor{\paramcolor}{\theta})  W^{\imagss}_{13}  + \textcolor{\paramcolor}{r}\sin(\textcolor{\paramcolor}{\theta}) W^{\realss}_{13})^2 \nonumber\\
\!\!\leq W^{\realss}_{11} \left(W^{\realss}_{22}+ \textcolor{\paramcolor}{r}^2 W^{\realss}_{33} +2\textcolor{\paramcolor}{r} \!\left( \cos(\textcolor{\paramcolor}{\theta}) W^{\realss}_{23} -  \sin(\textcolor{\paramcolor}{\theta})  W^{\imagss}_{23} \right) \right)\!.
 \label{eq_kk_aa_relax}
 \end{IEEEeqnarray}
We also generate them for the other partitions, i.e.,
    \begin{IEEEeqnarray}{L}
(W^{\realss}_{12} + \textcolor{\paramcolor}{r}\cos(\textcolor{\paramcolor}{\theta}) W^{\realss}_{23} -  \textcolor{\paramcolor}{r}\sin(\textcolor{\paramcolor}{\theta}) W^{\imagss}_{23} )^2 \nonumber\\
  + (  W^{\imagss}_{12}  -\textcolor{\paramcolor}{r}\cos(\textcolor{\paramcolor}{\theta})  W^{\imagss}_{23}  -  \textcolor{\paramcolor}{r}\sin(\textcolor{\paramcolor}{\theta}) W^{\realss}_{23})^2 \nonumber\\
\!\!\leq W^{\realss}_{22}  
\left(W^{\realss}_{11}+ \textcolor{\paramcolor}{r}^2 W^{\realss}_{33} 
+2 \textcolor{\paramcolor}{r} \!\left( 
\cos(\textcolor{\paramcolor}{\theta}) W^{\realss}_{13} -  \sin(\textcolor{\paramcolor}{\theta}) W^{\imagss}_{13} 
\right) \right)\!. 
\label{eq_kk_bb_relax} \\
  (W^{\realss}_{13} + \textcolor{\paramcolor}{r}\cos(\textcolor{\paramcolor}{\theta}) W^{\realss}_{23} +  \textcolor{\paramcolor}{r}\sin(\textcolor{\paramcolor}{\theta}) W^{\imagss}_{23} )^2 \nonumber\\
  + (  W^{\imagss}_{13}  +\textcolor{\paramcolor}{r}\cos(\textcolor{\paramcolor}{\theta})  W^{\imagss}_{23}  -  \textcolor{\paramcolor}{r}\sin(\textcolor{\paramcolor}{\theta}) W^{\realss}_{23})^2 \nonumber\\
\!\!\!\!\leq W^{\realss}_{33} \! \left(W^{\realss}_{11}+ \textcolor{\paramcolor}{r}^2 W^{\realss}_{22} +2\textcolor{\paramcolor}{r}\! \left( \cos(\textcolor{\paramcolor}{\theta}) W^{\realss}_{12} -  \sin(\textcolor{\paramcolor}{\theta})  W^{\imagss}_{12} \right) \right)\!. \label{eq_kk_cc_relax}
 \end{IEEEeqnarray}
Finally, we still need to state the linear inequalities based on  the Frobenius inner product \eqref{eq_frob_def}, i.e.,
     \begin{IEEEeqnarray}{C?s}
    W^{\realss}_{22}+ \textcolor{\paramcolor}{r^2} W^{\realss}_{33} +2\textcolor{\paramcolor}{r} \left( \cos(\textcolor{\paramcolor}{\theta}) W^{\realss}_{23} -  \sin(\textcolor{\paramcolor}{\theta})  W^{\imagss}_{23} \right) \geq 0 \label{eq_kk_aa_frob}, \IEEEyesnumber\IEEEyessubnumber\\
W^{\realss}_{11}+ \textcolor{\paramcolor}{r^2} W^{\realss}_{33} +2 \textcolor{\paramcolor}{r}\left( \cos(\textcolor{\paramcolor}{\theta}) W^{\realss}_{13} -  \sin(\textcolor{\paramcolor}{\theta})  W^{\imagss}_{13} \right) \geq 0 \label{eq_kk_bb_frob}, \IEEEyessubnumber \\
W^{\realss}_{11}+ \textcolor{\paramcolor}{r^2} W^{\realss}_{22} +2\textcolor{\paramcolor}{r} \left( \cos(\textcolor{\paramcolor}{\theta}) W^{\realss}_{12} -  \sin(\textcolor{\paramcolor}{\theta})  W^{\imagss}_{12} \right) \geq 0 \label{eq_kk_cc_frob}. \IEEEyessubnumber
 \end{IEEEeqnarray}
 Inequalities   \eqref{eq_kk_aa_relax}-\eqref{eq_kk_cc_relax} are families of SOC constraints. 
One re-obtains a  PM SOC for $\textcolor{\paramcolor}{r} = 0$, and similarly if one scales the second element of $\textcolor{\complexparamcolor}{\mathbf{c}}$ instead.

\section{Interpretation with Kirchhoff's Voltage Law}
\subsection{Kirchhoff's Voltage Law for 3-Cycle in Lifted Voltage}
Hijazi et al. \cite{Coffrin2015a} show that Kirchhoff's voltage law (KVL) for a 3-cycle in the lifted voltage variable space is,
   \begin{IEEEeqnarray}{C?s}
(W^{\realss}_{13} + \imagnumber W^{\imagss}_{13})^*(W^{\realss}_{23} + \imagnumber W^{\imagss}_{23})   = {W^{\realss}_{33}}   (W^{\realss}_{12} + \imagnumber W^{\imagss}_{12}  )^*.
 \end{IEEEeqnarray}
 We note that this corresponds to a minor that involves only a single diagonal element, i.e. `type 2' in \cite{Kocuk2017}.
In the reals,  this becomes,
   \begin{IEEEeqnarray}{C?s}
W^{\realss}_{13}W^{\realss}_{23}  + W^{\imagss}_{13} W^{\imagss}_{23}   = W^{\realss}_{33}   W^{\realss}_{12} \IEEEyesnumber\IEEEyessubnumber \label{eq_cycle_c_real},\\
W^{\realss}_{13}W^{\imagss}_{23} - W^{\imagss}_{13} W^{\realss}_{23}   = -W^{\realss}_{33}   W^{\imagss}_{12}  \IEEEyessubnumber \label{eq_cycle_c_imag}.
 \end{IEEEeqnarray}
The real-value permuted 3-cycle KVL expressions are,
   \begin{IEEEeqnarray}{C?s}
W^{\realss}_{12} W^{\realss}_{13}  - W^{\imagss}_{12} W^{\imagss}_{13}   = W^{\realss}_{11}   W^{\imagss}_{23}\IEEEyesnumber\IEEEyessubnumber , \label{eq_cycle_a_real}\\
 W^{\realss}_{12}W^{\imagss}_{13} +  W^{\imagss}_{12} W^{\realss}_{13}   = W^{\realss}_{11}   W^{\imagss}_{23}  \IEEEyessubnumber,\\
W^{\realss}_{12}W^{\realss}_{23}  + W^{\imagss}_{12} W^{\imagss}_{23}   = W^{\realss}_{22}   W^{\realss}_{13} \IEEEyesnumber\IEEEyessubnumber, \\
W^{\realss}_{12}W^{\imagss}_{23}  + W^{\imagss}_{12} W^{\realss}_{23}   = -W^{\realss}_{22}   W^{\imagss}_{13}  \IEEEyessubnumber \label{eq_cycle_b_imag}.
 \end{IEEEeqnarray}
 
 \begin{figure*}[t]
\begin{multicols}{3}
\centering
    \includegraphics[width=0.75\columnwidth]{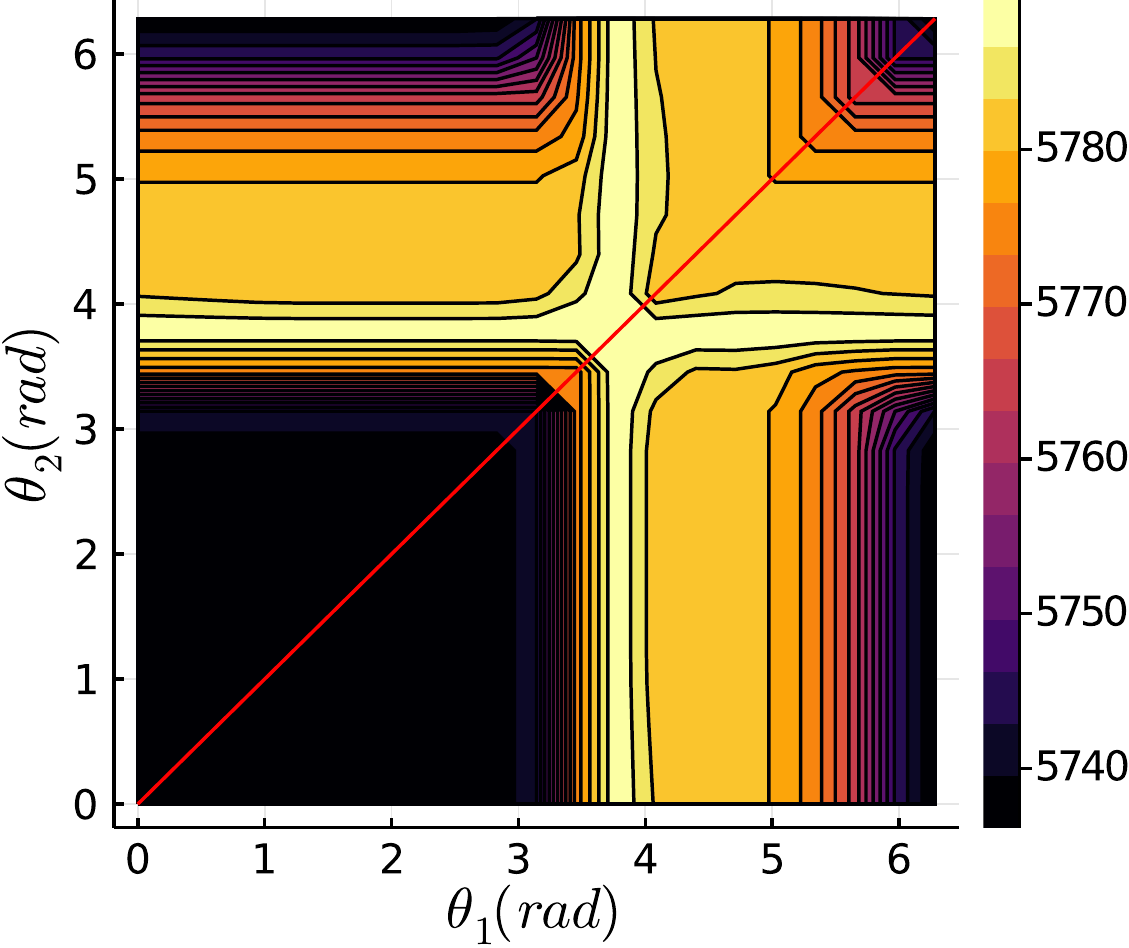}
    \caption{Contour plot for `case3\_lmbd'\label{fig:contour}}
    
    \includegraphics[width=0.75\columnwidth]{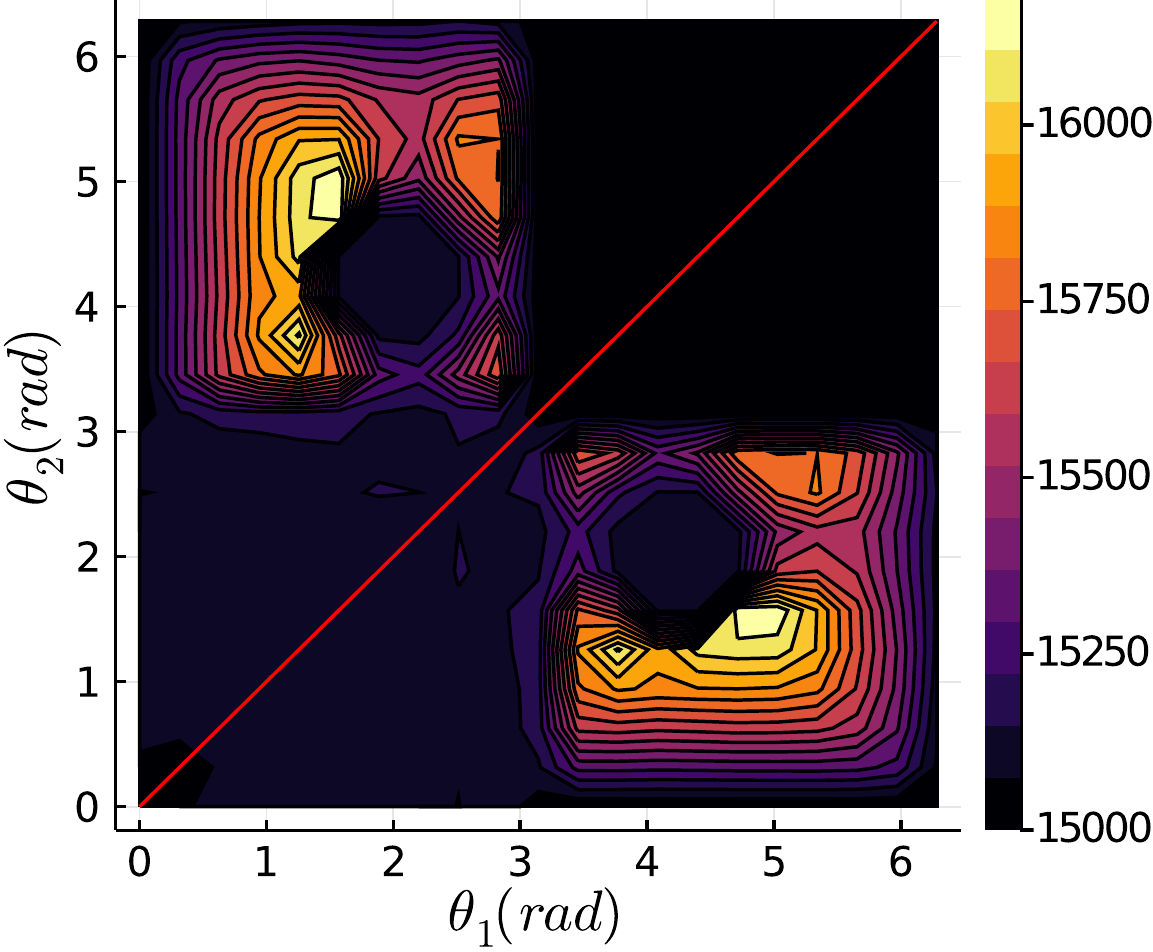} 
    \caption{Contour plot for `case5\_pjm'\label{fig:contour5}}
   
    \includegraphics[width=0.75\columnwidth]{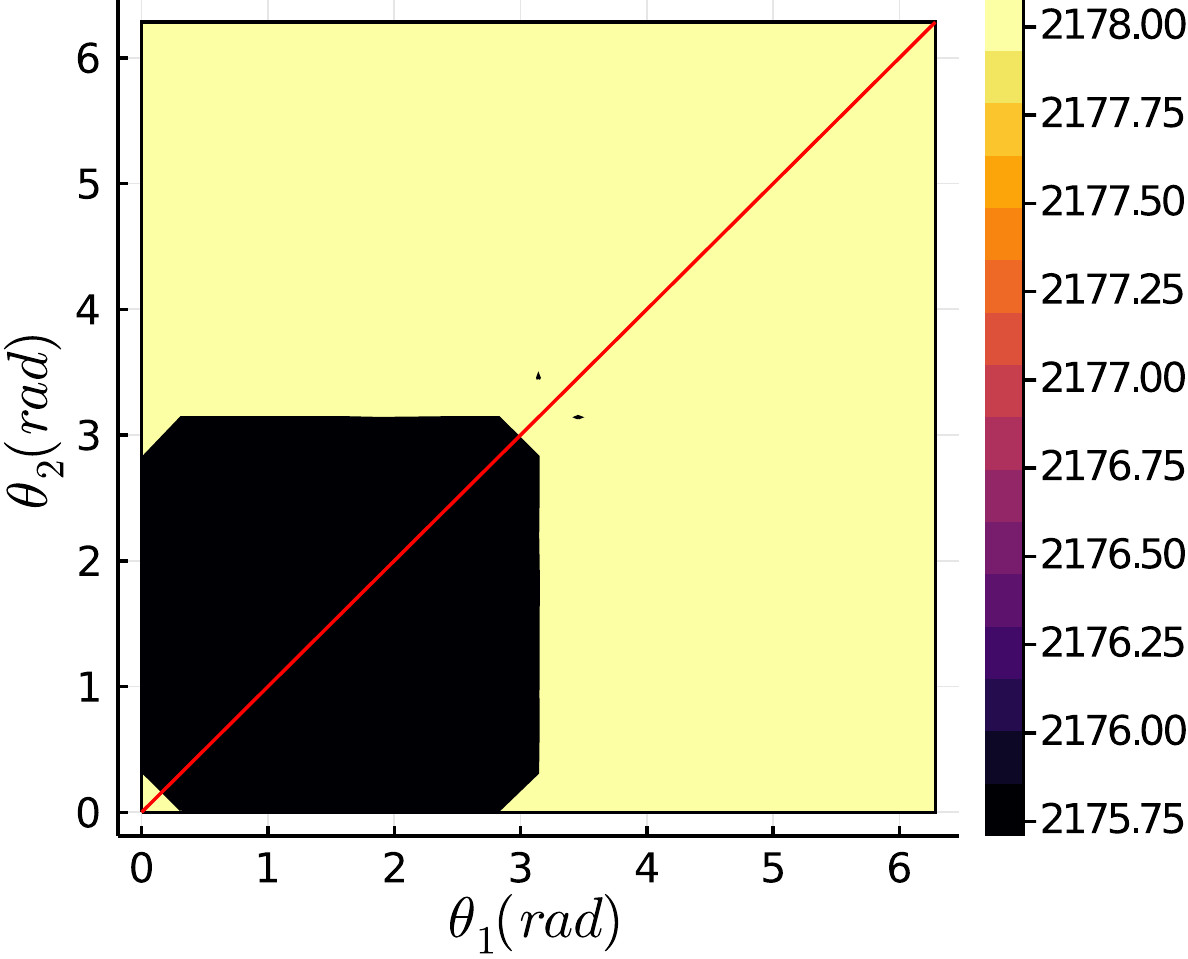} 
    \caption{Contour plot for `case14\_ieee'\label{fig:contour14}}
 \end{multicols}
 \vspace{-2mm}
\end{figure*}

\subsection{Kim et al. SOCs with Solution in Equality}
If $\rank(\textcolor{\complexcolor}{\mathbf{M}}) = 1$,   PM SOCs \eqref{eq_pm_ccaa_leq}-\eqref{eq_pm_ccbb_leq}   are active, 
   \begin{IEEEeqnarray}{C?s}
 (W^{\realss}_{13})^2 + (W^{\imagss}_{13})^2 = W^{\realss}_{33} W^{\realss}_{11},   \IEEEyesnumber\IEEEyessubnumber \label{eq_pm_ccaa}\\
   (W^{\realss}_{23})^2 + (W^{\imagss}_{23})^2 = W^{\realss}_{33}  W^{\realss}_{22}. \IEEEyessubnumber \label{eq_pm_ccbb}
 \end{IEEEeqnarray}
Under the same conditions, substituting  $\textcolor{\paramcolor}{r_1} = \textcolor{\paramcolor}{r_2} = 1,$ and $\textcolor{\paramcolor}{\theta_1} = 0, \textcolor{\paramcolor}{\theta_2} = \pi/2$, in \eqref{eq_kk_aa_relax} we obtain,
   \begin{IEEEeqnarray}{C?s}
    (  W^{\realss}_{13} + W^{\realss}_{23} )^2 + ( W^{\imagss}_{13}  +  W^{\imagss}_{23} )^2 \nonumber \\
    = W^{\realss}_{33} (W^{\realss}_{11} + W^{\realss}_{22} + 2 W^{\realss}_{12}  ), \IEEEyesnumber\IEEEyessubnumber \label{eq_kk_re_cc}\\
     (W^{\realss}_{13} + W^{\imagss}_{23} )^2 + ( W^{\imagss}_{13}  -  W^{\realss}_{23} )^2 \nonumber \\
      = W^{\realss}_{33} (W^{\realss}_{11}+ W^{\realss}_{22} - 2 W^{\imagss}_{12}   ). \IEEEyessubnumber \label{eq_kk_im_cc}
 \end{IEEEeqnarray}
Therefore the following identities hold when $\rank(\textcolor{\complexcolor}{\mathbf{M}}) = 1$,
   \begin{IEEEeqnarray}{C?s}
 \eqref{eq_pm_ccaa} + \eqref{eq_pm_ccbb} +2 \cdot \eqref{eq_cycle_c_real}\equiv  \eqref{eq_kk_re_cc}, \nonumber\\
  \eqref{eq_pm_ccaa} +  \eqref{eq_pm_ccbb} +2\cdot \eqref{eq_cycle_c_imag}\equiv  \eqref{eq_kk_im_cc}, \nonumber
   \end{IEEEeqnarray}
   as well as the equivalents for \eqref{eq_cycle_a_real}-\eqref{eq_cycle_b_imag} (i.e. a total of 6 SOCs). This establishes that the Kim et al. SOCs include a stronger representation of KVL than just the PM SOCs. 
Note that this \emph{completing the square} approach actually embeds nonconvex equalities \eqref{eq_cycle_c_real}-\eqref{eq_cycle_b_imag} into a  convex quadratic form.

Note that by themselves, \eqref{eq_kk_re_cc}-\eqref{eq_kk_im_cc} do not prove \eqref{eq_cycle_c_real}-\eqref{eq_cycle_c_imag} hold;  \eqref{eq_pm_ccaa}-\eqref{eq_pm_ccbb} are needed. 
Therefore observing that \eqref{eq_kk_aa_relax}-\eqref{eq_kk_cc_relax}  are active in the solution is insufficient to prove that $\rank(\textcolor{\complexcolor}{\mathbf{M}}) = 1$.
The PM SOCs \eqref{eq_pm_ccaa}-\eqref{eq_pm_ccbb} must also be active in the solution to give that confidence.  


\section{Numerical Illustration} \label{ch_results}
We extend the BIM implementation in PowerModels \cite{Coffrin2017} and use 3-cycle PGLib test cases \cite{Babaeinejadsarookolaee2019} with two $\textcolor{\paramcolor}{\theta}$-parameterized SOCs per 3-cycle permutation. 
We use Ipopt \cite{Wachter2006} as the solver for NLP `AC' OPF; Hypatia \cite{Coey2020} for the BIM SDP and  the Kim et al. SOC forms, and again Ipopt for the PM SOC. 
Fig.~\ref{fig:contour} - Fig.~\ref{fig:contour14} show sweeps for two SOC families $\textcolor{\paramcolor}{\theta_1}, \textcolor{\paramcolor}{\theta_2}$, w.r.t. objective value, for the 3, 5 and 14 bus test cases. 
The optimal $\textcolor{\paramcolor}{\theta}$ for the 3-bus case is around $3.7$
, with an objective value of $5790$. 
For `case5\_pjm', we want two distinct angles, i.e. $\textcolor{\paramcolor}{\theta_1}\approx 1.6, \textcolor{\paramcolor}{\theta_2}\approx4.9$ with an objective of $16181$;
for `case14\_ieee', one suffices, i.e. any $\textcolor{\paramcolor}{\theta} \in [\pi; 2\pi ]$ gives objective $2178$. It is interesting that, to obtain the tightest relaxation, one case requires two fine-tuned angles, another requires one fine-tuned angle, and the third only requires one angle in a broader range. 
A condition for recognizing that a single angle is sufficient is that the maximum is on the diagonal of the plot. 
Table~\ref{tab_results} summarizes numerical results in comparison with the NLP, SDP and PM SOC relaxations.



\vspace{-1mm}
  \begin{table}[tbh]
  \centering 
  \caption{Objective values; Kim et al. SOCs  for $\textcolor{\paramcolor}{\theta_1}=0, \textcolor{\paramcolor}{\theta_2}=3\pi/2$.}\label{tab_results}
    \begin{tabular}{l  r   r   r    r   }
\hline
		& AC NLP        & BIM SDP             & \textbf{Kim+PM SOC}  & PM SOC   \\ 
		&\$/h       &\% gap           & \% gap &  \% gap  \\ 
		\hline
case3\_lmbd	& 5812.64	 & 0.38 	&0.54  &  1.32  \\ 
case5\_pjm	& 17551.89	 & 5.22 	&14.47  & 14.54  \\ 
case14\_ieee	& 2178.08	 & 0.00 	&0.00  & 0.11 \\ 
\hline
\end{tabular}
\end{table}

We note that in general, for the convergence of the interior point method to be strong, you want to avoid having a large amount of inequalities active in the solution.
In general we expect this to be the case  as many OPF problems are naturally low-rank \cite{Lavaei2012}. 
Therefore, it may be beneficial to use the Kim et al. SOCs algorithmically as nonlinear cuts. 

\section{Conclusions} \label{ch_conclusions}

We develop a SOC relaxation does not entail the `angle relaxation'.
For a $3\times3$ Hermitian matrix, the Kim et al. SOCs sampled in two directions each ($2\times3$ unique ones), when combined with the $2\times 2$ principal minor-derived SOCs (3 unique ones), clearly offers an tighter SOC relaxation. 
Combining with results by Hijazi et al. \cite{Coffrin2015a}, we know now that when all 9 SOCs are satisfied in equality, it furthermore certifies that the original matrix is PSD and rank-1, and that the power flow satisfies Kirchhoff's voltage law in a 3-cycle. 

Nevertheless, these additional 6 SOCs may not be sufficient to obtain a relaxation that is equivalent to the SDP relaxation, as there are actually families of parameteric SOCs from which one needs to sample. The  results do not suggest an obvious sampling strategy.
Profiling on  benchmark libraries \cite{Babaeinejadsarookolaee2019} can aid in strategy development.


Real transmissions networks have loops in the  topology that go beyond 3-cycles. 
Nevertheless, using chordal completion similar to \cite{Kocuk2015}, \cite{jabr2012a}, 
we can re-cast the network graph in terms of intersected 3-cycles, and can apply the technique as such. 
It would be interesting to combine the Kim et al. SOCs with the QC formulation \cite{Duben2015}, to obtain a more scalable  form than QC+SDP, and more accurate than QC+PM SOC. 
In future work, we plan to  explore  generating relaxations to  $4\times4$ and higher-dimensional PSD variable matrices directly. 
A prime use case is  \emph{unbalanced} SDP BIM OPF, as it relies on a large amount of $6\times6$ complex PSD matrices, and solvers struggle with scalability and numerical stability \cite{Gan2014}.



\end{document}